\documentclass[a4paper,11pt]{article}
\usepackage[latin1]{inputenc}
\usepackage[T1]{fontenc}
\usepackage[francais]{babel}
\usepackage{amssymb}
\usepackage{amsfonts, amsmath}
\usepackage{graphicx}
\newtheorem{defn}{D\'efinition}
\newtheorem{Lem}[defn]{Lemma}
\newtheorem{Theo}{Theorem}

\newtheorem{Exmp}{Example}
\newtheorem{Rem}{Remark}

\title{Symplectic manifolds and Hamiltonian dynamical systems}

\author{\textbf{A. Lesfari}
\\\emph{Department of Mathematics}
\\\emph{Faculty of Sciences}
\\\emph{University of Choua\"{i}b Doukkali}
\\\emph{B.P. 20, 24000 El Jadida, Morocco}.
\\\emph{E. mail : lesfari.a@ucd.ac.ma, lesfariahmed@yahoo.fr.}}

\date{}
\begin{document}
\maketitle

\emph{Abstract}. This paper is devoted to the study of symplectic
manifolds and their connection with Hamiltonian dynamical systems.
We review some properties and operations on these manifolds and
see how they intervene when studying the complete integrability of
these systems, with detailed proofs. Several explicit calculations
for which references are not immediately available are given.
These results are exemplified by applications to some
Hamiltonian dynamical systems.

\vskip 0.5cm

\section{Introduction}

It is well known that symplectic manifolds play a crucial role in
classical mechanics, geometrical optics and thermodynamics and
currently has conquered a rich territory, asserting himself as a
central branch of differential geometry and topology. In addition
to its activity as an independent subject, symplectic manifolds
are strongly stimulated by important interactions with many
mathematical and physical specialties among others. The aim of
this paper is to study some properties of symplectic manifolds and
Hamiltonians dynamical systems and to review some operations on
these manifolds, with detailed proofs. This paper is organized as
follows: the first section is an introduction to the subject. In
section 2 we begin by briefly recalling some notions about
symplectic vector spaces. Section 3 develops the explicit
calculation of symplectic structures on a differentiable manifold.
Section 4 is devoted to the study of some properties of
one-parameter groups of diffeomorphisms or flow, Lie derivative,
interior product and Cartan's formula. We review some interesting
properties and operations on differential forms, with detailed
proofs. Section 5 deals with the study of a central theorem of
symplectic geometry namely Darboux's theorem: the symplectic
manifolds $(M, \omega)$ of dimension $2m$ are locally isomorphic
to $(\mathbb{R}^{2m},\omega)$. The classic proof [4] given by
Darboux of his theorem is by recurrence on the dimension of the
variety. We give a preview and see another demonstration [40] due
to Weinstein based on a result of Moser [31]. Section 6 contains
some technical statements concerning Hamiltonian vector fields.
The latter form a Lie subalgebra of the space vector field and we
show that the matrix associated with a Hamiltonian system forms a
symplectic structure. Several properties concerning Hamiltonian
vector fields, their connection with symplectic manifolds, Poisson
manifolds or Hamiltonian manifolds as well as interesting examples
are studied in section 7. We will see in section 8, how to define
a symplectic structure on the orbit of the coadjoint
representation of a Lie group. The remainder is dedicated to the
explicit computation of symplectic structures on adjoint and
coadjoint orbits of a Lie group with particular attention given to
the groups $SO(3)$ and $SO(4)$. Integrable Hamiltonian systems are
nonlinear ordinary differential equations described by a
Hamiltonian function and possessing sufficiently many independent
constants of motion in involution. By the Arnold-Liouville theorem
[2, 4, 9, 24],  the regular compact level manifolds defined by the
intersection of the constants of motion are diffeomorphic to a
real torus on which the motion is quasi-periodic as a consequence
of the following differential geometric fact; a compact and
connected $n$-dimensional manifold on which there exist $n$ vector
fields which commute and are independent at every point is
diffeomorphic to an $n$-dimensional real torus and there is a
transformation to so-called action-angle variables, mapping the
flow into a straight line motion on that torus. Outline in section
9 we give a proof as direct as possible of the Arnold-Liouville
theorem and we make a careful study of its connection with the
concept of completely integrable systems and finally, in section
10, apply it to concrete situations : the problem of the rotation
of a rigid body about a fixed point and Yang-Mills fields valued
in the Lie algebra associated to the Lie group $SU(2)$.

\section{Symplectic Vector Spaces}

First, remember that a symplectic vector space $(E,\omega)$ is a
vector space $E$  over a field equipped with a bilinear form
$\omega:E\times E\longrightarrow\mathbb{R}$ which is alternating
(or antisymmetric, i.e, $\omega(x,y)=-\omega(y,x)$, $\forall
x,y\in E$) and non degenerate (i.e., $\omega(x,y)=0$, $\forall
y\in E\Longrightarrow x=0$). The form $\omega$ is called
symplectic form (or symplectic structure). The dimension of a
symplectic vector space is necessarily even. We show (using a
reasoning similar to the Gram-Schmidt orthogonalization process)
that any symplectic vector space $(E,\omega)$ has a base $(e_1,
..., e_{2m})$ called symplectic basis\index{symplectic basis} (or
canonical basis\index{canonical basis}), satisfying the following
relations : $$\omega(e_{m+i},e_j)=\delta_{ij},
\quad\omega(e_i,e_j)=\omega(e_{m+i},e_{m+j})=0.$$ Note that each
$e_{m+i}$ is orthogonal to all base vectors except $e_i$. In terms
of symplectic basic vectors $(e_1,...,e_{2m})$, the matrix
$(\omega_{ij})$ where $\omega_{ij}\equiv\omega(e_i, e_j)$ has the
form
$$\left(\begin{array}{ccc}
\omega_{11}&...&\omega_{12m}\\
\vdots&\ddots&\vdots\\
\omega_{2m1}&...&\omega_{2m2m}
\end{array}\right)=
\left(\begin{array}{cc}
0&-I_m\\
I_m&0
\end{array}\right),$$
where $I_m$ denotes the $m\times m$ unit matrix.

\begin{Exmp}
The vector space $\mathbb{R}^{2m}$ with the form
$$\omega(x,y)=\sum_{k=1}^m(x_{m+k}y_k-x_ky_{m+k}),\quad x\in\mathbb{R}^{2m},\quad
y\in\mathbb{R}^{2m},$$ is a symplectic vector space. Let
$(e_1,...,e_m)$ be an orthonormal basis of $\mathbb{R}^m$. Then,
$((e_1,0),...,(e_m, 0), (0, e_1),...,(0, e_m))$ is a symplectic
base of $\mathbb{R}^{2m}$.
\end{Exmp}

Let $(E, \omega)$ be a symplectic vector space and $F$ a vector
subspace of $E$. Let $F^\bot$ be the the orthogonal (symplectic)
of $F$, i.e., the vector subspace of $E$ defined by
$$F^\bot=\{x\in E:\forall y\in F,\omega(x,y)=0\}.$$ The subspace
$F$ is isotropic if $F\subset F^\bot$, coisotropic if
$F^\bot\subset F$, Lagrangian if $F=F^\bot$ and symplectic if
$F\cap F^\bot=\{0\}$. If $F$, $F_1$ and $F_2$ are subspaces of a
symplectic space $(E, \omega)$, then $$\dim F+\dim F^\bot=\dim E,
\qquad(F^\bot)^\bot=F,$$ $$F_1\subset F_2\Longrightarrow
F_2^\bot\subset F_1^\bot,$$ $$(F_1\cap
F_2)^\bot=F_1^\bot+F_2^\bot, \qquad F_1^\bot\cap
F_2^\bot=(F_1+F_2)^\bot,$$ $F$ coisotropic if and only if $F^\bot$
isotropic and $F$ Lagrangian if and only if$ F$ isotropic and
coisotropic.

\section{Symplectic Manifolds}

We will define a symplectic structure on a differentiable manifold
and study some properties. A symplectic structure (or symplectic
form) on an even-dimensional differentiable manifold $M$ is a
closed non-degenerate differential $2$-form $\omega$ defined
everywhere on $M$. The non-degeneracy condition means that :
$$\forall x\in M,\quad\forall \xi \neq 0,\quad\exists \eta : \omega \left(\xi
,\eta \right) \neq 0,\left(\xi,\eta \in T_{x}M\right).$$ The pair
$(M,\omega)$ (or simply $M$) is called a symplectic manifold.
Hence, at a point $p\in M$, we have a non-degenerate antisymmetric
bilinear form on the tangent space $T_pM$, which explains why the
dimension of the $M$ manifold is even.

\begin{Exmp}
The space $M=\mathbb{R}^{2m}$ with the $2$-form
$$\omega=\sum_{k=1}^mdx_k\wedge dy_k,$$ is a symplectic manifold.
The vectors $$\left(\frac{\partial}{\partial
x_1}\right)_p,...,\left(\frac{\partial}{\partial x_m}\right)_p,
\left(\frac{\partial}{\partial
y_1}\right)_p,...,\left(\frac{\partial}{\partial y_m}\right)_p,
\quad p\in M,$$ constitute a symplectic basis of the tangent space
$T_pM$. Similarly, space $\mathbb{C}^{m}$ with the form
$$\omega=\frac{i}{2}\sum_{k=1}^mdz_k\wedge d\overline{z}_k,$$ is a
symplectic manifold. Note that this form coincides with that of
the preceding example by means of the identification
$\mathbb{C}^m\simeq\mathbb{R}^{2m}$, $z_k=x_k+iy_k$. The Riemann
surfaces are symplectic manifolds. Other examples of symplectic
manifolds which will not be considered here (and for which I refer
for example to [2, 26]) are the k\"{a}hlerian manifolds as well as
complex projective manifolds. Another important class of
symplectic manifolds consists of the coadjoint orbits
$\mathcal{O}\subset \mathcal{G}^*,$ where $\mathcal{G}$ is the
algebra of a Lie group $\mathcal{G}$ and
$$\mathcal{G}_\mu=\{Ad^*_g \mu : g\in \mathcal{G}\},$$ is the orbit
of $\mu \in \mathcal{G}^*$ under the coadjoint representation (to
see further).
\end{Exmp}

We will see that the cotangent bundle $T^{*}M$ (that is, the union
of all cotangent spaces at $M$) admits a natural symplectic
structure. The phase spaces of the Hamiltonian systems studied
below are symplectic manifolds and often they are cotangent
bundles equipped with the canonical structure.

\begin{Theo}
Let $M$ be a differentiable manifold of dimension $m$ and let
$T^{*}M$ be its cotangent bundle. Then $T^{*}M$ possesses in a
natural way a symplectic structure and in a local coordinate
$\left(x_{1},\ldots,x_{m},y_{1},\ldots,y_{m}\right)$, the form
$\omega$ is given by $$\omega=\sum_{k=1}^{m}dx_{k}\wedge dy_{k}.$$
\end{Theo}
\emph{Proof}. Let $(U, \varphi)$  be a local chart in the
neighborhood of $p\in M$,
$$\varphi : U\subset M\longrightarrow\mathbb{R}^m,\quad
p\longmapsto\varphi(p)=\sum_{k=1}^mx_ke_k,$$ where $e_k$ are the
vectors basis of $\mathbb{R}^m$. Consider the canonical
projections $TM\longrightarrow M$, and $T(T^*M)\longrightarrow
T^*M$, of tangent bundles respectively to $M$ and $T^*M$ on their
bases. We notice $$\pi^* : T^*M\longrightarrow M,$$ the canonical
projection and $$d\pi^* : T(T^*M)\longrightarrow TM,$$ its linear
tangent application. We have
$$\varphi^* : T^*M\longrightarrow \mathbb{R}^{2m},\quad
\alpha\longmapsto\varphi^*(\alpha)=\sum_{k=1}^m(x_ke_k+y_k\varepsilon_k),$$
where $\varepsilon_k$ are the basic forms of $T^*\mathbb{R}^m$ and
$\alpha$ denotes $\alpha_p\in T^*M$. So, if $\alpha$ is a $1$-form
on  $M$ and $\xi_\alpha$ is a vector tangent to $T^*M$, then
$$d\varphi^* : T(T^*M)\longrightarrow T\mathbb{R}^{2m}=\mathbb{R}^{2m},\quad
\xi_\alpha\longmapsto
d\varphi^*(\xi_\alpha)=\sum_{k=1}^m(\beta_ke_k+\gamma_k\varepsilon_k),$$
where $\beta_k$, $\gamma_k$ are the components of $\xi_\alpha$ in
the local chart of $\mathbb{R}^{2m}$. Let
$$\lambda_\alpha(\xi_\alpha)=\alpha(d\pi^*\xi_\alpha)=\alpha(\xi),$$
where $\xi$ is a tangent vector to $M$. Let $(x_1,...,x_m,
y_1,...,y_m)$ be a system of local coordinates compatible with a
local trivialization of the tangent bundle $T^*M$. Let's show
that:
$$
\lambda_\alpha(\xi_\alpha)=\alpha\left(\sum_{k=1}^m\beta_ke_k\right)=
\sum_{k=1}^m(x_ke_k+y_k\varepsilon_k)\left(\sum_{j=1}^m\beta_je_j\right)=
\sum_{k=1}^m\beta_ky_k.$$ Indeed, remember that if $(x_1,...,x_m)$
is a system of local coordinates around $p\in M$, like all
$\alpha\in T^*M$ can be written in the basis $(dx_1,...,dx_m)$
under the form $$\alpha=\sum_{k=1}^m\alpha_kdx_k,$$ then by
defining local coordinates $y_1,...,y_m$ by $y_k(\alpha)=y_k$,
$k=1,...,m$, the $1$-form $\lambda$ is written
$$\lambda=\sum_{k=1}^my_kdx_k.$$ The form $\lambda$ on the cotangent
bundle $T^*M$ doing correspondence $\lambda_\alpha$ to $\alpha$ is
called Liouville form. We have
\begin{eqnarray}
\lambda(\alpha)&=&\sum_{k=1}^my_k(\alpha)dx_k(\alpha),\nonumber\\
\lambda(\alpha)(\xi_\alpha)&=&\sum_{k=1}^my_k(\alpha)dx_k(\alpha)
\left(\sum_{j=1}^m\beta_je_j+\gamma_j\varepsilon_j\right),\nonumber
\end{eqnarray}
i.e.,
$$
\lambda(\alpha)(\xi_\alpha)=\sum_{k=1}^my_k\beta_k=\lambda_\alpha(\xi_\alpha),\qquad
\lambda=\sum_{k=1}^my_kdx_k.$$ The symplectic structure of $T^*M$
is given by the exterior derivative of $\lambda$, i.e., the
$2$-form $\omega=-d\lambda$. The forms $\lambda$ and $\omega$ are
called canonical forms on $T^*M$. We can visualize all this with
the help of the following diagram :
$$
\begin{array}{ccccccc}
&&&&T^*(T^*M)&&\\
&&&&\left\uparrow \lambda \right.&&\\
\mathbb{R}&\overset{\lambda_\alpha(\xi)}{\longleftarrow
}&T(T^*M)&\longrightarrow&T^*M&
\overset{\varphi^* }{\longrightarrow }&\mathbb{R}^{2m}\\
&&\left\downarrow d\pi ^{*}\right.&&\left\downarrow \pi^*\right.&&\\
\mathbb{R}&\overset{\alpha \left( \xi \right) }{\longleftarrow
}&TM&\longrightarrow&M& \overset{\varphi}{\longrightarrow
}&\mathbb{R}^m
\end{array}
$$
The form $\omega$ is closed : $d\omega=0$ since $d\circ d=0$ and
it is non degenerate. To show this last property, just note that
the form is well defined independently of the chosen coordinates
but we can also show it using a direct calculation. Indeed, let
$\xi=(\xi_1,...,\xi_{2m})\in T_pM$ and
$\eta=(\eta_1,...,\eta_{2m})\in T_pM$. We have
$$
\omega(\xi,\eta)=\sum_{k=1}^{m}dx_{k}\wedge
dy_{k}(\xi,\eta)=\sum_{k=1}^{m}\left(dx_{k}(\xi)
dy_{k}(\eta)-dx_{k}(\eta)dy_{k}(\xi)\right).$$ Since
$dx_{k}(\xi)=\xi_{m+k}$ is the $(m+k)^{th}$-component of $\xi$ and
$dy_{k}(\xi)=\xi_{k}$ is the $k^{th}$-component of $\xi$, then
$$
\omega(\xi,\eta)=\sum_{k=1}^{m}(\xi_{m+k} \eta_k-\eta_{m+k}\xi_k)
=(\xi_1...\xi_{2m})\left(\begin{array}{cc}
O&-I\\
I&O
\end{array}\right)\left(\begin{array}{c}
\eta_1\\
\vdots\\
\eta_{2m}
\end{array}\right),
$$
with $O$ the null matrix and $I$ the unit matrix of order $m$.
Then, for all $x\in M$ and for all $\xi=(\xi_1,...,\xi_{2m})\neq
0$, it exists $\eta =(\xi_{m+1},...,\xi_{2m},-\xi_1,...,-\xi_m)$
such that : $$\omega(\xi,\eta)=\sum_{k=1}^{m}\left(\xi^2_{m+k}-
\xi^2_k\right)\neq 0,$$ because $\xi_k\neq 0$, for any integer
$k=1,...,2m$. In the local coordinate system
$(x_1,...,x_m,y_1,...,y_m)$, this symplectic form is written
$$\omega=\sum_{k=1}^ndx_k\wedge dy_k,$$ which completes the proof.
$\square$

A manifold $M$, is said to be orientable if there exists on $M$ an
atlas such that the Jacobian of any change of chart is strictly
positive or if $M$ has a volume form (i.e., a differential form
that does not vanish anywhere). For example, $\mathbb{R}^n$ is
oriented by the volume form $dx_1\wedge...\wedge dx_n$. The circle
$S^1$ is oriented by $d\theta$. The torus $T^2=S^1\times S^1$ is
oriented by the volume form $d\theta\wedge d\varphi$. All
holomorphic manifolds are orientable.

\begin{Theo}
a) A closed differential $2$-form $\omega$ on a differentiable
manifold $M$ of dimension $2m$ is symplectic, if and only if,
$\omega^m$ is a volume form.

b) Any symplectic manifold is orientable.

c) Any orientable manifold of dimension two is symplectic. On the
other hand in even dimensions larger than $2$, this is no longer
true.
\end{Theo}
\emph{Proof}. a) Indeed, this is due to the fact that the
non-degeneracy of $\omega$ is equivalent to the fact that
$\omega^m$ is never zero.

b) In a system of symplectic charts $(x_1,...,x_{2m})$, we have
$$\omega=dx_1\wedge dx_{m+1}+\cdots+dx_m\wedge dx_{2m}.$$
Therefore,
\begin{eqnarray}
\omega^m&=&dx_1\wedge dx_{m+1}\wedge...\wedge dx_m\wedge
dx_{2m},\nonumber\\
&=&(-1)^{\frac{m(m-1)}{2}}dx_1\wedge dx_2\wedge...\wedge
dx_{2m},\nonumber
\end{eqnarray}
which means that the $2m$-form $\omega^m$ is a volume form on the
manifold $M$ and therefore this one is orientable. The orientation
associated with the differential form $\omega$ is the canonical
orientation of $\mathbb{R}^{2m}$.

c) This results from the fact that any  differential $2$-form on a
$2$-manifold is always closed. $\square$

\begin{Theo}
Let $\alpha$ be a differential $1$-form on the manifold $M$ and
denote by $\alpha^*\lambda$ the reciprocal image of the Liouville
form $\lambda$ on the cotangent bundle $T^*M$. Then, we have
$\alpha^*\lambda=\alpha$.
\end{Theo}
\emph{Proof}. Since $\alpha : M\longrightarrow T^*M$, we can
consider the reciprocal image that we note $\alpha^* :
T^*T^*M\longrightarrow T^*M$, of $\lambda : T^*M\longrightarrow
T^*T^*M$ (Liouville form), such that, for any vector $\xi$ tangent
to $M$, we have the following relation
$$\alpha^*\lambda(\xi)=\lambda(\alpha)(d\alpha\xi).$$ Since
$d\alpha$ is an application $TM\longrightarrow TT^*M$, then
\begin{eqnarray}
\alpha^*\lambda(\xi)&=&\lambda(\alpha)(d\alpha\xi),\nonumber\\
&=&\lambda_\alpha(d\alpha\xi),\nonumber\\
&=&\alpha d\pi^*d\alpha(\xi),\nonumber\\
&=&\alpha d(\pi^*\alpha)(\xi),\nonumber\\
&=&\alpha(\xi),\nonumber
\end{eqnarray}
because $\pi^*\alpha(p)=p$ where $p\in M$ and the result follows.
$\square$

A submanifold $\mathcal{N}$ of a symplectic manifold $M$ is called
Lagrangian if for all $p\in\mathcal{N}$, the tangent space
$T_p\mathcal{N}$ coincides with the following configuration space
: $\{\eta\in T_pM : \omega_p(\xi,\eta)=0, \forall \xi\in
T_p\mathcal{N}\}$. On this space the $2$-form $\sum dx_k\wedge
dy_k$  that defines the symplectic structure is identically zero.
Lagrangian submanifolds are considered among the most important
submanifolds of symplectic manifolds. Note that
$\dim\mathcal{N}=\frac{1}{2}\dim M$  and that for all vector
fields $X$, $Y$ on $\mathcal{N}$, we have $\omega(X,Y)=0$.

\begin{Exmp}
If $(x_1, ..., x_m, y_1, ..., y_m)$ is a local coordinate system
on an open $U\subset M$, then the subset of $U$ defined by
$y_1=\cdots=y_m=0$ is a Lagrangian submanifold of  $M$. The
submanifold $\alpha(M)$ is Lagrangian in $T^*M$ if and only if the
form $\alpha$ is closed because
$$0=\alpha^*\omega=\alpha^*(-d\lambda)=-d(\alpha^*\lambda)=-d\alpha.$$
\end{Exmp}

Let $M$ be a differentiable manifold, $T^*M$ its cotangent bundle
with the symplectic form $\omega$, and $$s_\alpha: U
\longrightarrow T^*M, \quad p\longmapsto\alpha(p),$$ a section on
an open $U\subset M$. From the local expression of $\omega$
(theorem 1), we deduce that the null section of the bundle $T^*M$
is a Lagrangian submanifold of $T^*M$. If $s_\alpha(U)$ is a
Lagrangian submanifold of $T^*M$, then $s_\alpha$  is called
Lagrangian section. We have (theorem 3),
$s^*_\alpha\lambda=\alpha$, and according to example 3,
$s_\alpha(U)$ is a Lagrangian submanifold of $T^*M$ if and only if
the form $\alpha$ is closed. Let $(M, \omega)$ and $(N, \eta)$ be
two symplectic manifolds of the same dimension and $f:
M\longrightarrow N$, a differentiable application. We say that $f$
is a symplectic morphism if it preserves the symplectic forms,
i.e., $f$ satisfies $f^*\eta=\omega$. When $f$ is a
diffeomorphism, we say that $f$ is a symplectic diffeomorphism or
$f$ is a symplectomorphism.

\begin{Theo}
a) A symplectic morphism is a local diffeomorphism.

b) A symplectomorphism  preserve the orientation.

\end{Theo}
\emph{Proof}. a) Indeed, the $2$-form $\Omega$ being non
degenerate then the differential $$df(p) : T_pM\longrightarrow
T_pM, \quad p\in M,$$ is a linear isomorphism and according to the
local inversion theorem, $f$ is a local diffeomorphism. Another
proof is to note that $$f^*\eta^m=(f^*\eta)^m=\omega^m.$$ The map
$f$ has constant rank $2m$ because $\omega^m$ and $\eta^m$ are
volume forms on $M$ et $N$ respectively. And the result follows.

b) It is deduced from a) that the symplectic diffeomorphisms or
symplectomorphisms  preserve the volume form and therefore the
orientation. The Jacobian determinant of the transformation is
$+1$. $\square$

\begin{Rem}
Note that the inverse $f^{-1}:N\longrightarrow M$ of a
symplectomorphism $f:M\longrightarrow N$ is also a
symplectomorphism.
\end{Rem}

Let $(M, \omega)$, $(N, \eta)$ be two symplectic manifolds,
$$pr_1:M\times N\longrightarrow M, \qquad pr_2:M\times N\longrightarrow
N,$$ the projections of $M\times N$ on its two factors. The two
forms $pr_1^*\omega+pr_2^*\eta$ and $pr_1^*\omega-pr_2^*\eta$ on
the product manifold $M\times N$, are symplectic forms. Take the
case where $\dim M=\dim N=2m$ and consider a differentiable map
$f: M\longrightarrow N$, as well as its graph defined by the set
$$A=\{(x,y)\in M\times N: y=f(x)\}.$$ Note that the application $g$
defined by $$g:M\longrightarrow A, \quad x\longmapsto (x,f(x)),$$
is a diffeomorphism. We show that $A$ is a $2m$-dimensional
Lagrangian submanifold of $(M\times N, pr_1^*\omega-pr_2^*\eta)$
if and only if the reciprocal image of $pr_1^*\omega-pr_2^*\eta$
by the application $g$ is the identically zero form on $M$.
Therefore, for the differentiable map $f$ to be a symplectic
morphism, it is necessary and sufficient that the graph of $f$ is
a Lagrangian submanifold of the product manifold $(M\times N,
pr_1^*\omega-pr_2^*\eta)$.

\begin{Theo}
a) Let $f:M\longrightarrow M$ be a diffeomorphism. Then, the
application $f^*: T^*M\longrightarrow T^*M$, is a
symplectomorphism.

b) Let $g: T^*M\longrightarrow T^*M$ be a diffeomorphism such that
: $g^*\lambda=\lambda$. Then, there is a diffeomorphism
$f:M\longrightarrow M$ such that : $g=f^*$.
\end{Theo}
\emph{Proof}. a) Let's show that $f^{**}\omega=\omega$. We have
\begin{eqnarray}
f^{**}\lambda(\alpha)(\xi_\alpha)&=&\lambda(f^*(\alpha))(df^*\xi_\alpha),\nonumber\\
&=&f^*(\alpha)d\pi^*df^*(\xi_\alpha),\nonumber\\
&=&\alpha (dfd\pi^*df^*(\xi_\alpha)),\nonumber
\end{eqnarray}
and therefore,
$$f^{**}\lambda(\alpha)(\xi_\alpha)=\alpha (d(f\circ\pi^*\circ
f^*)(\xi_\alpha)).$$ Since $f^*\alpha=\alpha_{f^{-1}(p)}$ and
 $\pi^*f^*\alpha=f^{-1}(p)$, then $$f\circ\pi^*\circ
f^*(\alpha)=p=\pi^*\alpha,$$ i.e.,
\begin{equation}\label{eqn:euler}
f\circ\pi^*\circ f^*=\pi^*
\end{equation}
and
$$
f^{**}\lambda(\alpha)(\xi_\alpha)=\alpha(d\pi^*(\xi_\alpha))=
\lambda_\alpha(\xi_\alpha)=\lambda(\alpha)(\xi_\alpha).$$
Consequently, $f^{**}\lambda=\lambda$, and $f^{**}\omega=\omega$.

b) Since $g^*\lambda=\lambda$, then
$$
g^*\lambda(\eta)=\lambda(dg\eta)=\omega(\xi,
dg\eta)=\lambda(\eta)=\omega(\xi, \eta).$$ Moreover, we have
$g^*\omega=\omega$, hence
$$\omega(dg\xi, dg\eta)=\omega(\xi, \eta)=\omega(\xi, dg\eta),$$$$
\omega(dg\xi-\xi, dg\eta)=0,\quad \forall\eta.$$ Since the form
$\omega$ is non-degenerate, we deduce that $dg\xi=\xi$ and that
$g$ preserves the integral curves of $\xi$. On the null section of
the tangent bundle (i.e., on the manifold), we have $\xi=0$ and
then $g|_M$ is an application $f: M\longrightarrow M$. Let's show
that : $$f\circ\pi^*\circ g=\pi^*=f\circ\pi^*\circ f^*.$$ Indeed,
taking the differential, we get $$df\circ d\pi^*\circ
dg(\xi)=df\circ d\pi^*(\xi)=df(\xi_p),$$ because $dg(\xi)=\xi$ and
$\xi_p\equiv d\pi^*(\xi))$, hence, $$df\circ d\pi^*\circ
dg(\xi)=\xi_p=d\pi^*(\xi).$$ Therefore, $$df\circ d\pi^*\circ
dg=d\pi^*,$$ and $$f\circ\pi^*\circ g=\pi^*.$$ Since
$f\circ\pi^*\circ f^*=\pi^*$ (according to (1)), so $g=f^*$.
$\square$

\begin{Theo}
Let $$I:T_{x}^{*}M\longrightarrow T_{x}M,
\quad\omega_{\xi}^{1}\longmapsto \xi,$$ be a map defined by the
relation $$\omega_{\xi}^{1}\left(\eta \right)=\omega \left(\eta
,\xi\right), \quad\forall \eta \in T_{x}M.$$ Then $I$ is an
isomorphism generated by the symplectic form $\omega$.
\end{Theo}
\emph{Proof}. Denote by $I^{-1}$ the map
$$I^{-1}:T_{x}M\longrightarrow T_{x}^{*}M, \quad\xi \longmapsto
I^{-1}(\xi)\equiv\omega_{\xi}^{1},$$ with
$$I^{-1}(\xi)(\eta)=\omega_{\xi}^{1}(\eta)=\omega(\eta,\xi),
\quad\forall \eta \in T_{x}M.$$ The fact that the form $\omega$ is
bilinear implies that
$$I^{-1}(\xi_{1}+\xi_{2})(\eta)=I^{-1}(\xi_{1})(\eta)+I^{-1}(\xi_{2})(\eta),
\quad\forall\eta\in T_{x}M.$$ Since $\dim T_{x}M=\dim T_{x}^{*}M,$
to show that $I^{-1}$ is bijective, it suffices to show that is
injective. The form $\omega$ is non-degenerate, it follows that
$KerI^{-1}=\{0\}$. Hence $I^{-1}$ is an isomorphism and
consequently $I$ is also an isomorphism (the inverse of an
isomorphism is an isomorphism). $\square$

\section{Flows, Lie derivative, inner product and Cartan's formula}

Let $M$ be a differentiable manifold of dimension $m$. Let $TM$ be
the tangent bundle to $M$, i.e., the union of spaces tangent to
$M$ at all points $x$, $TM=\bigcup_{x\in M}T_{x}M$. This bundle
has a natural structure of differentiable variety of dimension
$2m$ and it allows us to convey immutably to the manifolds the
whole theory of ordinary differential equations. Let
$X:M\longrightarrow TM$, be a vector field assumed to be different
from the zero vector of $TM_x$ only on a compact subset $K$ of the
manifold $M$.

Given a point $x\in M$, we write $g_{t}^{X}(x)$ (or simply
$g_{t}(x)$) the position of $x$ after a displacement of a duration
$t\in \mathbb{R}$. There is thus an application
$$g_{t}^{X} : M\longrightarrow M, \quad t\in \mathbb{R},$$
which is a diffeomorphism (a one-to-one differentiable mapping
with a differentiable inverse), by virtue of the theory of
differential equations (see theorem below). The vector field $X$
generates a one-parameter group of diffeomorphisms $g_{t}^{X}$ on
$M$, i.e., a differentiable application ($\mathcal{C}^{\infty }$)
: $M\times \mathbb{R}\longrightarrow M$, satisfying a group law :

$i)$ $\forall t\in \mathbb{R},\text{ }g_{t}^{X} : M\longrightarrow
M$ is a diffeomorphism.

$ii)$ $\forall t,s\in \mathbb{R},\text{
}g_{t+s}^{X}=g_{t}^{X}\circ g_{s}^{X}$.

The condition $ii)$ means that the mapping $t\longmapsto
g_{t}^{X}$, is a homomorphism of the additive group $\mathbb{R}$
into the group of diffeomorphisms of $M$ in $M$. It implies that
$g_{-t}^{X}=\left(g_{t}^{X}\right)^{-1}$, because
$g_{0}^{X}=id_{M}$ is the identical transformation that leaves
every point invariant.

The one-parameter group of diffeomorphisms $g_{t}^{X}$ on $M$,
which we have just described is called a flow and it admits the
vector field $X$  for velocity fields
$$\frac{d}{dt}g_{t}^{X}(x)=X\left(g_{t}^{X}(x)\right),$$ with the initial condition : $g_{0}^{X}(x)=x$.
Obviously
$$\left.\frac{d}{dt}g_{t}^{X}(x)\right|_{t=0}=X(x).$$
Hence by these formulas $g_{t}^{X}(x)$ is the curve on the
manifold that passes through $x$ and such that the tangent at each
point is the vector $X\left(g_{t}^{X}(x)\right)$.

We will show how to construct the $g_{t}^{X}$ flow on the manifold
$M$.

\begin{Theo}
The vector field $X$ generates a unique group of diffeomorphisms
of the compact manifold $M$. In addition, every solution of the
differential equation
$$\frac{dx(t)}{dt}=X(x(t)), \quad x\in M$$ with the initial
condition $x$ (for $t=0$), can be extended indefinitely. The value
of the solution $g_{t}^{X}(x)$ at time $t$ is differentiable with
respect to $t$ and the initial condition $x$.

\end{Theo}
\emph{Proof}. For the construction of $g_{t}^{X}$ for $t$ small,
we proceed as follows : for $x$ fixed, the differential equation
$$\frac{d}{dt}g_{t}^{X}(x)=X\left(g_{t}^{X}\right),$$
function of $t$ with the initial condition : $g_{0}^{X}(x)=x$,
admits a unique solution $g_{t}^{X}$ defined in the neighborhood
of the point $x_{0}$ smoothly ($\mathcal{C}^{\infty }$) depending
on the initial condition. Then $g_{t}^{X}$ is locally a
diffeomorphism. For each point $x_{0}\in M$, we can find a
neighborhood $U\left(x_{0}\right) \subset M$, a real positive
number $\varepsilon\equiv\varepsilon\left(x_{0}\right)$ such that
for all $t\in\left]-\varepsilon ,\varepsilon \right[$, the
differential equation in question with its initial condition has a
unique differentiable solution $g_{t}^{X}(x)$ defined in
$U\left(x_{0}\right)$ and satisfying the group relation
$$g_{t+s}^{X}(x)=g_{t}^{X}\circ g_{s}^{X}(x),$$
with $t,s,t+s\in\left]-\varepsilon ,\varepsilon\right[$. Indeed,
put $x_{1}=g_{t}^{X}(x)$, $t$ fixed and consider the solution of
the differential equation satisfying in the neighborhood of the
point $x_ {0}$ to the initial condition : $g_{s=0}^{X}=x_{1}$.
This solution satisfies the same differential equation and
coincides in a point $g_{t}^{X}(x)=x_{1}$, with the function
$g_{t+s}^{X}$. Therefore, by uniqueness of the solution of the
differential equation, the two functions are locally equal.
Therefore, the application $g_{t}^{X}$ is locally a
diffeomorphism. The vector field $X$ is assumed to be
differentiable (of class $\mathcal{C}^{\infty }$) and with compact
support $K$. Since $K$ is compact, then from the open covering
$U(x)$ of $K$, we can extract a finite sub-covering
$\left(U_{i}\right)$. Let us denote by $\varepsilon_{i}$ the
numbers $\varepsilon$ corresponding to $U_{i}$ and let
$$\varepsilon _{0}=\inf \left(\varepsilon _{i}\right), \quad
g_{t}^{X}(x)=x,\quad x\notin K.$$ The equation in question admits
a unique solution $g_{t}^{X}$ on
$M\times\left]-\varepsilon_{0},\varepsilon_{0}\right[$ satisfying
the relation of the group above, the inverse of $g_{t}^{X}$ being
$g_{-t}^{X}$ and therefore $g_{t}^{X}$ is a diffeomorphism for $ t
$ sufficiently small. We will now see how to construct $g_{t}^{X}$
for every $t$ $\in \mathbb{R}$. From what precedes, just construct
$g_{t}^{X}$ for
$t\in\left]-\infty,-\varepsilon_{0}\right[\cup\left]\varepsilon_{0},\infty\right[$.
We will see that the applications $g_{t}^{X}$ are defined
according to the multiplication law of the group. Note that $t$
can be written as
$$t=k\frac{\varepsilon_{0}}{2}+r, \quad k\in
\mathbb{Z},\quad r\in\left[0,\frac{\varepsilon_{0}}{2}\right[.$$
Let
$$g_{t}^{X}=\underset{k-\text{times}}{\underbrace{g_{\frac{\varepsilon _{0}}{2}}^{X}\circ
\cdots \circ g_{\frac{\varepsilon _{0}}{2}}^{X}}}\circ \text{
}g_{r}^{X},\quad t\in \mathbb{R}_{+}^{*},$$
$$g_{t}^{X}=\underset{k-\text{times}}{\underbrace{g_{-\frac{\varepsilon _{0}}{2}}^{X}\circ
\cdots \circ g_{-\frac{\varepsilon _{0}}{2}}^{X}}}\circ \text{
}g_{r}^{X},\quad t\in \mathbb{R}_{-}^{*}.$$ The diffeomorphisms
$g_{\pm\frac{\varepsilon_{0}}{2}}^{X}$ and $g_{r}^{X}$ have been
defined above. Therefore, for all real $t$, $g_{t}^{X}$ is a
diffeomorphism defined globally on $M$ and the result is deduced
immediately. $\square$

With every vector field $X$ we associate the first-order
differential operator $L_{X}$. This is the differentiation of
functions in the direction of the vector field $X$. We have
$$L_{X}:\mathcal{C}^{\infty }\left( M\right)
\longrightarrow \mathcal{C}^{\infty }\left( M\right) ,\text{
}F\longmapsto L_{X}F,$$ where
$$L_{X}F(x)=\left. \frac{d}{dt}F\left( g_{t}^{X}(x)\right) \right| _{t=0},
\text{ }x\in M.$$ $\mathcal{C}^{\infty }\left( M\right)$ being the
set of functions $F:M\longrightarrow \mathbb{R}$, of class
$\mathcal{C}^{\infty }$. The operator $L_{X}$ is linear :
$$L_{X}\left(\alpha_{1}F_{1}+\alpha_{2}F_{2}\right)=\alpha_{1}L_{X}F_{1}+\alpha_{2}L_{X}F_{2},$$
where $\alpha_{1},\alpha_{2}\in \mathbb{R}$, and satisfies the
Leibniz formula :
$$L_{X}\left(F_{1}F_{2}\right)=F_{1}L_{X}F_{2}+F_{2}L_{X}F_{1}.$$
Since $L_{X}F(x)$ only depends on the values of $F$ in the
neighborhood of $x$, we can apply the operator $L_{X}$ to the
functions defined only in the neighborhood of a point, without the
need to extend them to the full variety $M$. Let
$\left(x_{1},...,x_{m}\right)$ a local coordinate system on $M$.
In this system the vector field $X$ has components $f_{1}, \ldots,
f_{m}$ and the flow $g_{t}^{X}$ is defined by a system of
differential equations. Therefore, the derivative of the function
$F=F\left(x_{1}, ...,x_{m}\right)$ in the direction of $X$ is
written $$L_{X}F=f_{1}\frac{\partial F}{\partial
x_{1}}+\cdots+f_{m}\frac{\partial F}{\partial x_{m}}.$$ In other
words, in the coordinates $\left( x_{1},...,x_{m}\right)$ the
operator $L_{X}$ is written $$L_{X}=f_{1}\frac{\partial }{\partial
x_{1}}+\cdots +f_{m}\frac{\partial }{\partial x_{m}},$$ this is
the general form of the first order linear differential operator.

Let $X$ be a vector field on a differentiable manifold $M$. We
have shown (theorem 7) that the vector field $X$ generates a
unique group of diffeomorphisms $g_t^X$ (that we also note $g_t$)
on $M$, solution of the differential equation
$$\frac{d}{dt}g_t^X(p)=X(g_t^X(p)),\quad p\in M,$$
with the initial condition $g_0^X(p)=p$. Let $\omega$ be a
$k$-differential form. The Lie derivative of $\omega$ with respect
to $X$ is the $k$-differential form defined by
$$
L_X\omega=\left.\frac{d}{dt}g_t^*\omega\right|_{t=0}=\lim_{t\rightarrow
0}\frac{g_t^*(\omega(g_t(p)))-\omega(p)}{t}.$$ In general, for
$t\neq0$, we have
\begin{equation}\label{eqn:euler}
\frac{d}{dt}g_t^*\omega=\left.\frac{d}{ds}g_{t+s}^*\omega\right|_{s=0}=
g_t^*\left.\frac{d}{ds}g_s^*\omega\right|_{s=0}=g_t^*(L_X\omega).
\end{equation}
We can easily verify that for the differential $k$-form
$\omega(g_t (p))$ at the point $g_t(p)$, the expression
$g_t^*\omega(g_t(p))$ is indeed a differential $k$-form in $p$.

For all $t\in \mathbb{R}$, the application
$g_t:\mathbb{R}\longrightarrow\mathbb{R}$ being a diffeomorphism
then $dg_t$ and $dg_{-t}$ are applications
\begin{eqnarray}
dg_t&:&T_pM\longrightarrow T_{g_t(p)}M,\nonumber\\
dg_{-t}&:&T_{g_t(p)}\longrightarrow T_pM.\nonumber
\end{eqnarray}
The Lie derivative of a vector field $Y$ in the direction $X$ is
defined by
$$
L_XY=\left.\frac{d}{dt}g_{-t}Y\right|_{t=0}=\lim_{t\rightarrow
0}\frac{g_{-t}(Y(g_t(p)))-Y(p)}{t}.$$ In general, for $t\neq0$, we
have
$$
\frac{d}{dt}g_{-t}Y=\left.\frac{d}{ds}g_{-t-s}Y\right|_{s=0}
=g_{-t}\left.\frac{d}{ds}g_{-s}Y\right|_{s=0}=g_{-t}(L_Y).$$ An
interesting operation on differential forms is the inner product
that is defined as follows : the inner product of a differential
$k$-form $\omega$ by a vector field $X$ on a differentiable
manifold $M$ is a differential $(k-1)$-form, denoted $i_X\omega$,
defined by
$$(i_X\omega)(X_1,...,X_{k-1})=\omega(X,X_1,...,X_{k-1}),$$
where $X_1,...,X_{k-1}$ are vector fields. It is easily shown that
if $\omega$ is a differential $k$-form, $\lambda$ a differential
form of any degree, $X$ and $Y$ two vector fields, $f$ a linear
map and $a$ a constant, then
\begin{eqnarray}
i_{X+Y}\omega&=&i_X\omega+i_Y\omega,\qquad\qquad
i_{aX}\omega=ai_X\omega,\nonumber\\
i_Xi_Y\omega&=&-i_Yi_X\omega,\qquad\qquad
i_Xi_X\omega=0,\nonumber\\
i_X(f\omega)&=&f(i_X\omega),\qquad\qquad
i_Xf^*\omega=f^*(i_{fX}\omega),\nonumber\\
i_X(\omega\wedge \lambda)&=&(i_X\omega)\wedge
\lambda+(-1)^k\omega\wedge (i_X\lambda).\nonumber
\end{eqnarray}

\begin{Exmp}
Let's calculate the expression of the inner product in local
coordinates. If
$$X=\sum_{j=1}^mX_j(x)\frac{\partial}{\partial x_j},$$
is the local expression of the vector field on the variety $M$ of
dimension $m$ and
$$\omega=\sum_{i_1<i_2<...<i_k}f_{i_1...i_k}(x)dx_{i_1}\wedge...\wedge dx_{i_k},$$
then,
\begin{eqnarray}
i_X\omega&=&\omega (X,\cdot),\nonumber\\
&=&\sum_{i_2<i_3<...<i_k}\sum_{j=1}^mf_{ji_2...i_k}X_jdx_{i_2}\wedge...\wedge dx_{i_k}\nonumber\\
&&\quad-\sum_{i_1<i_3<...<i_k}\sum_{j=1}^mf_{i_1j...i_k}X_jdx_{i_1}\wedge dx_{i_3}\wedge...\wedge dx_{i_k}\nonumber\\
&&\quad+\cdots+(-1)^{k-1}
\sum_{i_1<i_2<...<i_{k-1}}\sum_{j=1}^mf_{i_1i_2...j}X_jdx_{i_1}\wedge dx_{i_2}\wedge...\wedge dx_{i_{k-1}},\nonumber\\
&=&k\sum_{i_2<i_3<...<i_k}\sum_{j=1}^mf_{ji_2...i_k}X_jdx_{i_2}\wedge...\wedge
dx_{i_k}.\nonumber
\end{eqnarray}
Hence,
$$i_{\frac{\partial}{\partial
x_j}}\omega=\frac{\partial}{\partial(dx_j)}\omega,$$ where we put
$dx_j$ in first position in $\omega$.
\end{Exmp}

The following properties often occur when solving practical
problems using Lie derivatives.

\begin{Theo}
a) If $f:M\longrightarrow\mathbb{R}$ is a differentiable function,
so the Lie derivative of $f$ is the image of $X$ by the
differential of$f$,
$$L_Xf=df(X)=X.f.$$

b) $L_X$ and $d$ commute, $L_X\circ d=d\circ L_X$.

c) Let $X, X_1, ...,X_k$ be vector fields on $M$ and $\omega$ a
$k$-form differential. So
$$(L_X\omega)(X_1,...,X_k)=L_X(\omega(X_1,...,X_k))-\sum_{j=1}^k\omega(X_1,...,L_XX_j,...,X_k).$$
d) For all differential forms $\omega$ and $\lambda$, we have
$$L_X(\omega\wedge\lambda)=L_X\omega\wedge\lambda+\omega\wedge
L_X\lambda.$$
\end{Theo}
\emph{Proof}. a) Indeed, we have
$$
L_Xf=\left.\frac{d}{dt}g_t^*f\right|_{t=0}=\left.\frac{d}{dt}f\circ
g_t\right|_{t=0}
=\left.df\left(\frac{dg_t}{dt}\right)\right|_{t=0}=df(X),$$ and
(see further theorem 10),
$$L_Xf=i_Xdf=X.f.$$

b) Indeed, as the differential and the reciprocal image commute,
then
$$
d\circ
L_X\omega=\left.d\circ\frac{d}{dt}g_t^*\omega\right|_{t=0}=\left.\frac{d}{dt}g_t^*\circ
d\omega\right|_{t=0}=L_X\circ d\omega.$$

c) We have
\begin{eqnarray}
(L_X\omega)(X_1,...,X_k)&=&\left.\frac{d}{dt}g_t^*\omega (X_1,...,X_k)\right|_{t=0},\nonumber\\
&=&\left.\frac{d}{dt}\omega(g_t)(dg_tX_1,...,dg_tX_k)\right|_{t=0},\nonumber\\
&=&\left.L_X\omega (g_t)(dg_tX_1,...,dg_tX_k)\right|_{t=0}\nonumber\\
&&+\left.\sum_{j=1}^k\omega(g_t)\left(dg_tX_1,...,\frac{d}{dt}dg_tX_j,...,dg_tX_k\right)\right|_{t=0},\nonumber
\end{eqnarray}
and the result is deduced from the fact that
$$
\left.\frac{d}{dt}dg_tX_j\right|_{t=0}=-\left.\frac{d}{dt}dg_{-t}X_j\right|_{t=0}=-L_XX_j.$$

d) Just consider $\omega$ and $\lambda$ of the form
$$
\omega=fdx_{i_1}\wedge...\wedge dx_{i_k},$$$$
\lambda=gdx_{j_1}\wedge...\wedge dx_{j_l}.$$ We have
$$\omega\wedge\lambda=fgdx_{i_1}\wedge...\wedge dx_{i_k}\wedge dx_{j_1}\wedge...\wedge
dx_{j_l},$$ and
\begin{eqnarray}
&&L_X(\omega\wedge\lambda)(X_1,...,X_k,X_{k+1},...,X_{k+l})\nonumber\\
&&\quad=(L_Xf).gdx_{i_1}\wedge...\wedge dx_{i_k}\wedge
dx_{j_1}\wedge...\wedge
dx_{j_l}(X_1,...,X_k,X_{k+1},...,X_{k+l})\nonumber\\
&&\quad+f(L_Xg)dx_{i_1}\wedge...\wedge dx_{i_k}\wedge
dx_{j_1}\wedge...\wedge
dx_{j_l}(X_1,...,X_k,X_{k+1},...,X_{k+l}),\nonumber\\
&&\quad=((L_X\omega)\wedge\lambda+\omega\wedge
(L_X\lambda))(X_1,...,X_k,X_{k+1},...,X_{k+l}),\nonumber
\end{eqnarray}
and the result follows. $\square$

\begin{Theo}
Let $X$ and $Y$ be two vector fields on $M$. Then, the Lie
derivative of $L_XY$ is the Lie bracket $[X,Y]$.
\end{Theo}
\emph{Proof}. We have
$$
L_XY(f)=\lim_{t\rightarrow
0}\frac{dg_{-t}Y-Y}{t}(f)=\lim_{t\rightarrow
0}dg_{-t}\frac{Y-dg_tY}{t}(f),$$ hence
$$L_XY(f)=\lim_{t\rightarrow 0}\frac{Y(f)-dg_tY(f)}{t}=
\lim_{t\rightarrow 0}\frac{Y(f)-Y(f\circ g_t)\circ g_t^{-1}}{t}.
$$
Put $g_t(x)\equiv g(t, x)$, and apply to $g(t, x)$ the Taylor
formula with integral remainder. So there is $h(t, x)$ such that :
$$f(g(t,x))=f(x)+th(t,x),$$
with $$h(0,x)=\frac{\partial}{\partial t}f(g(t,x))(0,x).$$
According to the definition of the tangent vector, we have
$$X(f)=\frac{\partial}{\partial t}f\circ g_t(x)(0,x),$$
hence, $h(0,x)=X(f)(x)$. Therefore,
\begin{eqnarray}
L_XY(f)&=&\lim_{t\rightarrow 0}\left(\frac{Y(f)-Y(f)\circ g_t^{-1}}{t}-Y(h(t,x))\circ g_t^{-1}\right),\nonumber\\
&=&\lim_{t\rightarrow 0}\left(\frac{(Y(f)\circ g_t-Y(f))\circ
g_t^{-1}}{t}-Y(h(t,x))\circ g_t^{-1}\right).\nonumber
\end{eqnarray}
Since $$\lim_{t\rightarrow 0}g_t^{-1}(x)=g_0^{-1}(x)=id.,$$ we
deduce that :
\begin{eqnarray}
L_XY(f)&=&\lim_{t\rightarrow 0}\left(\frac{Y(f)\circ
g_t-Y(f)}{t}-Y(h(0,x))\right),\nonumber\\
&=&\frac{\partial}{\partial t}Y(f)\circ g_t(x)-Y(X(f)),\nonumber\\
&=&X(Y(f))-Y(X(f)),\nonumber\\
&=&[X,Y],\nonumber
\end{eqnarray}
which completes the demonstration. $\square$

We will now establish a fundamental formula for the Lie
derivative, which can be used as a definition.

\begin{Theo}
Let $X$ be a vector field on $M$ and $\omega$ a differential
$k$-form. Then
$$L_X\omega=d(i_X\omega)+i_X(d\omega).$$
In other words, we have the Cartan homotopy formula
$$L_X=d\circ i_X+i_X\circ d.$$
\end{Theo}
\emph{Proof}. We will reason by induction on the degree $k$ of the
differential form $\omega$. Let
$$D_X\equiv d\circ i_X+i_X\circ d.$$ For a differential $0$-form, i.e., a $f$ function, we have
$$D_Xf=d(i_Xf)+i_X(df).$$ Or $i_Xf=0$, hence
$$d(i_Xf)=0,\qquad i_Xdf=df(X),$$ and so $$D_Xf=df(X).$$
Moreover, we know (theorem 8, a)) that $L_Xf=df(X)=X.f$, so
$$D_Xf=L_Xf.$$ Assume that the formula in question is true for a differential $(k-1)$-form
and is proved to be true for a differential $k$-form. Let
$\lambda$ be a differential $(k-1)$-form and let $\omega=df\wedge
\lambda$, where $f$ is a function. We have
\begin{eqnarray}
L_X\omega&=&L_X(df\wedge\lambda),\nonumber\\
&=&L_Xdf\wedge\lambda+df\wedge L_X\lambda,\quad (\mbox{theorem 8, d)}),\nonumber\\
&=&dL_Xf\wedge\lambda+df\wedge L_X\lambda,\quad (\mbox{because }L_Xdf=dL_Xf, \mbox{theorem 8, b)}),\nonumber\\
&=&d(df(X))\wedge\lambda+df\wedge L_X\lambda,\quad (\mbox{because
}L_Xf=df(X), \mbox{theorem 8, a))}.\nonumber
\end{eqnarray}
By hypothesis of recurrence, we have
$$L_X\lambda=d(i_X\lambda)+i_X(d\lambda).$$ Or $i_Xdf=df(X)$, then
\begin{equation}\label{eqn:euler}
L_X\omega=d(i_Xdf)\wedge\lambda+df\wedge d(i_X\lambda)+df\wedge
i_X(d\lambda).
\end{equation}
Moreover, we have
$$
i_Xd\omega=i_Xd(df\wedge\lambda)=-i_X(df\wedge
d\lambda)=-(i_Xdf)d\wedge d\lambda+df\wedge i_X(d\lambda),
$$
and
\begin{eqnarray}
d(i_X\omega)&=&di_X(df\wedge\lambda),\nonumber\\
&=&d\left((i_Xdf)\wedge\lambda-df\wedge (i_X\lambda)\right),\nonumber\\
&=&d(i_Xdf)\wedge\lambda+(i_Xdf)\wedge d\lambda+df\wedge
d(i_X\lambda), \quad(\mbox{because } d(df)=0),\nonumber
\end{eqnarray}
hence,
$$di_X(df\wedge\lambda)+i_Xd(df\wedge\lambda)=d(i_Xdf)\wedge\lambda+df\wedge
d(i_X\lambda)+df\wedge i_X(d\lambda).$$ Comparing this expression
with that obtained in (3), we finally obtain
$$L_X\omega=d(i_X\omega)+i_X(d\omega),$$ and the theorem is proved.
$\square$

\begin{Exmp}
For a differential form $\omega$, we have
$$i_XL_X\omega=L_Xi_X\omega.$$
Indeed, we have
\begin{eqnarray}
i_XL_X\omega&=&i_X(di_X\omega)+i_X(i_Xd\omega),\quad (\mbox{theorem 10})\nonumber\\
&=&i_X(di_X\omega), \quad(\mbox{because } i_Xi_X=0),\nonumber
\end{eqnarray}
and
$$
L_Xi_X\omega=(d\circ i_X+i_X\circ d)i_X\omega=i_X(di_X\omega),$$
hence, $i_XL_X\omega-L_Xi_X\omega=0$.
\end{Exmp}

\begin{Theo}
Let $X$ and $Y$ be two vector fields on $M$ and $\omega$ a
differential form. Then,
\begin{eqnarray}
L_{X+Y}\omega&=&L_X\omega+L_Y\omega,\nonumber\\
L_{fX}\omega&=&fL_X\omega+df\wedge i_X\omega,\nonumber
\end{eqnarray}
where $f:M\longrightarrow \mathbb{R}$ is a differentiable
function.
\end{Theo}
\emph{Proof}. Indeed, just use the theorem 10,
\begin{eqnarray}
L_{X+Y}\omega&=&d(i_{X+Y}\omega)+i_{X+Y}(d\omega),\nonumber\\
&=&d(i_X\omega+i_Y\omega)+i_X(d\omega)+i_Y(d\omega),\nonumber\\
&=&d(i_X\omega)+i_X(d\omega)+d(i_Y\omega)+i_Y(d\omega),\nonumber\\
&=&L_X\omega+L_Y\omega.\nonumber
\end{eqnarray}
Similarly, we have
\begin{eqnarray}
L_{fX}\omega&=&d(i_{fX}\omega)+i_{fX}(d\omega),\nonumber\\
&=&d(fi_X\omega)+fi_X(d\omega),\nonumber\\
&=&df\wedge i_X\omega+fd(i_X\omega)+fi_X(d\omega),\nonumber\\
&=&df\wedge i_X\omega+fL_X\omega,\nonumber
\end{eqnarray}
which completes the demonstration. $\square$

\begin{Exmp}
Let's calculate the expression of the Lie derivative of the
differential form
$$\omega=\sum_{i_1<...<i_k}f_{i_1...i_k}dx_{i_1}\wedge...\wedge
dx_{i_k},$$ in local coordinates. If
$$X=\sum_{j=1}^mX_j(x)\frac{\partial}{\partial x_j},$$
is the local expression of the vector field on the $m$-dimensional
manifold $M$, then
$$
L_{X}\omega=\sum_{j=1}^mL_{X_j\frac{\partial}{\partial
x_j}}\omega=\sum_{j=1}^m\left(dX_j\wedge
i_{\frac{\partial}{\partial
x_j}}\omega+X_jL_{\frac{\partial}{\partial x_j}}\omega\right).$$
According to example 4, we know that
$$
i_{\frac{\partial}{\partial x_j}}\omega =\frac{\partial}{\partial
(dx_j)}\omega=k\sum_{i_2<i_3<...<i_k}f_{ji_2...i_k}dx_{i_2}\wedge
dx_{i_3}\wedge...\wedge dx_{i_k},$$ hence,
$$dX_j\wedge i_{\frac{\partial}{\partial x_j}}\omega=k\sum_{i_1<i_2<...<i_k}f_{ji_2...i_k}
\frac{\partial X_j}{\partial x_{i_1}}dx_{i_1}\wedge...\wedge
dx_{i_k}.$$ Similarly, using theorem 8, c), we obtain
$$L_{\frac{\partial}{\partial x_j}}\omega=\sum_{i_1<...<i_k}
\frac{\partial f_{i_1...i_k}}{\partial x_j}dx_{i_1}\wedge...\wedge
dx_{i_k}.$$ Since $\left[\frac{\partial}{\partial x_j},
\frac{\partial}{\partial x_l}\right]=0$, we finally get
$$L_X\omega=\sum_{i_1<...<i_k}\sum_{j=1}^m\left(\frac{\partial f_{i_1...i_k}}{\partial x_j}X_j
+kf_{ji_2...i_k}\frac{\partial X_j}{\partial
x_{i_1}}\right)dx_{i_1}\wedge...\wedge dx_{i_k}.$$
\end{Exmp}

\begin{Theo}
If $X$ and $Y$ are two vector fields on $M$, then

a) $\left[L_X, i_Y\right]=i_{[X, Y]}$.

b) $\left[L_X, L_Y\right]=L_{[X, Y]}$.
\end{Theo}
\emph{Proof}. The proof is to show that for a a differential
$k$-form $\omega$, we have $$\left[L_X, i_Y\right]\omega=i_{[X,
Y]}\omega,$$ and $$\left[L_X, L_Y\right]\omega=L_{[X, Y]}\omega.$$

$a)$ We reason by induction assuming first that $k=1$, that is,
$\omega=df$. We have
\begin{eqnarray}
[L_X, i_Y]df&=&L_Xi_Ydf-i_YL_Xdf,\nonumber\\
&=&L_X(Y.f)-i_YdL_Xf,\quad(\mbox{because } L_X\circ d=d\circ L_X)\nonumber\\
&=&X.(Y.f)-i_Yd(X.f),\quad(\mbox{because } L_Xf=X.f)\nonumber\\
&=&X.(Y.f)-Y.(X.f),\nonumber\\
&=&[X, Y].f,\nonumber\\
&=&i_{[X, Y]}df.\nonumber
\end{eqnarray}
Suppose the formula in question is true for a $\omega$ form of
degree less than or equal to $k-1$. Let $\lambda$ and $\theta$ be
two forms of degree less than or equal to $k-1$, so that
$\omega=\lambda\wedge\theta$ is a form of degree $k$. We have
\begin{eqnarray}
[L_X, i_Y]\omega&=&L_Xi_Y\omega-i_YL_X\omega,\nonumber\\
&=&L_Xi_Y(\lambda\wedge \theta)-i_YL_X(\lambda\wedge \theta),\nonumber\\
&=&L_X(i_Y\lambda\wedge \theta+(-1)^{deg \lambda}\lambda\wedge
i_Y\theta)
-i_Y(L_X\lambda\wedge \theta+\lambda\wedge L_X\theta),\nonumber\\
&=&L_Xi_Y\lambda\wedge \theta+i_Y\lambda\wedge L_X\theta
+(-1)^{deg \lambda}L_X\lambda\wedge i_Y\theta\nonumber\\
&&+(-1)^{deg \lambda}\lambda\wedge
L_Xi_Y\theta-i_YL_X\lambda\wedge \theta-(-1)^{deg
\lambda}L_X\lambda\wedge i_Y\theta\nonumber\\
&&
-i_Y\lambda\wedge L_X\theta-(-1)^{deg \lambda}\lambda\wedge i_YL_X\theta,\nonumber\\
&=&(L_Xi_Y\lambda-i_YL_X\lambda)\wedge \theta+(-1)^{deg \lambda}\lambda\wedge (L_Xi_Y\theta-i_YL_X\theta),\nonumber\\
&=&i_{[X, Y]}\lambda\wedge \theta+(-1)^{deg \lambda}\lambda\wedge i_{[X, Y]}\theta,\nonumber\\
&=&i_{[X, Y]}(\lambda\wedge \theta),\nonumber\\
&=&i_{[X, Y]}\omega.\nonumber
\end{eqnarray}

$b)$ We have
\begin{eqnarray}
[L_X, L_Y]\omega&=&L_XL_Y\omega-L_YL_X\omega,\nonumber\\
&=&L_Xdi_Y\omega+L_Xi_Yd\omega-di_YL_X\omega-i_YdL_X\omega,\nonumber\\
&=&dL_Xi_Y\omega+L_Xi_Yd\omega-di_YL_X\omega-i_YL_Xd\omega,(L_Xd\omega=dL_X\omega)\nonumber\\
&=&di_{[X, Y]}\omega+i_{[X, Y]}d\omega,(\mbox{according to a)})\nonumber\\
&=&L_{[X, Y]}\omega,\nonumber
\end{eqnarray}
and the demonstration ends. $\square$

\begin{Exmp}
Using the results seen above, we give a quick proof of
Poincar\'{e} lemma : in the neighborhood of a point of a manifold,
any closed differential form is exact. Indeed, consider the
differential equation
$$\dot{x}=X(x)=\frac{x}{t},$$ as well as its solution
$g_t(x_0)=x_0t$. The latter is defined in the neighborhood of the
point $x_0$, depends on $\mathcal{C}^\infty$ of the initial
condition and is a one-parameter group of diffeomorphisms. We have
$$g_0(x_0)=0, \quad g_1(x_0)=x_0,\quad g_0^*\omega=0,\quad
g_1^*\omega=\omega.$$ Hence,
$$
\omega=g_1^*\omega-g_0^*\omega=\int_0^1\frac{d}{dt}g_t^*\omega
dt=\int_0^1g_t^*(L_X\omega)dt\quad(\mbox{according to (2)}),$$ and
according to the theorem 10 and the fact that $d\omega=0$, we have
$$
\omega=\int_0^1g_t^*(di_X\omega)dt=\int_0^1dg_t^*i_X\omega
dt,\quad(\mbox{because } df^*\omega=f^*d\omega).$$ We can
therefore find a differential form $\lambda$ such that :
$\omega=d\lambda$, where $$\lambda=\int_0^1g_t^*i_X\omega dt.$$
\end{Exmp}

\section{The Darboux-Moser-Weinstein theorem}

\begin{Theo}
Let $\{\omega_t\}$, $0\leq t\leq 1$, be a family of symplectic
forms, differentiable in $t$. Then, for all $p\in M$, there exists
a neighborhood $\mathcal{U}$ of $p$ and a function $g_t :
\mathcal{U}\longrightarrow\mathcal{U}$, such that :
$g_0^*=\mbox{identity}$ et $g_t^*\omega_t=\omega_0$.
\end{Theo}
\emph{Proof}: Looking for a family of vector fields $X_t$ on
$\mathcal{U}$ such that these fields generate locally a
one-parameter group of diffeomorphisms $g_t$ with
$$
\frac{d}{dt}g_t(p)=X_t(g_t(p)),\quad g_0(p)=p.
$$
First note that the form $\omega_t$ is closed (i.e.,
$d\omega_t=0$) as the form $\frac{d}{d t}\omega_t$ (since
$d\frac{d}{d t}\omega_t=\frac{d}{d t}d\omega_t=0$). Therefore, by
deriving the relationship $g_t^*\omega_t=\omega_0$ and using the
Cartan homotopy formula (theorem 10) :
$L_{X_t}=i_{X_t}d+di_{X_t}$, taking into account that $\omega_t$
depends on time, we obtain the expression
$$
\frac{d}{dt}g_t^*\omega_t=g_t^*\left(\frac{d}{d
t}\omega_t+L_{X_t}\omega_t\right)=g_t^*\left(\frac{d}{d
t}\omega_t+di_{X_t}\omega_t\right).
$$
By Poincar\'{e}'s lemma (in the neighborhood of a point, any
closed differential form is exact), the form
$\frac{\partial}{\partial t}\omega_t$ is exact in the neighborhood
of $p$. In other words, we can find a form $\lambda_t$ such that :
$\frac{d}{dt}\omega_t=d\lambda_t$. Hence,
\begin{equation}\label{eqn:euler}
\frac{d}{dt}g_t^*\omega_t=g_t^*d(\lambda_t+i_{X_t}\omega_t).
\end{equation}
We want to show that for all $p\in M$, there exists a neighborhood
$\mathcal{U}$ of $p$ and a function $g_t :
\mathcal{U}\longrightarrow\mathcal{U}$, such that :
$g_0^*=\mbox{identity}$ and $g_t^*\omega_t=\omega_0$, therefore
such that : $\frac{d}{dt}g_t^*\omega_t=0$. And by (4), the problem
amounts to finding $X_t$ such that :
$\lambda_t+i_{X_t}\omega_t=0$. Since the form $\omega_t$ is non
degenerate, then the above equation is solvable with respect to
the vector field $X_t$ and defines the family $\{g_t\}$ for $0\leq
t\leq 1$. In local coordinates $(x_k)$ of the $2m$-dimensional
manifold $M$, with $\left(\frac {\partial}{\partial x_k}\right)$ a
basis of $TM$ and $(dx_k)$ the dual basis of
$\left(\frac{\partial}{\partial x_k}\right)$, $k=1,...,2m$, we
have
\begin{eqnarray}
\lambda_t&=&\sum_{k=1}^{2m}\lambda_k(t,x)dx_k,\nonumber\\
X_t&=&\sum_{k=1}^{2m}X_k(t,x)\frac{\partial}{\partial x_k},\nonumber\\
\omega_t&=&\sum_{\underset{k<l}{k,l=1}}^{2m}\omega_{k,l}(t,x)dx_k\wedge
dx_l,\nonumber\\
i_{X_t}\omega_t&=&2\sum_{l=1}^{2m}\left(\sum_{k=1}^{2m}\omega_{k,l}X_k\right)dx_l.\nonumber
\end{eqnarray}
We therefore solve the system of equations in $ x_k(t, x)$
according to :
$$\lambda_l(t,x)+2\sum_{k=1}^{2m}\omega_{k,l}(t,x)X_k(t,x)=0.$$
The form $\omega_t$ is non degenerate and the matrix $(\omega_
{k}(t,x))$ is nonsingular. Then the above system has a unique
solution. This determines the vector field $X_t$ and thus
functions $g_t^*$ such that : $g_t^*\omega_t=\omega_0$, which
completes the proof. $\square$

Using the above theorem (well known as Moser's lemma [31]), we
give a proof (Weinstein [39, 40]) of the Darboux theorem, which
states that every point in a symplectic manifold has a
neighborhood with Darboux coordinates. The Darboux theorem plays a
central role in symplectic geometry; the symplectic manifolds $(M,
\omega)$ of dimension $2m$ are locally isomorphic to
$(\mathbb{R}^{2m}, \omega)$. More precisely, if $(M,\omega)$ is a
symplectic manifold of dimension $2m$, then in the neighborhood of
each point of this manifold, there exist local coordinates $(x_1,
..., x_ {2m})$ such that : $$\omega=\sum_{k=1}^mdx_k\wedge
dx_{m+k}.$$ In particular, there is no local invariant in
symplectic geometry, analogous to the curvature in Riemannian
geometry. The classical proof given by Darboux is by induction on
the dimension of the manifold (see below and [4]).

\begin{Theo}
Any symplectic form on a manifold $M$ of dimension $2m$ is locally
diffeomorphic to the standard form on $\mathbb{R}^{2m}$. In other
words, if $(M,\omega)$ is a symplectic manifold of dimension $2m$,
then in the neighborhood of each point of $M$, there exist local
coordinates $(x_1, ..., x_ {2m})$ such that :
$$\omega=\sum_{k=1}^mdx_k\wedge dx_{m+k}.$$
\end{Theo}
\emph{Proof 1}. Let $\{\omega_t\}$, $0\leq t\leq 1$, be a family
of $2$-differential forms which depends differentiably on $t$ and
let $$\omega_t=\omega_0+t(\omega-\omega_0),
\qquad\omega_0=\sum_{k=1}^mdx_k\wedge d_{m+k},$$ where $(x_1, ...,
x_ {2m})$ are local coordinates on $M$. Note that these $2$-forms
are closed. At $p\in M$, we have
$$\omega_t(p)=\omega_0(p)=\omega(p).$$ By continuity, we can find a
small neighborhood of $p$ where the form $\omega_t(p)$ is non
degenerate. So the $2$-forms $\omega_t$ are non degenerate in a
neighborhood of $p$ and independent of $t$ at $p$. In other words,
$\omega_t$ are symplectic forms and by theorem 13, for all $p\in
M$, there exists a neighborhood $\mathcal{U}$ of $p$ and a
function $g_t:\mathcal{U}\longrightarrow\mathcal{U}$ such that :
$g_t^*=\mbox{identity}$ and $g_t^*\omega_t=\omega_0$.
Differentiating this relation with respect to $t$, we obtain (as
in the proof of theorem 13),
\begin{eqnarray}
\frac{d}{dt}g_t^*\omega_t&=&0,\nonumber\\
g_t^*\left(\frac{d}{dt}\omega_t+L_{X_t}\omega_t\right)&=&0,\nonumber\\
g_t^*\left(\frac{d}{dt}\omega_t+di_{X_t}\omega_t\right)&=&0.\nonumber
\end{eqnarray}
Therefore, $$di_{X_t}\omega_t=-\frac{d}{d t}\omega_t,$$ and since
the form $\frac{d}{dt}\omega_t$ is exact in the neighborhood of
$p$ (Poincar\'{e}'s lemma), then $$di_{X_t}\omega_t=d\theta_t,$$
where $\theta_t$ is a $1$-differential form. In addition,
$\omega_t$ being non degenerate, the equation
$i_{X_t}\omega_t=\theta_t$ is solvable and determines uniquely the
vector field $X_t$ depending on $t$. Note that for $ t=1$,
$\omega_1=\omega$ and for $t=0$, $\omega_0=\omega_0$ and also we
can find $g_1^*$ such that : $g_1^*\omega=\omega_0$. Vector fields
$X_t$ generate one-parameter families of diffeomorphisms
$\{g_t\}$, $0\leq t\leq 1$. In other words, you can make a change
of coordinates as : $$\omega=\sum_{k=1}^mdx_k\wedge d_{m+k},$$ and
the proof is completed.

\emph{Proof 2}. We proceed by induction on $m$. Suppose the result
true for $m-1\geq0$ and show that it is also for $m$. Fix $x$ and
let $x_{m +1}$ be a differentiable function on $M$ whose
differential $dx_ {m}$ is a nonzero point $x$. Let $X$ be the
unique differentiable vector field satisfying the relation
$i_X\omega=dx_{m+1}$. As this vector field does not vanish at $x$,
then we can find a function $x_1$ in a neighborhood $\mathcal{U}$
of $x$ such that $X(x_1)=1$. Consider a vector field $Y$ on
$\mathcal{U}$ satisfying the relation $i_Y\omega=-dx_1$. Since
$d\omega=0$, then $L_X\omega=L_Y\omega=0$, according to the Cartan
homotopy formula. therefore
$$i_{[X,Y]}\omega=L_Xi_Y\omega
=L_X(i_Y\omega)-i_Y(L_X\omega)=L_X(-dx_1)=-d(X(x_1))=0,$$ from
which we have $[X,Y]=0$, since any point in the form $\omega$ is
of rank equal to $2m$. By the Recovery theorem \footnote{Let $X_1,
...,X_r$ be differentiable vector fields on a manifold $M$ and
$x\in M$. Assume that for all $k,l=1,...,r$, $[x_k,X_l]=0$ and
$X_1(x), ..., X_r(x)$ are linearly independent. We show that there
is an open $\mathcal{U}$ of $ M $ containing $x$ and a local
coordinate system on $\mathcal{U}$ such that :
$X_1|_\mathcal{U}=\frac{\partial}{\partial
x_1},...,X_r|_\mathcal{U}=\frac{\partial}{\partial x_r}$.}, it
follows that there exist local coordinates $x_1, x_ {m +1}, z_1,
Z_2, ..., z_ {2m-2}$ on a neighborhood
$\mathcal{U}_1\subset\mathcal{U}$ of $x$ such that :
$X=\frac{\partial}{\partial x_1}$, $Y=\frac{\partial}{\partial
x_{m+1}}$. Consider the differential form
$$\lambda=\omega-dx_1\wedge dx_{m+1}.$$ We have $d\lambda=0$ and
$$i_X\lambda=L_X\lambda=i_Y\lambda=L_Y\lambda=0.$$ So $\lambda$ is
expressed as a $2$-differential form based only on variables $z_1,
z_2,...,z_{2m-2}$. In particular, we have $\lambda^{m+1}=0$.
Furthermore, we have $$0\neq\omega^m=mdx_1\wedge
dx_{m+1}\wedge\lambda^{m-1}.$$ The $2$-form $\lambda$ is closed
and of maximal rank (rank half) $m-1$ on an open set of
$\mathbb{R}^{2m-2}$. It is therefore sufficient to apply the
induction hypothesis to $\lambda$, which completes the proof.
$\square$

\begin{Rem}
If the variety $M$ is compact, connected and
$$\int_M\omega_t=\int_M\omega_0,$$ where $\{\omega_t\}$, $0\leq
t\leq 1$ is a family of volume forms, then one can find a family
of diffeomorphisms $g_t : M\longrightarrow M$, such that :
$g_0^*=\mbox{identity}$ and $g_t^*\omega_t=\omega_0$. Indeed, just
use a reasoning similar to theorem 7, provided to replace the
Poincar\'{e}'s lemma which is local, by the De Rham's theorem
which is global. This means that a volume form $\omega$ on $M$ is
exact if and only if $\int_M\omega=0$.
\end{Rem}

\section{Poisson brackets on symplectic manifolds and Hamiltonian systems}

As a consequence of the foregoing, the symplectic form $\omega$
induces a Hamiltonian vector field $$IdH:M\longrightarrow T_{x}M,
\quad x\longmapsto IdH(x),$$ where $H:M\longrightarrow
\mathbb{R}$, is a differentiable function (Hamiltonian). In others
words, the differential system defined by
$$\dot x(t)=X_{H}(x(t))=IdH(x),$$
is a Hamiltonian vector field associated to the function $H$. The
Hamiltonian vector fields form a Lie subalgebra of the vector
field space. The flow $g^t_X$ leaves invariant the symplectic form
$\omega$.

\begin{Theo}
The matrix that is associated to an Hamiltonian system determine a
symplectic structure.
\end{Theo}
\emph{Proof}. Let $\left(x_{1},\ldots ,x_{m}\right) $ be a local
coordinate system on $M,$ $(m=\dim M).$ We have
\begin{equation}\label{eqn:euler}
\dot x(t)=\sum_{k=1}^{n}\frac{\partial H} {\partial
x_{k}}I\left(dx_{k}\right)=\sum_{k=1}^{n}\frac{\partial
H}{\partial x_{k}}\xi^{k},
\end{equation}
where $I\left( dx_{k}\right)=\xi^{k}\in T_{x}M$ is defined such
that :
$$\forall \eta \in T_{x}M,\text{ }\eta _{k}=dx_{k}\left( \eta \right)=
\omega \left( \eta ,\xi ^{k}\right),\quad(k^{th}\text{-component
of }\eta).$$ Define $\left( \eta _{1},\ldots ,\eta _{m}\right) $
and $\left(\xi_{1}^{k},\ldots,\xi_{m}^{k}\right) $ to be
respectively the components of $\eta $ and $\xi^{k}$, then
$$
\eta_{k}=\sum_{i=1}^{m}\eta_{i}\left(\frac{\partial}{\partial
x_{i}}, \frac{\partial }{\partial x_{j}}\right)
\xi_{j}^{k}=\left(\eta_{1},\ldots ,\eta_{m}\right) J^{-1}
\left(\begin{array}{c}
\xi_{1}^{k}\\
\vdots \\
\xi_{m}^{k}
\end{array}
\right),
$$
where $J^{-1}$ is the matrix defined by
$$J^{-1}\equiv
\left(\omega\left(\frac{\partial}{\partial x_{i}},
\frac{\partial}{\partial x_{j}}\right)\right)_{1\leq i,j\leq m}.$$
Note that this matrix is invertible. Indeed, it suffices to show
that the matrix $J^{-1}$ has maximal rank. Suppose this were not
possible, i.e., we assume that $rank(J^{-1})\neq m$. Hence
$$\sum_{i=1}^{m}a_{i}\omega \left(\frac{\partial}{\partial
x_{i}},\frac{\partial}{\partial x_{j}}\right)=0, \quad\forall
1\leq j\leq m,$$ with $a_{i}$ not all null and
$$\omega\left(\sum_{i=1}^{m}a_{i}\frac{\partial }{\partial
x_{i}},\frac{\partial}{\partial x_{j}}\right)=0, \quad\forall
1\leq j\leq m.$$ In fact, since $\omega$ is non-degenerate, we
have $\sum_{i=1}^{m}a_{i}\frac{\partial}{\partial x_{i}}=0$. Now
$\left(\frac{\partial}{\partial
x_{1}},\ldots,\frac{\partial}{\partial x_{m}}\right)$ is a basis
of $T_{x}M,$ then $a_{i}=0$, $\forall i$, contradiction. Since
this matrix is invertible, we can search $\xi^{k}$ such that :
$$J^{-1}\left(\begin{array}{c}
\xi_{1}^{k}\\
\vdots\\
\xi_{m}^{k}
\end{array}
\right)=\left(\begin{array}{cc}
0&\\
\vdots \\
0&\\
1&\leftrightsquigarrow k^{th}\text{-place}\\
0&\\
\vdots \\
0&
\end{array}
\right).$$ The matrix $J^{-1}$ is invertible, which implies
$$\left(\begin{array}{c}
\xi_{1}^{k}\\
\vdots \\
\xi_{m}^{k}
\end{array}
\right)=J\left(\begin{array}{c}
0\\
\vdots \\
0\\
1\\
0\\
\vdots \\
0
\end{array}
\right),$$ from which $\xi^{k}=(k^{th}$-column of $J$), i.e.,
$\xi_{i}^{k}=J_{ik}$, $1\leq i\leq m$, and consequently
$$\xi^{k}=\sum_{i=1}^{m}J_{ik}\frac{\partial}{\partial
x_{i}}.$$ It is easily verified that the matrix $J$ is
skew-symmetric\footnote{Indeed, since
$\omega\left(\frac{\partial}{\partial
x_{i}},\frac{\partial}{\partial x_{j}}\right)=-\omega
\left(\frac{\partial}{\partial x_{j}},\frac{\partial}{\partial
x_{i}}\right)$, i.e., $\omega$ is symmetric, it follows that
$J^{-1}$ is skew-symmetric. Then,
$I=J.J^{-1}=\left(J^{-1}\right)^{\top }.J^{\top }=-J^{-1}.J$ and
consequently $J^{\top}=J$}. From (3) we deduce that
$$\dot x(t)=\sum_{k=1}^{m}\frac{\partial
H}{\partial x_{k}}\sum_{i=1}^{m}J_{ik}\frac{\partial}{\partial
x_{i}}=\sum_{i=1}^{m}\left(\sum_{k=1}^{m}J_{ik}\frac{\partial
H}{\partial x_{k}}\right)\frac{\partial}{\partial x_{i}}.$$
Writing $$\dot
x(t)=\sum_{i=1}^{m}\frac{dx_{i}(t)}{dt}\frac{\partial}{\partial
x_{i}},$$ it is seen that $$\dot
x_{i}(t)=\sum_{k=1}^{m}J_{ik}\frac{\partial H}{\partial x_{k}},
\quad1\leq i\leq j\leq m$$ which can be written in more compact
form
$$\dot x(t)=J(x)\frac{\partial H}{\partial x},$$ this is the Hamiltonian
vector field associated to the function $H$. $\square$

Let $(M, \omega)$ be a symplectic manifold. To any pair of
differentiable functions $(F, G)$ over $M$, we associate the
function
$$\{F,G\}=d_{u}F(X_{G})=X_{G}F(u)=\omega(X_{G},X_{F}),$$
where $X_{F}$ and $X_{G}$ are the Hamiltonian vector fields
associated with the functions $F$ and $G$ respectively. We say
that $\{F,G\}$ is a Poisson bracket (or Poisson structure) of the
functions $F$ and $G$. It is easily verified that the Poisson
bracket on the space $\mathcal{C}^{\infty }$, i.e., the bilinear
application
$$\{ ,\} :\mathcal{C}^{\infty }(M)
\times \mathcal{C}^{\infty }(M) \longrightarrow
\mathcal{C}^{\infty}(M),\text{ }(F,G) \longmapsto \{F,G\},$$
defined above (where $\mathcal{C}^{\infty }(M)$ is the commutative
algebra of regular functions on $M$) is skew-symmetric
$\{F,G\}=-\{G,F\}$, obeys the Leibniz rule
$$\{FG,H\}=F\{G,H\}+G\{F,H\},$$ and satisfies the Jacobi identity
$$\{\{H,F\},G\}+\{\{F,G\},H\}+\{\{G,H\},F\}=0.$$ The variety $M$ is
called a Poisson manifold or a Hamiltonian variety. The Leibniz
formula ensures that the mapping $G\longmapsto \{G,F\}$ is a
derivation. The antisymmetry and identity of Jacobi ensure that
$\{, \}$ is a Lie bracket, they provide $\mathcal{C}^{\infty}(M)$
of an infinite-dimensional Lie algebra structure. When this
Poisson structure is non-degenerate, we obtain the symplectic
structure discussed above.

Consider now $M=\mathbb{R}^{n}\times \mathbb{R}^{n}$ and let $p\in
M$. By Darboux's theorem, there exists a local coordinate system
$\left(x_{1},\ldots,x_{n}, y_{1},\ldots,y_{n}\right)$ in a
neighborhood of $p$ such that
$$\{H,F\}=\sum_{i=1}^{n}\left(\frac{\partial H}{\partial
x_{i}}\frac{\partial F}{\partial y_{i}}-\frac{\partial H}{\partial
y_{i}}\frac{\partial F}{\partial x_{i}}\right)=X_{H},$$ and
$$X_{H}F=\{H,F\}, \quad\forall F\in \mathcal{C}^{\infty}(M)$$ The
manifold $M$ with the local coordinates $y_{1},\ldots, y_{n},
x_{1},\ldots, x_{n}$ and the the above mentioned canonical Poisson
bracket is a Poisson manifold. The Hamiltonian systems form a Lie
algebra. A nonconstant function $F$ is called an integral (first
integral or constant of motion) of $X_F$, if $X_{H}F=0$; this
means that $ F $ is constant on the trajectories of $X_{H}$. In
particular, $H$ is integral. Two functions $F$ and $G$ are said to
be in involution or to commute, if $\{F,G\}=0$. An interesting
result is given by the following Poisson theorem :

\begin{Theo}
If $F$ and $G$ are two first integrals of a Hamiltonian system,
then $\{F, G\}$ is also a first integral.
\end{Theo}
\emph{Proof}. Jacobi's identity is written
$$\left\{\left\{H,F\right\},G\right\}+\left\{\left\{F,G\right\},
H\right\}+\left\{\left\{G,H\right\},F\right\}=0,$$ where $H$ is
the Hamiltonian. Since $\{H,F\}=\{H,G\}=0$, then we have
$\left\{\left\{F,G\right\}, H\right\}=0$, which shows that $\{F,
G\}$ is a first integral. $\square$

\begin{Rem}
If we know two first integrals, we can, according to Poisson's
theorem, find new integrals. But let's mention that we often fall
back on known first integrals or a constant.
\end{Rem}

Let $M$ and $N$ two differentiable manifolds and
$f\in\mathcal{C}^\infty (M,N)$. The linear tangent map to $f$ at
the point $p$ is the induced mapping between the tangent spaces
$T_pM $ and $T_{f(p)}N$, defined by
$$f_* :
T_pM\longrightarrow T_{f(p)}N,\quad f_*v(\varphi)=v(\varphi\circ
f),$$ where $v\in T_pM$ and $\varphi\in\mathcal{C}^\infty (N,
\mathbb{R})$. Let $L : TM\longrightarrow \mathbb{R}$ be a
differentiable function (Lagrangian\index{Lagrangian}) on the
tangent bundle $TM$. We say that $(M, L)$ is invariant under the
differentiable application $g: M\longrightarrow M$ if for all
$v\in TM$, we have $$L(g_*v)=L(v).$$ The theorem of Noether below,
expresses the existence of a first integral associated with a
symmetry of the Lagrangian. In other words, each parameter of a
group of transformations corresponds to a conserved quantity. One
of the consequences of the invariance of the Lagrangian with
respect to a group of transformations is the conservation of the
generators of the group. For example, the first integral
associated with rotation invariance is the kinetic moment.
Similarly, the first integral associated with the invariance with
respect to the translations is the pulse. The Noether theorem
applies to certain classes of theories, described either by a
Lagrangian or a Hamiltonian. We will give below the theorem in its
original version, which applies to the theories described by a
Lagrangian. There is also a version that applies to theories
described by a Hamiltonian.

\begin{Theo}
If $(M, L)$ is invariant under a parameter group of
diffeomorphisms $g_{s}: M\longrightarrow M$,  $s\in\mathbb{R}$,
$g_0=E$, then the system of Lagrange equations\index{system of
Lagrange equations}, $$\frac{d}{dt}\frac{\partial L}{\partial
\overset{.}{q}}=\frac{\partial L}{\partial q},$$ corresponding to
$ L $ admits a first integral $I: TM\longrightarrow\mathbb{R}$
with $$\left.I(q, \overset{.}{q})=\frac{\partial L}{\partial
\overset{.}{q}}\frac{dg_{s}(q)}{ds}\right|_{s=0},$$ the $q$ being
local coordinates on $M$.
\end{Theo}
\emph{Proof}. The first integral $I$ is independent of the choice
of local coordinates $q$ over $M$ and so we can just consider the
case $M=\mathbb{R}^n$. Let $$f : \mathbb{R}\longrightarrow M,
\quad t\longmapsto q=f(t),$$ be a solution of the system of
Lagrange equations above. By hypothesis, $g_{*s}$ leaves $L$
invariant, so $$g_s\circ f : \mathbb{R}\longrightarrow M, \quad
t\longmapsto g_s\circ f(t),$$ also satisfies the system of
Lagrange equations. We translate the solution $f(t)$ considering
the application $$F :
\mathbb{R}\times\mathbb{R}\longrightarrow\mathbb{R}^n, \quad
(s,t)\longmapsto q=g_s(f(t)).$$ The fact that $g_s$ leaves
invariant $L$ implies that : $$0=\frac{\partial L (F,
\overset{.}{F})}{\partial s}=\frac{\partial L}{\partial
q}\frac{\partial F}{\partial s}+\frac{\partial L}{\partial
\overset{.}{q}}\frac{\partial \overset{.}{F}}{\partial s},$$ i.e.,
\begin{equation}\label{eqn:euler}
\frac{\partial L}{\partial q}\frac{\partial q}{\partial
s}+\frac{\partial L}{\partial \overset{.}{q}}\frac{\partial
\overset{.}{q}}{\partial s}=0.
\end{equation}
Since $F$ is also a solution of the system of Lagrange equations,
i.e.,
$$\frac{d}{dt}\left(\frac{\partial L}{\partial \overset{.}{q}}\left(F(s,t), \overset{.}{F}(s,t)\right)\right)
= \frac{\partial L}{\partial
q}\left(F(s,t),\overset{.}{F}(s,t)\right),$$ so noting that :
$$\frac{\partial \overset{.}{q}}{\partial
s}=\frac{d}{dt}\frac{\partial q}{\partial s},$$ and equation (6)
is written in the form
\begin{eqnarray}
0&=&\frac{\partial q}{\partial s}\frac{d}{dt}\left(\frac{\partial
L}{\partial \overset{.}{q}}
\left(F(s,t),\overset{.}{F}(s,t)\right)\right)+\frac{\partial
L}{\partial\overset{.}{q}}\frac{d}{dt}\frac{\partial q}{\partial s},\nonumber\\
&=&\frac{d}{dt}\left(\frac{\partial L}{\partial
\overset{.}{q}}\left(F(s,t),\overset{.}{F}(s,t)\right)\frac{\partial
q}{\partial s}\right),\nonumber
\end{eqnarray}
which completes the proof of the theorem. $\square$

We now give the following definition of the Poisson bracket :
$$
\{F,G\}=\left\langle \frac{\partial F}{\partial x},
J\frac{\partial G}{\partial x}\right\rangle
=\sum_{i,j}J_{ij}\frac{\partial F}{\partial x_{i}}\frac{\partial
G}{\partial x_{j}}.$$ We will look for conditions on the matrix
$J$ for Jacobi's identity to be satisfied. This is the purpose of
the following theorem :

\begin{Theo}
The matrix $J$ satisfies the Jacobi identity, if
$$\sum_{k=1}^{2n}\left(J_{kj}\frac{\partial J_{li}}
{\partial x_{k}}+J_{ki}\frac{\partial J_{jl}} {\partial
x_{k}}+J_{kl}\frac{\partial J_{ij}}{\partial
x_{k}}\right)=0,\text{ }\forall 1\leq i,j,l\leq 2n.$$
\end{Theo}
\emph{Proof}. Consider the Jacobi identity :
$$\{\{H,F\},G\}+\{\{F,G\},H\}+\{\{G,H\},F\}=0.$$ We have
$$\{\{H,F\},G\}=\left\langle \frac{\partial \left\{H,F\right\}
}{\partial x},J\frac{\partial G}{\partial x}\right\rangle
=\sum_{k,l}J_{kl}\frac{\partial \{H,F\} }{\partial
x_{k}}\frac{\partial G}{\partial x_{l}},$$ hence
\begin{eqnarray}
\{\{H,F\},G\}&=&\sum_{k,l}\sum_{i,j}J_{kl}\frac{\partial
J_{ij}}{\partial x_{k}}\frac{\partial H}{\partial
x_{i}}\frac{\partial F}{\partial x_{j}}\frac{\partial G}{\partial
x_{l}}+\sum_{k,l}\sum_{i,j}J_{kl}J_{ij}\frac{\partial
^{2}H}{\partial x_{k}\partial x_{i}}\frac{\partial F}{\partial
x_{j}}\frac{\partial G}{\partial x_{l}}\nonumber\\
&&+\sum_{k,l}\sum_{i,j}J_{kl}J_{ij}\frac{\partial H}{\partial
x_{i}}\frac{\partial ^{2}F}{\partial x_{k}\partial
x_{j}}\frac{\partial G}{\partial x_{l}}.\nonumber
\end{eqnarray}
By symmetry, we have immediately $\{\{F,G\},H\}$ and
$\{\{G,H\},F\}$. Then
\begin{eqnarray}
&&\{\{H,F\},G\}+\{\{F,G\},H\}+\{\{G,H\},F\}\nonumber\\
&&\qquad\qquad\qquad\qquad=\sum_{k,l}\sum_{i,j}J_{kl}\frac{\partial
J_{ij}}{\partial x_{k}}\frac{\partial H}{\partial
x_{i}}\frac{\partial F}{\partial
x_{j}}\frac{\partial G}{\partial x_{l}}\nonumber\\
&&\qquad\qquad\qquad\qquad+\sum_{k,l}\sum_{i,j}J_{kl}J_{ij}\frac{\partial
^{2}H}{\partial x_{k}\partial x_{i}}\frac{\partial F}{\partial
x_{j}}\frac{\partial G}{\partial x_{l}}\\
&&\qquad\qquad\qquad\qquad+\sum_{k,l}\sum_{i,j}J_{kl}J_{ij}\frac{\partial
H}{\partial x_{i}}\frac{\partial ^{2}F}
{\partial x_{k}\partial x_{j}}\frac{\partial G}{\partial x_{l}}\\
&&\qquad\qquad\qquad\qquad+\sum_{k,l}\sum_{i,j}J_{kl}\frac{\partial
J_{ij}}{\partial x_{k}}\frac{\partial G}
{\partial x_{i}}\frac{\partial H}{\partial x_{j}}\frac{\partial F}{\partial x_{l}}\nonumber\\
&&\qquad\qquad\qquad\qquad+\sum_{k,l}\sum_{i,j}J_{kl}J_{ij}\frac{\partial
^{2}G}{\partial x_{k}\partial x_{i}}
\frac{\partial H}{\partial x_{j}}\frac{\partial F}{\partial x_{l}}\\
&&\qquad\qquad\qquad\qquad+\sum_{k,l}\sum_{i,j}J_{kl}J_{ij}\frac{\partial
G}{\partial x_{i}}
\frac{\partial ^{2}H}{\partial x_{k}\partial x_{j}}\frac{\partial F}{\partial x_{l}}\\
&&\qquad\qquad\qquad\qquad+\sum_{k,l}\sum_{i,j}J_{kl}\frac{\partial
J_{ij}}{\partial x_{k}}\frac{\partial F}
{\partial x_{i}}\frac{\partial G}{\partial x_{j}}\frac{\partial H}{\partial x_{l}}\nonumber\\
&&\qquad\qquad\qquad\qquad+\sum_{k,l}\sum_{i,j}J_{kl}J_{ij}\frac{\partial
^{2}F}{\partial x_{k}\partial x_{i}}
\frac{\partial G}{\partial x_{j}}\frac{\partial H}{\partial x_{l}}\\
&&\qquad\qquad\qquad\qquad+\sum_{k,l}\sum_{i,j}J_{kl}J_{ij}\frac{\partial
F}{\partial x_{i}}\frac{\partial ^{2}G} {\partial x_{k}\partial
x_{j}}\frac{\partial H}{\partial x_{l}}.
\end{eqnarray}
Notice that the indices $i,j,k$ and $l$ play a symmetric roll.
Applying in the term (10) the permutation $i\leftarrow l$,
$j\leftarrow k$, $k\leftarrow i$, $l\leftarrow j$, and add the
term (7), with the understanding that $J_{lk}=-J_{kl}$, we get
$$\sum_{k,l}\sum_{i,j}\left(J_{ij}J_{lk}+J_{kl}J_{ij}\right)
\frac{\partial G}{\partial x_{l}}\frac{\partial^{2}H} {\partial
x_{i}\partial x_{k}}\frac{\partial F}{\partial x_{j}}=0,$$ as a
consequence of the Schwarz's lemma. Again applying in the term
(11) the permutation $i\leftarrow k$, $j\leftarrow l$,
$k\leftarrow j$, $l\leftarrow i$, and add the term (8), yields
$$\sum_{k,l}\sum_{i,j}\left(J_{ji}J_{kl}+J_{kl}J_{ij}\right)
\frac{\partial^{2}F}{\partial x_{j}\partial x_{k}}\frac{\partial
G}{\partial x_{l}}\frac{\partial H}{\partial x_{i}}=0.$$ By the
same argument as above, applying in the term (12) the permutation
$i\leftarrow l$, $j\leftarrow k$, $k\leftarrow i$, $l\leftarrow
j$, and add the term (9), we obtain
$$\sum_{k,l}\sum_{i,j}\left(J_{ij}J_{lk}+J_{kl}J_{ij}\right)
\frac{\partial F}{\partial x_{l}}\frac{\partial^{2}G} {\partial
x_{i}\partial x_{k}}\frac{\partial H}{\partial x_{j}}=0,$$ and
thus
\begin{eqnarray}
&&\{\{ H,F\},G\}+\{\{F,G\},H\}+\{\{G,H\},F\}\nonumber\\
&&\qquad\qquad\qquad\qquad=\sum_{k,l}\sum_{i,j}J_{kl}\frac{\partial
J_{ij}}{\partial x_{k}}\frac{\partial H}{\partial
x_{i}}\frac{\partial F}{\partial
x_{j}}\frac{\partial G}{\partial x_{l}}\\
&&\qquad\qquad\qquad\qquad+\sum_{k,l}\sum_{i,j}J_{kl}\frac{\partial
J_{ij}}{\partial x_{k}}\frac{\partial G}{\partial
x_{i}}\frac{\partial H}{\partial
x_{j}}\frac{\partial F}{\partial x_{l}}\\
&&\qquad\qquad\qquad\qquad+\sum_{k,l}\sum_{i,j}J_{kl}\frac{\partial
J_{ij}}{\partial x_{k}}\frac{\partial F}{\partial
x_{i}}\frac{\partial G}{\partial x_{j}}\frac{\partial H}{\partial
x_{l}}.\nonumber
\end{eqnarray}
Under permuting the indices $i\leftarrow l,$ $j\leftarrow i$,
$k\leftarrow k$, $l\leftarrow j$, for (13) and $i\leftarrow j$,
$j\leftarrow l$, $k\leftarrow k$, $l\leftarrow i$, for (14), we
obtain the following :
\begin{eqnarray}
&&\{\{H,F\},G\}+\{\{F,G\},H\}+\{\{G,H\},F\}\nonumber\\
&&\qquad\qquad=\sum_{i,j,l}\left[\sum_{k}\left(J_{kj}\frac{\partial
J_{li}}{\partial x_{k}}+J_{ki}\frac{\partial J_{jl}}{\partial
x_{k}}+J_{kl}\frac{\partial J_{ij}}{\partial x_{k}}\right)\right]
\frac{\partial H}{\partial x_{l}}\frac{\partial F}{\partial
x_{i}}\frac{\partial G}{\partial x_{j}}.\nonumber
\end{eqnarray}
Since the Jacobi identity must be identically zero, then the
expression to prove follows immediately, ending the proof of
theorem. $\square$

Consequently, we have a complete characterization of Hamiltonian
vector field
\begin{equation}\label{eqn:euler}
\dot x(t)=X_{H}(x(t))=J\frac{\partial H}{\partial x},\text{ }x\in
M,
\end{equation}
where $H:M\longrightarrow \mathbb{R}$, is the Hamiltonian and
$J=J(x)$ is a skew-symmetric matrix, for which the corresponding
Poisson bracket satisfies the Jacobi identity :
$$\{\{H,F\},G\}+\{\{F,G\},H\}+\{\{G,H\},F\}=0,$$ with
$$\{H,F\}=\left\langle
\frac{\partial H}{\partial x},J\frac{\partial F}{\partial
x}\right\rangle =\sum_{i,j}J_{ij}\frac{\partial H}{\partial x_{i}}
\frac{\partial F}{\partial x_{j}},\mbox{ (Poisson bracket)}.$$

\section{Examples}

\begin{Exmp}
An important special case is when $$J=\left(\begin{array}{cc}
O&-I\\
I&O
\end{array}\right),$$
where $I$ is the $n\times n$ identity matrix. The condition on $J$
is trivially satisfied. Indeed, here the matrix $J$ do not depend
on the variable $x$ and we have
$$
\left\{H,F\right\} =\sum_{i=1}^{2n}\frac{\partial H}{\partial
x_{i}}\sum_{j=1}^{2n}J_{ij}\frac{\partial F}{\partial
x_{j}}=\sum_{i=1}^{n}\left(\frac{\partial H}{\partial
x_{n+i}}\frac{\partial F}{\partial x_{i}}-\frac{\partial
H}{\partial x_{i}}\frac{\partial F}{\partial x_{n+i}}\right) .$$
Moreover, equations (15) are transformed into
$$\dot q_{1}=\frac{\partial H}{\partial p_{1}},\ldots ,\dot q_{n}
=\frac{\partial H}{\partial p_{n}},\text{ }\dot
p_{1}=-\frac{\partial H} {\partial q_{1}},\ldots ,\dot
p_{n}=-\frac{\partial H}{\partial q_{n}},$$ o\`{u}
$q_{1}=x_{1},\ldots ,q_{n}=x_{n},p_{1}=x_{n+1},\ldots
,p_{n}=x_{2n}$. These are exactly the well known differential
equations of classical mechanics in canonical form. They show that
it suffices to know the Hamiltonian function $H$ to determine the
equations of motion. They are often interpreted by considering
that the variables $p_{k}$ and $q_{k}$ are the coordinates of a
point that moves in a space with $2n$ dimensions, called phase
space. The flow associated with the system above obviously leaves
invariant each hypersurface of constant energy $H=c$. The Hamilton
equations above, can still be written in the form
$$
\dot{q}_{i}=\{H, q_i\}=\frac{\partial H}{\partial p_{i}},\qquad
\dot{p}_{i}=\{H, p_i\}=-\frac{\partial H}{\partial q_{i}},$$ where
$1\leq i\leq n$. Note that the functions $1, q_i, p_i$ ($1\leq
i\leq n$), verify the following commutation relations :
$$\{q_i, q_j\}=\{p_i, p_j\}=\{q_i, 1\}=\{p_i, 1\}=0,\quad \{p_i,
q_j\}=\delta_{ij},\quad 1\leq i,j\leq n.$$ These functions
constitute a basis of a real Lie algebra (Heisenberg algebra), of
dimension $2n+1$.
\end{Exmp}

\begin{Exmp}
The H\'{e}non-Heiles differential equations are defined by
\begin{eqnarray}
\dot{y}_{1}&=&x_{1},\qquad \dot{x}_{1}=-Ay_{1}-2y_{1}y_{2},\nonumber\\
\dot{y}_{2}&=&x_{2},\qquad\dot{x}_{2}=-By_{2}-y_{1}^{2}-\varepsilon
y_{2}^{2},\nonumber
\end{eqnarray}
where $A, B, \varepsilon$ are constants. The above equations can
be rewritten as a Hamiltonian vector field
$$\dot{x}=J\frac{\partial H}{\partial x},
\quad x=(y_{1},y_{2},x_{1},x_{2})^\top,$$ where
$$
H=\frac{1}{2}(x_{1}^{2}+x_{2}^{2}+Ay_{1}^{2}+By_{2}^{2})
+y_{1}^{2}y_{2}+\frac{\varepsilon}{3}y_{2}^{3},\quad\mbox{(Hamiltonian)}$$
and $J=\left(\begin{array}{cc}
0&-I\\
I&0
\end{array}\right)$,
is the matrix associated with the vector field.
\end{Exmp}

\begin{Exmp}
The Euler equations of the rotation motion of a solid around a
fixed point, taken as the origin of the reference bound to the
solid, when no external force is applied to the system, can be
written in the form :
\begin{eqnarray}
\dot{m}_{1}&=&\left(\lambda_{3}-\lambda_{2}\right)m_{2}m_{3},\nonumber\\
\dot{m}_{2}&=&\left( \lambda_{1}-\lambda _{3}\right) m_{1}m_{3},\nonumber\\
\dot{m}_{3}&=&\left(\lambda_{2}-\lambda_{1}\right)m_{1}m_{2},\nonumber
\end{eqnarray}
where $(m_{1},m_{2},m_{3})$ is the angular momentum of the solid
and $\lambda _{i}\equiv I_{i}^{-1}$, $I_{1},I_{2}$ et $I_{3}$
being moments of inertia. These equations can be written in the
form of a Hamiltonian vector field : $$\dot{x}=J\frac{\partial
H}{\partial x}, \quad x=\left(m_{1},m_{2},m_{3}\right)^\top,$$
with
$$H=\frac{1}{2}\left(\lambda_{1}m_{1}^{2}+\lambda_{2}m_{2}^{2}+\lambda_{3}m_{3}^{2}\right),\quad\mbox{(Hamiltonian)}$$
To determine the matrix $J=(J_{ij})_{1\leq i,j\leq 3}$, we proceed
as follows : since $J$ is antisymmetric, then obviously
$J_{ii}=0$, $J_{ij}=-J_{ji}$, $1\leq i,j\leq 3$, hence
$$J=\left(\begin{array}{ccc}
0&J_{12}&J_{13}\\
-J_{12}&0&J_{23}\\
-J_{13}&-J_{23}&0
\end{array}\right).
$$
Therefore,
\begin{eqnarray}
\left(\begin{array}{c}
\dot{m}_{1}\\
\dot{m}_{2}\\
\dot{m}_{3}
\end{array}\right)&=&\left(\begin{array}{ccc}
0&J_{12}&J_{13}\\
-J_{12}&0&J_{23}\\
-J_{13}&-J_{23}&0
\end{array}\right)
\left(\begin{array}{c}
\lambda_1m_1\\
\lambda_2m_2\\
\lambda_3m_3
\end{array}\right),\\
&=&\left(\begin{array}{c}
\left(\lambda_{3}-\lambda _{2}\right)m_{2}m_{3}\\
\left(\lambda_{1}-\lambda _{3}\right)m_{1}m_{3}\\
\left(\lambda_{2}-\lambda _{1}\right)m_{1}m_{2}
\end{array}\right).
\end{eqnarray}
Comparing (16) and (17), we deduce that : $J_{12}=-m_3$,
$J_{13}=m_2$ and $J_{23}=-m_1$. Finally,
$$J=\left(\begin{array}{ccc}
0&-m_{3}&m_{2}\\
m_{3}&0&-m_{1}\\
-m_{2}&m_{1}&0
\end{array}\right)\in so(3),
$$
is the matrix of the Hamiltonian vector field. It is easy to
verify that it satisfies the Jacobi identity or according to
theorem 18, to the formula :
$$\sum_{k=1}^{3}\left( J_{kj}\frac{\partial J_{li}}{\partial m_{k}}+
J_{ki}\frac{\partial J_{jl}}{\partial m_{k}}+J_{kl}\frac{\partial
J_{ij}}{\partial m_{k}}\right)=0,\text{ }\forall 1\leq i,j,l\leq
3.$$
\end{Exmp}

\begin{Exmp}
The equations of the geodesic flow on the group $SO(4)$ can be
written in the form  :
\begin{eqnarray}
\dot{x}_{1}&=&\left( \lambda _{3}-\lambda _{2}\right)
x_{2}x_{3}+\left( \lambda _{6}-\lambda _{5}\right)
x_{5}x_{6},\nonumber\\
\dot{x}_{2}&=&\left( \lambda _{1}-\lambda _{3}\right) x_{1}x_{3}
+\left( \lambda _{4}-\lambda _{6}\right) x_{4}x_{6},\nonumber\\
\dot{x}_{3}&=&\left( \lambda _{2}-\lambda _{1}\right) x_{1}x_{2}
+\left( \lambda _{5}-\lambda _{4}\right) x_{4}x_{5},\\
\dot{x}_{4}&=&\left( \lambda _{3}-\lambda _{5}\right) x_{3}x_{5}
+\left( \lambda _{6}-\lambda _{2}\right) x_{2}x_{6},\nonumber\\
\dot{x}_{5}&=&\left( \lambda _{4}-\lambda _{3}\right) x_{3}x_{4}
+\left( \lambda _{1}-\lambda _{6}\right) x_{1}x_{6},\nonumber\\
\dot{x}_{6}&=&\left( \lambda _{2}-\lambda _{4}\right)
x_{2}x_{4}+\left( \lambda _{5}-\lambda _{1}\right)
x_{1}x_{5},\nonumber
\end{eqnarray}
where $\lambda_1,...,\lambda_6$ are constants. These equations can
be written in the form of a Hamiltonian vector field. We have
$$\dot{x}(t)=J\frac{\partial H}{\partial x}, \quad x\in
\mathbb{R}^{6},$$ with
$$
H=\frac{1}{2}\left(\lambda _{1}x_{1}^{2}+\lambda
_{2}x_{2}^{2}+\cdots +\lambda _{6}x_{6}^{2}\right).$$ By
proceeding in a similar way to the previous example, we obtain
$$
J=\left(\begin{array}{cccccc}
0&-x_{3}&x_{2}&0&-x_{6}&x _{5}\\
x_{3}&0&-x_{1}&x_{6}&0&-x _{4}\\
-x_{2}&x_{1}&0&-x _{5}&x _{4}&0\\
0&-x_{6}&x_{5}&0&-x_{3}&x_{2}\\
x_{6}&0&-x_{4}&x_{3}&0&-x_{1}\\
-x _{5}&x_{4}&0&-x_{2}&x_{1}&0
\end{array}\right).
$$
\end{Exmp}

\begin{Exmp}
The movement of the Kowalewski spinning top is governed by the
following equations (see subsection 10.1 for more information) :
\begin{eqnarray}
\dot{m}&=&m\wedge \lambda m+\gamma \wedge l,\nonumber\\
\dot{\gamma}&=&\gamma \wedge \lambda m,\nonumber
\end{eqnarray}
where $m, \gamma$ and $l$ denote respectively the angular
momentum, the direction cosine of the $z$ axis (fixed in space),
the center of gravity which can be reduced to
$l=\left(1,0,0\right)$ and $\lambda
m=\displaystyle{\left(\frac{m_{1}}{2},
\frac{m_{2}}{2},\frac{m_{3}}{2}\right)}$. These equations can be
written in the form of a Hamiltonian vector field. The system
above is written in the form of a Hamiltonian vector field
$$\dot{x}=J\frac{\partial H}{\partial x},
\quad
x=(m_{1},m_{2},m_{3},\gamma_{1},\gamma_{2},\gamma_{3})^\top,$$
with
$$H=\frac{1}{2}\left( m_{1}^{2}+m_{2}^{2}\right) +m_{3}^{2}+2\gamma
_{1},$$ the Hamiltonian and
$$
J=\left(\begin{array}{cccccc}
0&-m_{3}&m_{2}&0&-\gamma_{3}&\gamma _{2}\\
m_{3}&0&-m_{1}&\gamma_{3}&0&-\gamma _{1}\\
-m_{2}&m_{1}&0&-\gamma _{2}&\gamma _{1}&0\\
0&-\gamma_{3}&\gamma_{2}&0&0&0\\
\gamma_{3}&0&-\gamma_{1}&0&0&0\\
-\gamma _{2}&\gamma_{1}&0&0&0&0
\end{array}\right).
$$
\end{Exmp}

\begin{Exmp}
The motion of a solid in a perfect fluid is described using the
Kirchhoff equations :
\begin{equation}\label{eqn}
\dot{p}=p\wedge \frac{\partial H}{\partial l},\qquad
\dot{l}=p\wedge \frac{\partial H}{\partial p}+l\wedge
\frac{\partial H}{\partial l},
\end{equation}
where $p=(p_{1},p_{2},p_{3})\in \mathbb{R}^{3}$,
$l=(l_{1},l_{2},l_{3})\in \mathbb{R}^{3}$ and $H$ the Hamiltonian.
The problem of this movement is a limit case of the geodesic flow
on $SO(4)$. In the case of Clebsch, we have
$$H=\frac{1}{2}\sum_{k=1}^{3}\left(
a_{k}p_{k}^{2}+b_{k}l_{k}^{2}\right),$$ with the condition :
$$\frac{a_{2}-a_{3}}{b_{1}}+\frac{a_{3}-a_{1}}{b_{2}}+\frac{a_{1}
-a_{2}}{b_{3}}=0.$$ The system (19) is written in the form of a
Hamiltonian vector field : $$\dot{x}= J\frac{\partial H}{\partial
x}, \quad x=(p_{1},p_{2},p_{3},l_{1},l_{2},l_{3})^\top,$$ where
$$
J=\left(\begin{array}{cccccc}
0&0&0&0&-p_{3}&p _{2}\\
0&0&0&p_{3}&0&-p _{1}\\
0&0&0&-p _{2}&p _{1}&0\\
0&-p_{3}&p_{2}&0&-l_{3}&l_{2}\\
p_{3}&0&-p_{1}&l_{3}&0&-l_{1}\\
-p _{2}&p_{1}&0&-l_{2}&l_{1}&0
\end{array}\right).
$$
\end{Exmp}

\begin{Exmp}
a) Let
$$\frac{df}{dt}=\sum_{k=1}^n\left(\frac{\partial f}{\partial p_k}\dot{p}_k+
\frac{\partial f}{\partial q_k}\dot{q}_k\right)+\frac{\partial
f}{\partial t},$$ be the total derivative of a function $f(p, q,
t)$ with respect to $t$. We will determine a necessary and
sufficient condition for $f$ to be a first integral of a system
described by a Hamiltonian $H$. Taking into account Hamilton's
equations, we obtain the expression
$$\frac{df}{dt}=\{f,H\}+\frac{\partial f}{\partial t}.$$ We deduce
that $f$ is a first integral of a system described by a
Hamiltonian $H(p, q, t)$ explicitly dependent on $t$ if and only
if
\begin{equation}\label{eqn:euler}
\{f,H\}+\frac{\partial f}{\partial t}=0,
\end{equation}
and obviously if $f$ does not depend explicitly on $t$, we have
$\{f,H\}=0$.
\end{Exmp}

\begin{Exmp}
Consider a Hamiltonian $$H=\frac{1}{2m}(p_1^2+p_2^2+p_3^2)+V(r,
t), \quad r=\sqrt{q_1^2+q_2^2+q_3^2},$$ describing the motion of a
particle having a mass $m$ and immersed into a potential $V(r,
t)$. We will determine three first integrals of the system
described by this Hamiltonian. The two components of kinetic
moment are equal to $$H_1=q_2p_3-q_3p_2,\qquad
H_2=q_3p_1-q_1p_3.$$ They are obviously first integrals. According
to Poisson's theorem 16, we have
$$\{H_1,H_2\}=q_1p_2-q_2p_1=H_3,$$ which shows that $H_3$ is also
a first integral. Note also that : $$\{H_3,H_1\}=H_2, \qquad
\{H_2,H_3\}=H_1.$$ If in a system two components of kinetic moment
are first integrals, then the third component is also a first
integral.
\end{Exmp}

\begin{Exmp}
We have already seen that in a conservative system, the
Hamiltonian $H(p, q)$ is a first integral. We will show that if
$F(p, q, t)$ denotes another first integral explicitly dependent
on $t$, then $\frac{\partial^k F}{\partial t^k}$ is also a first
integral. We will apply this result to the case of the Hamiltonian
of the harmonic oscillator:
$$H=\frac{1}{2m}p^2+\frac{m\omega^2}{2}q^2.$$ According to the
Poisson theorem 10, $\{F, H\}$ is also a first integral.
Therefore, $$\frac{\partial F} {\partial t}=-\{F, H\},$$ is a
first integral under (20). Similarly, we have
$$\left\{\frac{\partial F}{\partial t},H\right\}+\frac{\partial^2
F}{\partial t^2}=0,$$ which shows that $$\frac{\partial^2
F}{\partial t^2}=-\left\{\frac{\partial F}{\partial
t},H\right\},$$ is also a first integral. And similarly, we show
that $\frac{\partial^k F}{\partial t^k}$ is a first integral. For
the Hamiltonian of the harmonic oscillator, we easily check that
$$F=q\cos\omega t-\frac{1}{m\omega}p\sin\omega t,$$ and
$$\frac{\partial F}{\partial t}=-\omega q\sin\omega
t-\frac{1}{m\omega}p\cos\omega t,$$ are first integrals of the
Hamiltonian system associated with $H$.
\end{Exmp}

\section{Coadjoint orbits and their symplectic structures}

We will first define the adjoint and coadjoint orbits of a Lie
group. Let $G$ be a Lie group and $g$ an element of $G$. The Lie
group $G$ operates on itself by left translation : $$L_g
:G\longrightarrow G, \quad h\longmapsto gh,$$ and by right
translation : $$R_g : G\longrightarrow G, \quad h\longmapsto hg.$$
By virtue of the associative law of the group, we have
$$L_gL_h=L_{gh}, \quad R_gL_h=R_{hg}, \quad L_{g^{-1}}=L_g^{-1},
\quad R_{g^{-1}}=R_g^{-1}.$$ In particular, the applications $R_g$
and $L_g$  are diffeomorphisms of $G$. Also, because of
associativity, $R_g$ and $L_g$ commute. Consider $$R_g^{-1}L_g :
G\longrightarrow G, \quad h\longmapsto ghg^{-1},$$ the
automorphism of the group $G$. It leaves the unit $e$ of the group
$G$ fixed, i.e., $$R_g^{-1}L_g(e)=geg^{-1}=e.$$ We can define the
adjoint representation of the group $G$ as the derivative of
$R_g^{-1}L_g$ in the unit $e$, that is, the induced application of
tangent spaces as follows
$$Ad_g : \mathcal{G}\longrightarrow \mathcal{G},\quad \xi
\longmapsto
\left.\frac{d}{dt}R_g^{-1}L_g(e^{t\xi})\right|_{t=0},$$ where
$\mathcal{G}=T_eG$ is the Lie algebra of the $G$ group ; it is the
tangent space at $G$ in its unit $e$. This definition has a
meaning because $R_g^{-1}L_g(e^{t\xi})$ is a curve in $G$ and
passes through the identity in $t=0$. Therefore, $g\xi
g^{-1}\in\mathcal{G}$.

\begin{Theo}
For any element $\xi\in \mathcal{G}$, we have $$Ad_g(\xi)=g\xi
g^{-1}, \quad g\in G,$$ and $$Ad_{gh}=Ad_g.Ad_h.$$ The application
$Ad_g$ is an algebra homomorphism, i.e.,
$$Ad_g[\xi,\eta]=[Ad_g\xi,Ad_g\eta], \quad(\xi, \eta\in
\mathcal{G}).$$
\end{Theo}
\emph{Proof}. We have
$$
Ad_g(\xi)=\left.\frac{d}{dt}R_g^{-1}L_g(e^{t\xi})\right|_{t=0}
=\left.\frac{d}{dt}ge^{t\xi}g^{-1}\right|_{t=0},
$$
hence
\begin{eqnarray}
Ad_g(\xi)&=&\left.\frac{d}{dt}g\left(\sum_{n=0}^\infty
\frac{t^n\xi^n}{n!}\right)g^{-1}\right|_{t=0},\nonumber\\
&=&\left.\frac{d}{dt}\sum_{n=0}^\infty\frac{t^n}{n!}g\xi^ng^{-1}\right|_{t=0},\nonumber\\
&=&\left.\frac{d}{dt}\sum_{n=0}^\infty\frac{t^n}{n!}\underset{n-times}{\underbrace{g\xi
g^{-1}.g\xi g^{-1}...g\xi g^{-1}}}\right|_{t=0},\nonumber\\
&=&\left.\frac{d}{dt}\sum_{n=0}^\infty\frac{t^n}{n!}(g\xi
g^{-1})^n\right|_{t=0},\nonumber
\end{eqnarray}
and finally
$$
Ad_g(\xi)=\left.\frac{d}{dt}e^{t(g\xi g^{-1})}\right|_{t=0}=g\xi
g^{-1}.
$$
We easily check that : $$Ad_{gh}=Ad_g.Ad_h.$$ Indeed, we have
\begin{eqnarray}
Ad_{gh}(\xi)&=&gh\xi(gh)^{-1}=gh\xi h^{-1}g^{-1},\nonumber\\
Ad_g.Ad_h(\xi)&=&Ad_g(h\xi h^{-1})=gh\xi h^{-1}g^{-1}.\nonumber
\end{eqnarray}
We have
\begin{eqnarray}
Ad_g[\xi,\eta]&=&Ad_g(\xi \eta-\eta \xi),\nonumber\\
&=&g(\xi \eta-\eta\xi)g^{-1},\nonumber\\
&=&g\xi \eta g^{-1}-g\eta \xi g^{-1},\nonumber\\
&=&g\xi g^{-1}g\eta g^{-1}-g\eta g^{-1}g\xi g^{-1},\nonumber\\
&=&[g\xi g^{-1},g\eta g^{-1}],\nonumber\\
&=&[Ad_g\xi, Ad_g\eta],\nonumber
\end{eqnarray}
which completes the demonstration. $\square$

The adjoint orbit of $\xi$ is defined by
$$\mathcal{O}_G(\xi)=\{Ad_g(\xi) : g\in G\}\subset \mathcal{G}.$$
Now consider the function
$$Ad : G\longrightarrow \mbox{End} (\mathcal{G}),\quad
g\longmapsto Ad(g)\equiv Ad_g,$$ where $\mbox{End} (\mathcal{G})$
is the space of the linear operators on the algebra $\mathcal{G}$.
The application $Ad$ is differentiable and its derivative
$Ad_{*e}$ in the unit of the group $G$ is a linear map from the
algebra $T_eG=\mathcal{G}$ to the vector space $T_I\mbox{End}
(\mathcal{G})=\mbox{End} (\mathcal{G})$. This application will be
noted
$$ad\equiv Ad_{*e} : \mathcal{G} \longrightarrow \mbox{End}
(\mathcal{G}),\quad \xi \longmapsto
ad_\xi=\left.\frac{d}{dt}Ad_{g(t)}\right|_{t=0},$$ where $g(t)$ is
a one-parameter group with
$\displaystyle{\left.\frac{d}{dt}g(t)\right|_{t=0}=\xi}$ and
$g(0)=e$.

\begin{Theo}
Let $\xi\in \mathcal{G}$ et $\eta\in \mbox{End} (\mathcal{G})$. By
setting $ad_\xi\equiv Ad_{*e}(\xi)$, then
$$ad_\xi(\eta)=[\xi,\eta].$$
\end{Theo}
\emph{Proof}. We have
\begin{eqnarray}
ad_\xi(\eta)&=&Ad_{*e}(\xi)(\eta),\nonumber\\
&=&\left.\frac{d}{dt}Ad_{g(t)}(\eta)\right|_{t=0},\nonumber\\
&=&\left.\frac{d}{dt}(g(t)\eta g^{-1}(t))\right|_{t=0},\nonumber\\
&=&\left.\dot{g}(t)\eta g^{-1}(t)\right|_{t=0}-\left.g(t)\eta
g^{-1}(t)\dot{g}(t)g^{-1}(t)\right|_{t=0},\nonumber\\
&=&\dot{g}(0)\eta-\eta\dot{g}(0),\nonumber\\
&=&\xi\eta-\eta\xi,\nonumber\\
&=&[\xi,\eta],\nonumber
\end{eqnarray}
which completes the proof. $\square$

Let $T^*_gG$ be the cotangent space to the group $G$ at $g$; it is
the dual to the tangent space $T_gG$. Then an element $\zeta \in
T^*_gG$ is a linear form on $T_gG$ and its value on $\eta \in
T_gG$ will be denoted by, $$\zeta(\eta)\equiv
\langle\zeta,\eta\rangle.$$ Let $\mathcal{G}^*=T^*_eG$ be the dual
vector space to the Lie algebra $\mathcal{G}$; it is the cotangent
space to the group $G$ in its unit $e$. The transpose operators
$Ad^*_g : \mathcal{G}^*\longrightarrow \mathcal{G}^*$, where $g$
runs through the Lie group $G$ are defined by $$\langle
Ad^*_g(\zeta),\eta\rangle=\langle\zeta,Ad_g\eta\rangle,
\quad\zeta\in \mathcal{G}^*, \quad\eta\in\mathcal{G}.$$ $Ad^*_g$
is called coadjoint representation of the Lie group $G$. The
coadjoint orbit (also called Kostant-Kirillov orbit) is defined at
the point $x\in \mathcal{G}^*$ by $$\mathcal{O}^*_G(x)=\{Ad^*_g(x)
: g\in G\}\subset \mathcal{G}^*.$$

\begin{Theo}
The transpose operators $Ad^*_g$ form a representation of the Lie
group $G$, i.e., they satisfy the relations :
$Ad^*_{gh}=Ad^*_{h}.Ad^*_{g}$.
\end{Theo}
\emph{Proof}. Indeed, let $\zeta\in \mathcal{G}^*$, $\eta\in
\mathcal{G}$. We have
$$
\langle
Ad^*_{gh}(\zeta),\eta\rangle=\langle\zeta,Ad_{gh}(\eta)\rangle
=\langle\zeta,Ad_{h}.Ad_{g}(\eta)\rangle,$$ hence,
$$
\langle Ad^*_{gh}(\zeta),\eta\rangle=\langle Ad^*_{g}(\zeta),
Ad_{h}(\eta)\rangle=\langle Ad^*_{h}.Ad^*_{g}(\zeta),
\eta\rangle,$$ which completes the demonstration. $\square$

Consider the map
$$Ad^* : G\longrightarrow \mbox{End}(\mathcal{G}^*),\quad
g\longmapsto Ad^*(g)\equiv Ad^*_g,$$ and its derivative in the
unity of the group
$$ad^*\equiv (Ad^*)_{*e} : \mathcal{G}\longrightarrow \mbox{End}(\mathcal{G}^*),\quad
\xi\longmapsto ad^*_\xi.$$

\begin{Theo}
By setting $$\langle ad^*_\xi(\zeta),\eta\rangle=\langle
\zeta,[\xi,\eta]\rangle=\langle\{\xi,\zeta\},\eta\rangle,$$ where
$$\{,\}: \mathcal{G}\times \mathcal{G}^*\longrightarrow
\mathcal{G}^*, \quad(\xi,\zeta)\longmapsto \{\xi,\zeta\},
\quad(\xi ,\eta \in \mathcal{G}, \zeta \in \mathcal{G}^*),$$ then
$$ad^*_\xi(\zeta)=\{\xi,\zeta\}.$$
\end{Theo}
\emph{Proof}. We have $$\langle
ad^*_{\xi}(\zeta),\eta\rangle=\langle
(Ad^*)_{*e}(\zeta),\eta\rangle=\left\langle
\left.\frac{d}{dt}Ad^*_{e^{t\xi}}(\zeta)\right|_{t=0},\eta\right\rangle,$$
with $\left.e^{t\xi}\right|_{t=0}=e$ and
$\left.\frac{d}{dt}e^{t\xi}\right|_{t=0}=\xi$. Hence,
\begin{eqnarray}
\langle ad^*_{\xi}(\zeta),\eta\rangle&=&\left.\frac{d}{dt}\langle
Ad^*_{e^{t\xi}}(\zeta),\eta\rangle\right|_{t=0},\nonumber\\
&=&\left.\frac{d}{dt}\langle\zeta,Ad_{e^{t\xi}}(\eta)\rangle\right|_{t=0},\nonumber\\
&=&\left\langle\zeta,\left.\frac{d}{dt}Ad_{e^{t\xi}}(\eta)\right|_{t=0}\right\rangle,\nonumber\\
&=&\langle\zeta,ad_\xi\eta\rangle,\nonumber\\
&=&\langle\zeta,[\xi,\eta]\rangle,\nonumber\\
&=&\langle\{\xi,\zeta\},\eta\rangle,\nonumber
\end{eqnarray}
which completes the proof. $\square$

We will show below, how to find the adjoint orbit and the
coadjoint orbit in the case of the group $SO(n)$. Recall that
$SO(n)$ is the special orthogonal group\index{special orthogonal
group} of order $n$, that is, the set of matrices $X$ of order
$n\times n$ such that : $X^\top.X=I$ (or $X^{-1}=X^\top$) and
$\det X=1$. $SO(n)$ is a Lie group. The tangent space to the
identity of the group $SO(n)$, which is denoted $so(n)$, consists
of the antisymmetric matrices of order $n\times n$, i.e., that is,
matrices $A$ such that : the commutator of two antisymmetric
matrices is still an antisymmetric matrix (if $A, B\in so(n)$,
then $[A, B]=AB-BA\in so(n)$). This product defines a Lie algebra
structure on $so(n)$; it is the Lie algebra of the group $SO(n)$.
In addition, we have $\dot{X}=AX$ with $A\in so(n)$ and therefore
the tangent space to the identity of $SO(n)$ is $T_ISO(n)=so(n)$.
Let $$R^{-1}_YL_Y : SO(n)\longrightarrow SO(n), \quad X\longmapsto
YXY^{-1}, \quad Y\in SO(n),$$ be the automorphism interior of the
group $SO(n)$. When looking for the coadjoint orbit, we have to
use the following obvious lemma :

\begin{Lem}
The Lie algebra $so(n)$ with the commutator $[,]$ of matrix is
isomorphic to the space $\mathbb{R}^{\frac{n(n-1)}{2}}$ with the
vector product $\wedge$. The isomorphism is given by $$a\wedge
b\longmapsto [A,B]=AB-BA,$$ where $a,b\in
\mathbb{R}^{\frac{n(n-1)}{2}}$ and $A,B \in so(n)$.
\end{Lem}

\begin{Theo}
The orbit of the adjoint representation of the group $SO(n)$ is
$$\mathcal{O}_{SO(n)}(A)=\{YAY^{-1}: Y\in SO(n)\}, \quad A\in
so(n).$$ Let $A\in so(n)$. The coadjoint orbit of the group
$SO(n)$ is
\begin{eqnarray}
\mathcal{O}^*_{SO(n)}(A)&=&\{Y^{-1}AY: Y\in SO(n)\},\nonumber\\
&=&\{C\in so(n): C=Y^{-1}AY, \mbox{spectrum of } C=\mbox{spectrum
of } A\}.\nonumber
\end{eqnarray}
With the notation of theorem 22, we have $$\{A,B\}=[B,A],
\quad(A,B\in so(n)).$$
\end{Theo}
\emph{Proof}. Let $Y\in SO(n)$, $A\in so(n)$. By definition, the
adjoint representation of the group $SO(n)$ is $$Ad_Y :
so(n)\longrightarrow so(n), \quad A\longmapsto YAY^{-1}.$$ We have
$$(YAY^{-1})^\top=(Y^{-1})^\top A^\top
Y^\top=-YAY^\top=-YAY^{-1}.$$ So $YAY^{-1}\in so(n)$. Let $$Ad :
SO(n)\longrightarrow \mbox{End} (so(n)), \quad Y\longmapsto
Ad_Y,$$ where $Ad_Y(A)=YAY^{-1}$, $A \in so(n)$, and let $$ad :
so(n)\longrightarrow \mbox{End} (so(n)), \quad
\dot{Y}(0)\longmapsto ad_{\dot{Y}(0)},$$ with
$$ad_{\dot{Y}(0)}\bullet =[\dot{Y}(0),\bullet] : so(n)\longrightarrow
so(n),\quad A\longmapsto [\dot{Y}(0),A],$$ where $Y(t)$ is a curve
in $SO(n)$ with $Y(0)=I$. Since $(R^{n\times n})^*\simeq
R^{n\times n}$, then according to the previous lemma, we also have
the isomorphism $(so(n))^*\simeq so(n)$. We can therefore define
$Ad^*$ by $Ad_Y^* : so(n)\longrightarrow so(n)$, with
$$
\langle Ad_Y^*(A),B\rangle=\langle A,Ad_YB\rangle=\langle
A,YBY^{-1}\rangle,\quad (A,B\in so(n)),$$ i.e.,
$$\langle Ad_Y^*(A),B\rangle=-\frac{1}{2}\mbox{tr} (AYBY^{-1})
=-\frac{1}{2}\mbox{tr} (Y^{-1}AYB) =\langle Y^{-1}AY,B\rangle,$$
hence, $$Ad^*_Y(A)=Y^{-1}AY.$$ We easily check that $Y^{- 1}AY\in
so(n)$. Indeed, we have $$(Y^{-1}AY)^\top =Y^\top A^\top
(Y^{-1})^\top =-Y^{-1}AY,$$ because $Y \in SO(n)$ and $A\in
so(n)$. Then $$\mathcal{O}^*_{SO(n)}(A)=\{Y^{-1}AY: Y\in
SO(n)\},$$ that we can write in the form
$$\mathcal{O}^*_{SO(n)}(A)=\{C\in so(n) : \exists Y \in SO(n),
C=Y^{-1}AY\}.$$ Note that $\det(C-\lambda I)=\det(A-\lambda I)$.
Then the matrices $C$ and $A$ have the same characteristic
polynomial, and consequently they have the same spectrum.
$$
\mathcal{O}^*_{SO(n)}(A)=\{C\in so(n): C=Y^{-1}AY, \mbox{spectrum
of } C=\mbox{spectrum of } A\}.
$$
Now apply theorem 16 to the case of the group $SO(n)$. Let's go
back to the linear form knowing that $(so(n))^*=so(n)$,
$$\{,\} : so(n)\times so(n) \longrightarrow so(n),\quad
(A,B)\longmapsto\{A,B\},$$ as well as the applications
$$Ad^* : SO(n)\longrightarrow \mbox{End} (so(n)),\quad
Y\longmapsto Ad^*_Y(B)=Y^{-1}BY,\quad B\in so(n),$$
$$ad^* : so(n)\longrightarrow \mbox{End} (so(n)),\quad
A\longmapsto ad^*_A,$$ where $$\langle ad^*_A(B),C\rangle=\langle
B,[A,C]\rangle=\langle \{A,B\},C\rangle.$$ We have
$$
\langle \{A,B\},C\rangle=\langle B,[A,C]\rangle
=-\frac{1}{2}\mbox{tr} (B.[A,C])=-\frac{1}{2}\mbox{tr}
(BAC-BCA),$$ hence, $$\langle
\{A,B\},C\rangle=-\frac{1}{2}\mbox{tr} ([B,A].C)=\langle
[B,A],C\rangle.$$ Then $\{A, B\}=[B, A]$, and the theorem is
proved. $\square$

We will see how to define a symplectic structure on the coadjoint
orbit with an application in the case of the groups $SO(3)$ and
$SO(4)$. Let $x\in\mathcal{G}^*$, $\xi$ the tangent vector in $x$
to the orbit. Since $\mathcal{G}^*$ is a vector space, then
obviously $\xi\in T_x\mathcal{G}^*=\mathcal{G}^*$. let's remember
that $$\mathcal{O}^*_G(x)=\{Ad^*_g(x) : g\in G\}\subset
\mathcal{G}^*.$$ For $x\in\mathcal{O}^*_G(x)$, there exists $g\in
G$ such that : $x=Ad^*_g$. Let $a\in \mathcal{G}$ and $e^{ta}$ be
a group with a parameter in $G$ with $\left.e^{ta}\right|_{t=0}=g$
and $\left.\frac{d}{dt}Ad^*_{e^{ta}}(x)\right|_{t=0}=\xi$. Since
$$\left.\frac{d}{dt}Ad^*_{e^{ta}}(x)\right|_{t=0}\equiv
ad^*_a(x)=\{a,x\},$$ therefore the vector $\xi$ can be represented
as the velocity vector of the motion of $x$ under the action of a
group $e^{ta}$, $a\in \mathcal{G}$. In other words, any vector
$\xi$ tangent to the orbit $\mathcal{O}^*_G(x)$ is expressed as a
function of $a\in\mathcal{G}$ by
\begin{equation}\label{eqn:euler}
\xi=\{a,x\},\quad a\in \mathcal{G},\quad x\in \mathcal{G}^*.
\end{equation}
Therefore, we can determine the value of a $2$-form $\Omega$ on
the orbit $\mathcal{O}^*_G(x)$ as follows : let $\xi_1$ and
$\xi_2$ be two vectors tangent to the orbit of $x$. From the
above, we have $$\xi_1=\{a_1,x\}, \quad \xi_2=\{a_2,x\},
\quad(a_1, a_2\in \mathcal{G}), \quad x\in \mathcal{G}^*.$$ We can
easily verify that the differential $2$-form
\begin{equation}\label{eqn:euler}
\Omega(\xi_1,\xi_2)(x)=\langle x,[a_1,a_2]\rangle,\quad a_1,
a_2\in \mathcal{G},\quad x\in \mathcal{G}^*,
\end{equation}
on $\mathcal{O}^*_G(x)$ is well defined; its value does not depend
on the choice of $a_1, a_2$. It is antisymmetric, non-degenerate
and closed. To determine the symplectic structure on
$\mathcal{O}^*_{SO(3)}(X)$, we proceed as follows: according to
(22), we have $$\Omega(\xi_1,\xi_2)(X)=\langle X,[A,B]\rangle,$$
where $A,B\in so(3)$, $X\in (so(3))^*=so(3)$ and according to
(21), $$\xi_1=\{A,X\}, \quad \xi_2=\{B,X\},$$ are two tangent
vectors to the orbit in $X$ or what the same according to theorem
23, $\xi_1=[X,A]$, $\xi_2=[X,B]$. Using the isomorphism between
$(so(3),[,])$ and $(\mathbb{R}^3,\wedge)$, we also have
$\xi_1=x\wedge a$, $\xi_2=x\wedge b$, with
$$\Omega(\xi_1,\xi_2)(x)=\langle x,a\wedge b\rangle.$$ According to
theorem 23, the coadjoint orbit of $SO(3)$ is
$$\mathcal{O}^*_{SO(3)}(A)=\{C\in so(3) : C=Y^{-1}AY, \mbox{ spectrum of } C=\mbox{spectrum of }
A\},$$ where $A\in so(3)$ et $Y\in SO(3)$. Let's determine the
spectrum of the matrix
$${A}=\left(\begin{array}{ccc}
0&-a_{3}&a_{2}\\
a_{3}&0&-a_{1}\\
-a_{2}&a_{1}&0
\end{array}\right)\in so(3).$$
We have $$\det (A-\lambda
I)=-\lambda(\lambda^2+a_1^2+a_2^2+a_3^2)=0,$$ hence, $\lambda=0$
and $\lambda=\pm i\sqrt{a_1^2+a_2^2+a_3^2}$. Then
$$\mathcal{O}^*_{SO(3)}(A)=\{C\in so(3) :
c_1^2+c_2^2+c_3^2=r^2\},$$ with
$${C}=\left(\begin{array}{ccc}
0&-c_{3}&c_{2}\\
c_{3}&0&-c_{1}\\
-c_{2}&c_{1}&0
\end{array}\right)\in so(3),$$
and $r^2=a_1^2+a_2^2+a_3^2$. Since the algebra $so(3)$ is
isomorphic to $\mathbb{R}^3$, we deduce that the orbit
$\mathcal{O}^*_{SO (3)}(A)$ is isomorphic to a sphere $S^2$ of
radius $r$. Like vectors $\xi_1$, $\xi_2$ belong to the tangent
plane $T_X\mathcal{O}^*_{SO (3)}$ to $X$, they also belong to the
tangent plane $T_xS^2$ in $x$. Let
$$S^2=\left\{\left(y_1,y_2,y_3\right)\in \mathbb{R}^3 :
y_1^2+y_2^2+y_3^2=r^2\right\},$$ be the sphere of radius $r$, then
the plane tangent to this sphere in $x$ of coordinates
$(x_1,x_2,x_3)$ is
\begin{eqnarray}
T_xS^2&=&\left\{\left(y_1,y_2,y_3\right)\in \mathbb{R}^3 :
y_1x_1+y_2x_2+y_3x_3=0\right\},\nonumber\\
&=&\left\{\left(y_1,y_2,-\frac{y_1x_1+y_2x_2}{x_3}\right)\right\}.
\end{eqnarray}
Let $z=(z_1,z_2,z_3)\in T_xS^2$ and determine $a=(a_1,a_2,a_3)$
such that : $x\wedge a=z$. The latter is equivalent to the system
$$\left(\begin{array}{ccc}
0&-a_{3}&a_{2}\\
a_{3}&0&-a_{1}\\
-a_{2}&a_{1}&0
\end{array}\right)\left(\begin{array}{c}
a_{1}\\
a_{2}\\
a_{3}
\end{array}\right)=
\left(\begin{array}{c}
z_{1}\\
z_{2}\\
-\frac{z_1x_1+z_2x_2}{x_3}
\end{array}\right),
$$
whose solution is
$$a=\left(\frac{x_1a_3+z_2}{x_3},\frac{x_2a_3-z_1}{x_3},a_3\right),\quad a_3
\in\mathbb{R}.$$ Since the symplectic form on $S^2$ that one wants
to determine is intrinsic, i.e., does not depend on the choice of
local coordinates, one can choose as local coordinates $x_1$,
$x_2$ and the same reasoning will be valid for the other cases,
i.e., $x_2$, $x_3$ and $x_3$, $x_1$. So we will calculate $a$ and
$b$ relative to the basis $\left(\frac{\partial}{\partial
x_1},\frac{\partial}{\partial x_2}\right)$ of $T_xS^2$ with
$$\frac{\partial}{\partial
x_1}=\left(1,0,-\frac{x_1}{x_3}\right),\qquad
\frac{\partial}{\partial x_2}=\left(0,1,-\frac{x_2}{x_3}\right).$$
We have
$$a=(a_1,a_2,a_3)=\left(\frac{x_1b_3+1}{x_3},\frac{x_2b_3}{x_3},b_3\right),$$
and
$$
a\wedge b=\left(a_2b_3-a_3b_2,a_3b_1-a_1b_3,a_1b_2-a_2b_1\right)
=\left(-\frac{b_3}{x_3},\frac{a_3}{x_3},\frac{x_1b_3-x_2a_3+1}{x_3^2}\right).
$$
Therefore, $$\Omega \left(\frac{\partial}{\partial
x_1},\frac{\partial}{\partial x_2}\right)=\left(x,a\wedge
b\right)=\frac{1}{x_3},$$ consequently $$\Omega=\frac{dx_1\wedge
dx_2}{x_3}.$$ The symplectic form being intrinsic, we will finally
have
$$\Omega=\frac{dx_1\wedge dx_2}{x_3}=\frac{dx_2\wedge dx_3}{x_1}
=\frac{dx_3\wedge dx_1}{x_2}.$$

\begin{Exmp}
The symplectic structure obtained here is equivalent to that
associated with the system (17). Indeed, we know that
$$J^{-1}=\left(\omega\left(\frac{\partial}{\partial x_{i}},
\frac{\partial}{\partial x_{j}}\right)\right)_{i,j=1,2},$$ so the
matrix associated with the form
$$\Omega=\frac{dx_1\wedge dx_2}{x_3},$$ is
$\left(\begin{array}{cc}
0&-x_{3}\\
x_{3}&0
\end{array}\right)$. Let's show that there is equivalence between
$$
\dot{x}(t)=J\frac{\partial H}{\partial x},\quad
\mbox{where}\quad \left\{\begin{array}{rl}x=&(m_1,m_2,m_3)^\top,\nonumber\\
H=&\frac{1}{2}\left( \lambda _{1}m_{1}^{2}+\lambda_{2}m_{2}^{2}+\lambda_{3}m_{3}^{2}\right),\nonumber\\
J=&\left(\begin{array}{ccc}
0&-m_{3}&m_{2}\\
m_{3}&0&-m_{1}\\
-m_{2}&m_{1}&0
\end{array}\right),
\end{array}\right.
$$
and
$$
\dot{x}(t)=\mathbf{J}\frac{\partial \mathbf{H}}{\partial x},\quad
\mbox{where}\quad
\left\{\begin{array}{rl} x=&(m_1,m_2,m_3)^\top,\nonumber\\
\mathbf{H}=&H(m_1,m_2,m_3),\nonumber\\
\mathbf{J}=&\left(\begin{array}{cc}
0&-m_{3}\\
m_{3}&0
\end{array}\right).
\end{array}\right.
$$
Indeed, we have
$$
\dot{m}_1=-m_3\frac{\partial\textbf{H}}{\partial
m_2}=-m_3\left(\frac{\partial H}{\partial m_2}+\frac{\partial
H}{\partial m_3}\frac{\partial m_3}{\partial m_2}\right),$$ and
$$
\dot{m}_2=m_3\frac{\partial\textbf{H}}{\partial
m_1}=m_3\left(\frac{\partial H}{\partial m_1}+\frac{\partial
H}{\partial m_3}\frac{\partial m_3}{\partial m_1}\right).$$
According to example 10, we have
$$dm_3=-\frac{m_1dm_1+m_2dm_2}{m_3},$$ hence
$$\frac{dm_3}{dm_2}=-\frac{m_2}{m_3},
\qquad\frac{dm_3}{dm_1}=-\frac{m_1}{m_3}.$$
Therefore, we have
\begin{eqnarray}
\dot{m}_{1}&=&\left(\lambda_{3}-\lambda_{2}\right)m_{2}m_{3},\nonumber\\
\dot{m}_{2}&=&\left(\lambda_{1}-\lambda_{3}\right)m_{1}m_{3},\nonumber
\end{eqnarray}
and the result follows.
\end{Exmp}

\begin{Exmp}
To determine the symplectic structure on the coadjoint orbit of
the Lie group $SO(4)$, we can follow the same method as in the
previous case but the calculation is longer. On the other hand,
one can easily obtain the result by using a geometric approach by
observing that $so(4)$ breaks down into two copies of $so(3)$ and
that the generic orbits are a product of two spheres. More
precisely, from $SO(4)=SO(3)\otimes SO(3)$, it is more interesting
to consider the coordinates $(x_1,x_2,x_3)$, $(x_4,x_5,x_6)$ with
$$(x_1,x_2,x_3)\oplus (x_4,x_5,x_6)\in so(4)\simeq so(3)\oplus
so(3).$$ We obtain
$$\Omega=-x_{3}dx_{1}\wedge dx_{2}-x_{6}dx_{1}\wedge dx_{5}
+x_{6}dx_{2}\wedge dx_{4}-x_{3}dx_{4}\wedge dx_{5}.$$
\end{Exmp}

\section{Arnold-Liouville theorem and completely integrable systems}

The so-called Arnold-Liouville theorem [4] play a crucial role in
the study of the integrability of Hamiltonian systems; the regular
compact level manifolds defined by the intersection of the
constants of motion are diffeomorphic to a real torus on which the
motion is quasi-periodic as a consequence of the following purely
differential geometric fact : a compact and connected
$n$-dimensional manifold on which there exist $n$ vector fields
which commute and are independent at every point is diffeomorphic
to an $n$-dimensional real torus and each vector field will define
a linear flow there. Consider the Hamiltonian system (15)
associated with the function $H$ (Hamiltonian) on a
$2n$-dimensional symplectic manifold $M$.

\begin{Theo}
Let $H_{1}=H,H_{2},...,H_{n}$, be $n$ first integrals on a
$2n$-dimensional symplectic manifold $M$ that are functionally
independent, i.e., $$dH_1\wedge ... \wedge dH_n\neq 0,$$ and
pairwise in involution, i.e., $$\left\{H_{i},H_{j}\right\}=0,
\quad 1\leq i,j\leq n.$$ For generic $c=(c_1,...,c_n)\in
\mathbb{R}^n$, the level set
$$M_{c}=\bigcap_{i=1}^{n}\left\{ x\in M:H_{i}\left( x\right)
=c_{i}\right\},$$ will be an $n$-manifold. If $M_{c}$ is compact
and connected, it is diffeomorphic to an $n$-dimensional torus
$$T^{n}=\mathbb{R}^{n}/lattice=\left\{\left(\varphi
_{1},...,\varphi_{n}\right) \text{ mod. }2\pi \right\}.$$ The
flows $g_{t}^{X_{1}}(x)$,...,$g_{t}^{X_{n}}(x)$ defined by the
vector fields $X_{H_{1}}$,...,$X_{H_{n}}$, are straight-line
motions on $T^n$ and determine on $T^{n}$ a quasi-periodic motion,
i.e., in angular coordinates $\varphi_{1},...,\varphi_{n}$, we
have
$$\dot{\varphi}_i=\omega_i(c),\quad\omega_i(c)=\text{constants},\quad\varphi_i(t)=\varphi_i(0)+\omega_it.$$
The equations (15) of the problem are integrable by
quadratures\index{equations integrable by quadratures}.
\end{Theo}
\emph{Proof}. \textbf{1)} Let us first show that a compact and
connected $n$-dimensional manifold $M$ on which there exist $m$
differential (of class $\mathcal{C}^{\infty }$) vector fields
$X_{1},...,X_{m}$ which commute and are independent at every point
is diffeomorphic to an $m$-dimensional real torus :
$$T^{m}=\mathbb{R}^{m}/lattice=\{(\varphi_{1},...,\varphi_{m})
\text{ mod. }2\pi\}.$$ Let us define the application
$$g:\mathbb{R}^{m}\longrightarrow M,
\quad\left(t_{1},...,t_{m}\right)\longmapsto
g\left(t_{1},...,t_{m}\right),$$ where
$$g\left( t_{1},...,t_{m}\right) =g_{t_{1}}^{X_{1}}
\circ \cdots \circ g_{t_{m}}^{X_{m}}\left( x\right) =
g_{t_{m}}^{X_{m}}\circ \cdots \circ g_{t_{1}}^{X_{1}}\left(
x\right),\text{ }x\in M.$$ $a)$ The application $ g $ is a local
diffeomorphism. Indeed, let
$$g_{r}\equiv g\mid_{_{U}}:U\longrightarrow M,\text{ }
\left(t_{1},...,t_{m}\right) \longmapsto
g_{r}\left(t_{1},...,t_{m}\right)=g_{t_{m}}^{X_{m}}\circ \cdots
\circ g_{t_{1}}^{X_{1}}(x),$$ be the restriction of $g$ on a
neighborhood $U$ of $(0,...,0)$ in $\mathbb{R}^{m}$ with
$x=g_{r}(0,...,0)$. Let us show that the map $g_{r}$ is
differentiable (of class $\mathcal{C}^{\infty }$). We have
$$\frac{\partial }{\partial
t_{1}}g_{t_{1}}^{X_{1}}=X_{1}(x)=(\dot{x}_{1},...,\dot{x}_{m}),$$
with
$$
\dot{x}_{1}=f_{1}(x_{1},...,x_{m}),...,\dot{x}_{m}=f_{m}(x_{1},...,x_{m}),$$
where $f_{1},...,f_{m}:M\longrightarrow \mathbb{R}$ are functions
on $M$. Similarly, we have
\begin{eqnarray}
&&\frac{\partial ^{2}}{\partial
t_{1}^{2}}g_{t_{1}}^{X_{1}}=(\ddot{x}_{1},...,
\ddot{x}_{m})=\left( \sum_{k=1}^{m}\frac{\partial f_{1}}{\partial
x_{k}}\dot{x}_{k},...,\sum_{k=1}^{m}\frac{\partial f_{m}}{\partial
x_{k}}\dot{x}_{k}\right),\nonumber\\
&&\frac{\partial ^{3}}{\partial
t_{1}^{3}}g_{t_{1}}^{X_{1}}=(\dddot{x}_{1},...,\dddot{x}_{m}),\nonumber\\
&&=\left( \sum_{k=1}^{m}\sum_{l=1}^{m}\frac{\partial
^{2}f_{1}}{\partial x_{k}\partial
x_{l}}\dot{x}_{k}\dot{x}_{l}+\frac{\partial f_{1}}{\partial
x_{k}}\ddot{x}_{k},...,\sum_{k=1}^{m}\sum_{l=1}^{m}\frac{\partial^{2}f_{m}}{\partial
x_{k}\partial x_{l}}\dot{x}_{k}\dot{x}_{l}+\frac{\partial
f_{m}}{\partial x_{k}}\ddot{x}_{k}\right),\nonumber
\end{eqnarray} etc. All these
expressions have a meaning because by hypothesis all the functions
$f_{1},...,f_{m}$ are $\mathcal{C}^{\infty }$. A similar reasoning
shows that $g_{t_{2}}^{X_{2}},...,g_{t_{m}}^{X_{m}}$ are also
$\mathcal{C}^{\infty }$. Since the composite of functions
$\mathcal{C}^{\infty }$ is $\mathcal{C}^{\infty }$, we deduce that
$g_{r}\left(t_{1},...,t_{m}\right)$ is $\mathcal{C}^{\infty }$.
Let us show that the Jacobian matrix of $g_{r}$ in
$\left(0,\ldots,0\right)$ is invertible. Consider
$$g_{r}\left(t_{1},...,t_{m}\right)\equiv\left(G_{1}\left(t_{1},...,t_{m}\right),...,G_{m}
\left(t_{1},...,t_{m}\right)\right).$$
We have
$$
\det\left(\begin{array}{ccc} \frac{\partial G_{1}}{\partial
t_{1}}&\cdots &\frac{\partial G_{m}}{\partial t_{1}}\\
\vdots &\ddots &\vdots \\
\frac{\partial G_{1}}{\partial t_{m}}&\cdots &\frac{\partial
G_{m}}{\partial t_{m}}
\end{array}\right)
=\det\left(\begin{array}{ccc} \frac{\partial }{\partial t_{1}}
g_{t_{m}}^{X_{m}}\circ \cdots \circ g_{t_{1}}^{X_{1}}\left( x\right) \\
\vdots \\
\frac{\partial }{\partial t_{m}}g_{t_{m}}^{X_{m}}\circ \cdots
\circ g_{t_{1}}^{X_{1}}\left( x\right)
\end{array}\right)
\neq 0,
$$
because the vector fields $X_{1},...,X_{m}$ are linearly
independent at each point of $M$. According to the local inversion
theorem, there exists a sufficiently small neighborhood $V\subset
U$ of $\left(0,\ldots,0\right)$ and a neighborhood $W$ of $x$ such
that $g_{r}$  induces a bijection of $V$ on $W$ whose inverse
$g_{r}^{-1}:W\rightarrow V$, is $\mathcal{C}^{\infty }$. In other
words, $g_{r}$ is a diffeomorphism of $V$ over $g_{r}(V)$. This
result is local because even if the above Jacobian matrix is
invertible for any $\left(t_{1},...,t_{m}\right)$, then the
inverse "global" of $g_{r}$ does not necessarily exist.\\
$b)$ The application $g$ is surjective. Indeed, let
$(t_{1},...,t_{m})\in \mathbb{R}^{m}$ such that :
$$g(t_{1},...,t_{m})=g_{t_{m}}^{X_{m}}\circ \cdots \circ
g_{t_{1}}^{X_{1}}(x)=y\in M.$$ We showed in the part $a)$ that $g$
is a local diffeomorphism. So for every point $x_{1}$ contained in
a neighborhood of $x$, there exists $\left(t_{1},...,t_{m}\right)
\in \mathbb{R}^{m}$ such that : $$g_{t_{m}}^{X_{m}}\circ \cdots
\circ g_{t_{1}}^{X_{1}}(x)=x_{1}.$$ Since the variety $M$ is
connected, we can connect the point $x$ to the point $y$ by a
curve $\mathcal{C}$. Let $B_{1}$ be an open ball in $M$ containing
the point $x_{1}$. This ball exists since $M$ is compact. Let
$x_{2}\in\mathcal{C}$ such that $x_{2}$ be contained in the ball
$B_{1}$. We reason as before, the map $g$ being a local
diffeomorphism, then there exists $\left(t_{1}^{\prime
},...,t_{m}^{\prime }\right) \in \mathbb{R}^{m}$ such that :
$$\left(g_{t_{m}}^{X_{m}}\right)^{\prime}\circ \cdots \circ
\left(g_{t_{1}}^{X_{1}}\right)^{\prime}(x_{1})=x_{2}.$$ Hence,
$$x_{2}=\left(g_{t_{m}}^{X_{m}}\right)^{\prime}+t_{m}\circ \cdots
\circ \left(g_{t_{1}}^{X_{1}}\right)^{\prime}+t_{1}(x).$$
Similarly, let $B_{2}$ be an open ball in $ M $ containing the
point $x_{2}$ and let $x_{3}\in\mathcal {C}$ such that $x_{3}$ be
either contained in the ball $B_{2}$. Since the application $g$ is
a local diffeomorphism, then there exists $\left(t_{1}^{\prime
\prime },...,t_{m}^{\prime \prime}\right) \in \mathbb{R}^{m}$ such
that : $$\left(g_{t_{m}}^{X_{m}}\right)^{\prime\prime}\circ \cdots
\circ
\left(g_{t_{1}}^{X_{1}}\right)^{\prime\prime}(x_{2})=x_{3}.$$ So
$$x_{3}=\left(g_{t_{m}}^{X_{m}}\right)^{\prime\prime}+t_{m}^{\prime
}+ t_{m}\circ \cdots \circ
\left(g_{t_{1}}^{X_{1}}\right)^{\prime\prime}+t_{1}^{\prime
}+t_{1}(x).$$ Continuing this way, we show (after a finite number
$k$ of steps) the existence of a point $\left(t_{1}^{\left(
k-1\right) },...,t_{m}^{\left( k-1\right)}\right) \in
\mathbb{R}^{m}$, such that :
$$\left(g_{t_{m}}^{X_{m}}\right)^{(k-1)}\circ \cdots \circ
\left(g_{t_{1}}^{X_{1}}\right)^{(k-1)}(x_{k-1})=x_{k},$$ where
$x_{k}\in $ $\mathcal{C}$, $x_{k}$ is contained in an open ball
$B_{k-1}$ of $M$, with $B_{k-1}\ni $ $x_{k-1}$. Therefore, for $k$
finite, we have
$$x_{k}=\left(g_{t_{m}}^{X_{m}}\right)^{(k-1)}+t_{m}^{(k-2)}
+\cdots +t_{m}^{\prime }+t_{m}\circ \cdots \circ
\left(g_{t_{1}}^{X_{1}}\right)^{(k-1)}+t_{1}^{(k-2)}+ \cdots
+t_{1}^{\prime }+t_{1}(x).$$ This construction shows that in a
finite number $k$ of steps, we can cover the curve $\mathcal{C}$
connecting the point $x$ to the point $y$ by neighborhoods of $x$;
the point $y$ playing the role of $x_{k}$.
%$$\includegraphics[scale=0.4]{Figure28.jpg}$$
Note that the application $g$ can not be injective. In fact, if
$g$ is injective, we would have, according to part a), a bijection
between a compact $M$ and a noncompact $\mathbb{R}^{m}$, which is
absurd.\\
$c)$ The stationary group
$$\Lambda=\left\{\left( t_{1},...,t_{m}\right) \in
\mathbb{R}^{m}:g\left(t_{1},...,t_{m}\right)=g_{t_{m}}^{X_{m}}\circ
\cdots \circ g_{t_{1}}^{X_{1}}\left( x\right)=x\right\},$$ is a
discrete subgroup of $\mathbb{R}^{m}$ independent of point $x\in
M$. Indeed, let us first note that $\Lambda \neq \emptyset$
because $(0,...,0)\in \Lambda$. Let $\left(t_{1},...,t_{m}\right)
\in \Lambda$, $\left( t_{1}^{\prime },...,t_{m}^{\prime}\right)
\in \Lambda$. We have
$$g\left(t_{1},...,t_{m}\right)=g\left(t_{1}^{\prime},...,t_{m}^{\prime}\right)=x.$$
As the vector fields $X_{1},...,X_{m}$ are commutative, then
\begin{eqnarray}
g\left( t_{1}+t_{1}^{\prime },...,t_{m}+t_{m}^{\prime }\right)
&=&g_{t_{m}+t_{m}^{\prime }}^{X_{m}}\circ \cdots \circ
g_{t_{1}+t_{1}^{\prime }}^{X_{1}}\left( x\right) ,\nonumber\\
&=&g_{t_{m}^{\prime }}^{X_{m}}\circ \cdots \circ g_{t_{1}^{\prime
}}^{X_{1}}\circ g_{t_{m}}^{X_{m}}\circ \cdots \circ
g_{t_{1}}^{X_{1}}\left( x\right) ,\nonumber\\
&=&g_{t_{m}^{\prime }}^{X_{m}}\circ \cdots \circ g_{t_{1}^{\prime
}}^{X_{1}}\left( x\right) ,\nonumber\\
&=&x,\nonumber
\end{eqnarray}
and
\begin{eqnarray}
g\left( -t_{1},...,-t_{m}\right)&=&g_{-t_{m}}^{X_{m}} \circ
\cdots \circ g_{-t_{1}}^{X_{1}}\left( x\right),\nonumber\\
&=&g_{-t_{m}}^{X_{m}}\circ \cdots \circ g_{-t_{1}}^{X_{1}} \circ
g_{t_{m}}^{X_{m}}\circ \cdots \circ g_{t_{1}}^{X_{1}}
\left( x\right),\nonumber\\
&=&g_{-t_{m}}^{X_{m}}\circ \cdots \circ g_{-t_{1}}^{X_{1}} \circ
g_{t_{1}}^{X_{1}}\circ \cdots \circ g_{t_{m}}^{X_{m}}
\left( x\right),\nonumber\\
&=&g_{-t_{m}}^{X_{m}}\circ \cdots \circ g_{-t_{2}}^{X_{2}} \circ
g_{t_{2}}^{X_{2}}\circ \cdots \circ g_{t_{m}}^{X_{m}}
\left( x\right),\nonumber\\
&\vdots& \nonumber\\
&=&g_{-t_{m}}^{X_{m}}\circ g_{t_{m}}^{X_{m}}\left( x\right),\nonumber\\
&=&x.\nonumber
\end{eqnarray}
Hence
$\left(t_{1}+t_{1}^{\prime},...,t_{m}+t_{m}^{\prime}\right)\in\Lambda$
and $\left(-t_{1},...,-t_{m}\right)\in\Lambda$. Therefore
$\Lambda$ is stable for addition, the inverse of
$\left(t_{1},...,t_{m}\right)$ is $\left(
-t_{1},...,-t_{m}\right)$ and consequently $\Lambda$ is a subgroup
of $\mathbb{R}^{m}$. We show that $\Lambda$ is independent of $x$.
Let
$$\Lambda^{\prime}=\left\{\left(t_{1}^{\prime},...,t_{m}^{\prime}\right)\in\mathbb{R}^{m} :
g\left(t_{1}^{\prime},...,t_{m}^{\prime}\right)=g_{t_{m}^{\prime}}^{X_{m}}\circ\cdots\circ
g_{t_{1}^{\prime}}^{X_{1}}(y)=y\right\}.$$ By surjectivity, one
can find $\left( s_{1},...,s_{m}\right) \in \mathbb{R}^{m}$ such
that :
$$g_{s_{m}}^{X_{m}}\circ \cdots \circ g_{s_{1}}^{X_{1}}(x)
=y.$$ Let $\left(t_{1}^{\prime },...,t_{m}^{\prime}\right)
\in\Lambda^{\prime}$. We have
\begin{eqnarray}
g_{t_{m}^{\prime }}^{X_{m}}\circ \cdots \circ g_{t_{1}^{\prime
}}^{X_{1}}(y)&=&y,\nonumber\\
g_{t_{m}^{\prime }}^{X_{m}}\circ \cdots \circ g_{t_{1}^{\prime
}}^{X_{1}}\circ g_{s_{m}}^{X_{m}}\circ \cdots \circ
g_{s_{1}}^{X_{1}}(x)&=&g_{s_{m}}^{X_{m}}\circ \cdots
\circ g_{s_{1}}^{X_{1}}(x),\nonumber\\
g_{-s_{m}+t_{m}^{\prime }+s_{m}}^{X_{m}}\circ \cdots\circ
g_{-s_{1}+t_{1}^{\prime }+s_{1}}^{X_{1}}(x)&=&x,\nonumber\\
g_{t_{m}^{\prime }}^{X_{m}}\circ \cdots \circ g_{t_{1}^{\prime
}}^{X_{1}}(x)&=&x.\nonumber
\end{eqnarray}
Therefore,
$\left(t_{1}^{\prime},...,t_{m}^{\prime}\right)\in\Lambda$ and
therefore $\Lambda$ does not depend on $x$. To show that $\Lambda$
is discrete, we consider a neighborhood $V$ sufficient small of
the point $(0,...,0)$ and a neighborhood $W$ of the point $x$.
From $a)$, the application $g$ is a local diffeomorphism, so $g :
V\longrightarrow W$, is bijective and consequently no point of $W
\backslash \left\{(0,...,0)\right\}$ is sent on $x$; the points of
the subgroup $\Lambda$ have no accumulation point in
$\mathbb{R}^{m}$.\\
$d)$ The variety $M$ is diffeomorphic to a $m$-dimensional real
torus. Indeed, let
$$T^{k}\times\mathbb{R}^{m-k}=\{(\varphi_1,...,\varphi_k;u_1,...,u_{m-k})\},
\quad (\varphi_1,...,\varphi_k)\mbox{ mod.}2\pi$$ be the direct
product of $k$ circles and $m-k$ straight lines and consider the
application $$\pi:\mathbb{R}^m\longrightarrow
T^{k}\times\mathbb{R}^{m-k},$$ defined by
$$\pi(\varphi_1,...,\varphi_k;u_1,...,u_{m-k})=((\varphi_1,...,\varphi_k)\mbox{ mod.}2\pi;(u_1,...,u_{m-k})).$$
The points $f_1,...,f_k\in\mathbb{R}^m$ where each $f_i$ has the
coordinates $$\varphi_i=2\pi, \quad\varphi_j=0, \quad
u_1=\cdots=u_{m-k}=0,$$ are sent in $0$ by this application. Let
us first note that the stationary group $\Lambda$ (see point $c)$
can be written in the form $$\Lambda=\mathbb{Z}e_{1}\oplus \cdots
\oplus \mathbb{Z}e_{k}, \quad1\leq k\leq m,$$ where
$e_{1},...,e_{m}$ are linearly independent vectors. Indeed, to fix
ideas, let us take $m=2$, i.e.,
$$\Lambda=\left\{\left(t_{1},t_{2}\right)\in\mathbb{R}^{2}:
g\left(t_{1},t_{2}\right)=g_{t_{2}}^{X_{2}}\circ
g_{t_{1}}^{X_{1}}(x)=x\right\}.$$ Three cases are possible : $i)$
$\Lambda=\{0\}$, $ii)$ $\Lambda=\mathbb{Z}e_{1}$, $iii)$
$\Lambda=\mathbb{Z}e_{1}\oplus \mathbb{Z}e_{2}$. The first case is
to be rejected because we have a diffeomorphism between a
non-compact $\mathbb{Z}^{2}/\{0\}$ and a compact $M$, which is
absurd. The second case $\mathbb{Z}^{2}/$ $\mathbb{Z}e_{1}$ (a
cylinder) is also to be rejected for the same reasons as in the
first case. It remains the last case, which is valid because
$\mathbb{Z}^{2}/$ $\mathbb{Z}e_{1}\oplus \mathbb{Z}e_{2}$ is a
$2$-dimensional torus. In general, for every discrete subgroup of
$\mathbb{R}^{m}$, there exist $k$ linearly independent vectors
such that this group is the set of all their integer linear
combinations. Let $e_{1},...,e_{k}\in\Lambda\subset\mathbb{R}^m$
be generators of the stationary group $\Lambda$. We now apply the
vector space
$\mathbb{R}^m=\{(\varphi_1,...,\varphi_k;u_1,...,u_{m-k})\}$ in a
surjective way over space vector $\mathbb{R}^m=\{(t_1,...,t_m)\}$
such that the vectors $f_i$ are transformed into $e_i$. Let
$h:\mathbb{R}^m\longrightarrow\mathbb{R}^m$ be such an isomorphism
and notice that
$\mathbb{R}^m=\{(\varphi_1,...,\varphi_k;u_1,...,u_{m-k})\}$
(resp. $\mathbb{R}^m=\{(t_1,...,t_m)\}$) determines charts of
$T^{k}\times\mathbb{R}^{m-k}$ (respectively of the variety $M$).
The application $h$ determines a diffeomorphism
$$\widetilde{h}:T^{k}\times\mathbb{R}^{m-k}\longrightarrow M,$$ and
since by hypothesis $M$ is compact, then $k=m$ and consequently
$M$ is a $m$-dimensional torus. Let's check this out in more
detail. Since $\Lambda$ is the kernel of $g$, there exists a
canonical surjection
$$\widetilde{h}:\mathbb{R}^{m}/\Lambda \rightarrow M,\text{ }
\left[\left( t_{1},...,t_{m}\right)\right] \mapsto \widetilde{h}
\left[\left( t_{1},...,t_{m}\right)\right]=g_{t_{m}}^{X_{m}}\circ
\cdots \circ g_{t_{1}}^{X_{1}}(x).$$ Indeed, let
$\left(t_{1},...,t_{m}\right)$ et $\left(s_{1},...,s_{m}\right)$
such that :
$$\widetilde{h}\left[\left(t_{1},...,t_{m}\right)\right]=\widetilde{h}\left[\left(s_{1},...,s_{m}\right)\right].$$
We have $$g_{t_{m}}^{X_{m}}\circ \cdots \circ
g_{t_{1}}^{X_{1}}(x)=g_{s_{m}}^{X_{m}}\circ \cdots \circ
g_{s_{1}}^{X_{1}}(x),$$ hence
\begin{eqnarray}
&&g_{-s_{1}}^{X_{1}}\circ \cdots \circ g_{-s_{m}}^{X_{m}}\circ
g_{t_{m}}^{X_{m}}\circ \cdots \circ g_{t_{1}}^{X_{1}}(x)\nonumber\\
&&\qquad\qquad=g_{-s_{1}}^{X_{1}}\circ \cdots \circ
g_{-s_{m}}^{X_{m}}\circ
g_{s_{m}}^{X_{m}}\circ \cdots \circ g_{s_{1}}^{X_{1}}(x),\nonumber\\
&&\qquad\qquad=g_{-s_{1}}^{X_{1}}\circ \cdots \circ
g_{-s_{m-1}}^{X_{m-1}} \circ g_{s_{m-1}}^{X_{m-1}}
\circ \cdots \circ g_{s_{1}}^{X_{1}}(x),\nonumber\\
&&\qquad\qquad\vdots \nonumber\\
&&\qquad\qquad=g_{-s_{1}}^{X_{1}}\circ g_{s_{1}}^{X_{1}}(x),\nonumber\\
&&\qquad\qquad=x.\nonumber
\end{eqnarray}
Since $X_{1},...,X_ {m}$ are commutative, then
$$g_{t_{m}-s_{m}}^{X_{m}}\circ \cdots \circ g_{t_{1}-s_{1}}^{X_{1}}
\left( x\right)=x.$$ Consequently, we have
$$\left[\left(t_{1}-s_{1},...,t_{m}-s_{m}\right) \right]=0,
\qquad \left[\left( t_{1},...,t_{m}\right)
-\left(s_{1},...,s_{m}\right) \right]=0,$$ and $$\left[ \left(
t_{1},...,t_{m}\right) \right]=\left[\left( s_{1},...,s_{m}\right)
\right].$$ So $\widetilde{h}$ is a diffeomorphism and the proof of
part 1) is complete.

\textbf{2)} By hypothesis the variety $M_{c}$ is compact and
connected. Therefore, from the part 1), it is enough to show that
$M_{c}$ is differentiable, of dimension $n$ and that it is
equipped with $n$ commutative vectors fields. The
differentiability of this variety arises from the implicit
function theorem since the vectors $J\frac{\partial
H_{1}}{\partial x}, \ldots, J\frac {\partial H_{n}}{\partial x}$
are assumed to be independent. As $m=2n$, then the first integrals
$H_{i}(x_{1},...,x_{2n})$ are functions of the variables
$x_{1},...,x_{n},x_{n+1},...,x_{2n}$. Therefore, $$\dim \left\{
x\in M : H_{i}=c_{i}\right\}=2n-1,$$ and
$$\dim \left( \left\{ x\in M:H_{i}=c_{i}\right\}\cap \left\{ x\in M:H_{j}
=c_{j}\right\} \right) =2n-2,\text{}i\neq j,$$ and so $\dim
M_{c}=n$. Let $X_{i}$ et $X_{j},$ $1\leq i,j\leq n$, be
differentiable ($\mathcal{C}^{\infty }$) vector fields on $M$, so
on the variety $M_{c}$ also. Let us define the differential
operator $L_{X}$ by $$L_{X}:\mathcal{C}^{\infty }\left(
M_{c}\right) \longrightarrow \mathcal{C}^{\infty }\left(
M_{c}\right), \quad F\longmapsto L_{X}F,$$ such that:
$$L_{X}F(x)=\left. \frac{d}{dt}F\left( g_{t}^{X}(x)\right) \right| _{t=0},
\text{ }x\in M_{c}.$$ We have $$L_{X_{i}}F=\left\{
F,H_{i}\right\}, \qquad L_{Xj}L_{X_{i}}F=\left\{ \left\{
F,H_{i}\right\},H_{j}\right\},$$ and
\begin{eqnarray}
L_{Xi}L_{X_{j}}F-L_{Xj}L_{X_{i}}F&=&\left\{ \left\{
F,H_{j}\right\} ,H_{i}\right\} -\left\{ \left\{ F,H_{i}\right\}
,H_{j}\right\}
,\nonumber\\
&=&-\left\{ \left\{ H_{j},F\right\} ,H_{i}\right\} -\left\{
\left\{F,H_{i}\right\} ,H_{j}\right\} ,\nonumber\\
&=&\left\{ \left\{ H_{i},H_{j}\right\} ,F\right\} ,\nonumber
\end{eqnarray}
according to the identity of Jacobi. Since $H_{i}$ and $H_{j}$ are
in involution, then $\left[L_{X_{i}}, L_{X_{j}}\right]=0$. The
construction of the angular coordinates $\varphi_1,...,\varphi_m$
mod. $2\pi$ on the variety $M$ is obviously valid on the invariant
variety $M_c$. note that
$$(\varphi_1,...,\varphi_m)=h^{-1}(t_1,...,t_m),$$ and that the
angular coordinates $\varphi_1,...,\varphi_m$ vary uniformly under
the action of the Hamiltonian flow $H$, i.e.,
$$\frac{d\varphi_k }{dt_i}=\{H_i,\varphi_k\}=\omega_i(c),\qquad\omega_i(c)=\text{constants.}$$
In other words, the motion is quasi-periodic on the invariant
torus $M_c$. Finally, to show that the equations of the problem
are integrable by quadratures as well as several information about
the variables called action-angle, one will consult with profit
[4]. The demonstration of theorem ends. $\square$

If we restrict ourselves to an invariant open set, we can always
assume that the fibers of $M_c$ (where $c$ is a regular value) are
connected. The tori obtained in the theorem are Lagrangian
sub-varieties. If $M_c$ is not compact but the flow of each of the
vector fields $X_{H_k}$ is complete on $M_c$ (a vector field is
called complete if every one of its flow curves exist for all
time), then $M_c$ is diffeomorphic to a cylinder
$\mathbb{R}^k\times T^{n-k}$ under which the vector fields
$X_{H_k}$ are mapped to linear vector fields.

\begin{Exmp}
The rank of the matrix $J$ is even. Indeed, let $\lambda$ be the
eigenvalue associated with the eigenvector $Z$. We have
$$JZ=\lambda Z, \quad Z\neq 0,$$ and $$Z^*JZ=\lambda Z^*Z, \qquad Z^*\equiv
\overline{Z}^\top,$$ where $\lambda=\frac{Z^*JZ}{Z^*Z}$. Since
$\overline{J}=J$ and $J^\top=-J$, then
$$\overline{Z^*JZ}=Z^\top\overline{JZ}=Z^\top J \overline{Z}=
(Z^\top J \overline{Z})^\top=Z^*J^\top Z=-Z^*JZ,$$ which implies
that $Z^*JZ$ is either zero or imaginary pure. Since $Z^*Z$ is
real, it follows that all the eigenvalues of $J$ are either null
or imaginary pure. Now $J\overline{Z}=\overline{\lambda Z}$, so if
$\lambda$ is an eigenvalue, then $\overline{\lambda}$ is also an
eigenvalue. Consequently, the eigenvalues (non-zero) of $J$ come
in pairs, hence the result.
\end{Exmp}

As a consequence, we obtain the concept of complete integrability
of a Hamiltonian system (15) with $x\in M=\mathbb{R}^m$. For the
sake of clarity, we shall distinguish two cases :

$\textbf{a})$ \emph{Case $1$} : $\det J\neq 0$. The rank of the
matrix $J$ is even (example 19), $m=2n$. A Hamiltonian system
(15), $x\in M=\mathbb{R}^m$, is completely integrable or
Liouville-integrable if there exist $n$ firsts integrals
$H_{1}=H,H_{2},\ldots,H_{n}$ in involution, i.e.,
$\{H_{k},H_{l}\}=0,\text{ }1\leq k,l\leq n,$ with linearly
independent gradients, i.e., $dH_{1}\wedge ...\wedge dH_{n}\neq
0.$ For generic $c=(c_1,...,c_n)$ the level set
$$M_{c}=\bigcap_{i=1}^{n}\left\{ x\in M:H_{i}(x)=c_{i},\text{
}c_{i}\in \mathbb{R}\right\},$$ will be an $n$-manifold. By the
Arnold-Liouville theorem, if $M_{c}$ is compact and connected, it
is diffeomorphic to an $n$-dimensional torus
$\mathbb{T}^n=\mathbb{R}^{n}/\mathbb{Z}^{n}$ and each vector field
will define a linear flow there. In some open neighborhood of the
torus there are coordinates
$s_{1},\ldots,s_{n},\varphi_{1},\ldots,\varphi_{n}$ in which
$\omega$ takes the form $$\omega=\sum_{k=1}^{n}ds_{k}\wedge
d\varphi_{k}.$$ Here the functions $s_{k}$ (called
action-variables) give coordinates in the direction transverse to
the torus and can be expressed functionally in terms of the first
integrals $H_{k}.$ The functions $\varphi_{k}$ (called
angle-variables) give standard angular coordinates on the torus,
and every vector field $X_{H_{k}}$ can be written in the form
$$\dot\varphi_{k}=h_{k}\left(s_{1},\ldots,s_{n}\right),$$ that is,
its integral trajectories define a conditionally-periodic motion
on the torus. In a neighborhood of the torus the Hamiltonian
vector field $X_{H_{k}}$ take the following form $$\dot
s_{k}=0,\qquad \dot
\varphi_{k}=h_{k}\left(s_{1},\ldots,s_{n}\right),$$ and can be
solved by quadratures.

$\textbf{b})$ \emph{Case $2$} : $\det J=0$. We reduce the problem
to $m=2n+k$ and we look for $k$ Casimir functions
$H_{n+1},...,H_{n+k},$ leading to identically zero Hamiltonian
vector fields $$J\frac{\partial H_{n+i}}{\partial x}=0,\quad 1\leq
i\leq k.$$ In other words, the system is Hamiltonian on a generic
symplectic manifold $$\bigcap_{i=n+1}^{n+k}\left\{x\in
\mathbb{R}^{m}:H_{i}(x)=c_{i}\right\},$$ of dimension $m-k=2n$. If
for most values of $c_i\in \mathbb{R}$, the invariant manifolds
$$\bigcap_{i=1}^{n+k}\left\{x\in
\mathbb{R}^{m}:H_{i}(x)=c_{i}\right\},$$ are compact and
connected, then they are $n$-dimensional tori
$\mathbb{T}^n=\mathbb{R}^n/\mathbb{Z}^n$ by the Arnold-Liouville
theorem and the Hamiltonian flow is linear in angular coordinates
of the torus.

\begin{Exmp}
The simple pendulum and the harmonic oscillator are trivially
integrable systems (any $2$-dimensional Hamiltonian system where
the set of non-fixed points is dense, is integrable. Let
$T^*\mathbb{R}^n$ with coordinates $q_1,...,q_n, p_1,...,p_n$. The
system corresponding to the Hamiltonian of the harmonic oscillator
is integrable
$$H=\frac{1}{2}\sum_{j=1}^n(p_j^2+\lambda_jq_j^2).$$ The
Hamiltonian structure is defined by the Poisson bracket $$\{F,
H\}=\sum_{j=1}^n\left(\frac{\partial F}{\partial
q_j}\frac{\partial H}{\partial p_j}-\frac{\partial F}{\partial
p_j}\frac{\partial H}{\partial q_j}\right).$$ The Hamiltonian
field corresponding to $H$ is written explicitly
$$
\dot{q}_j=p_j,\qquad \dot{p}_j=-2\lambda_jq_j, \quad j=1,...,n$$
and admits the following first $n$ integral :
$$H_j=\frac{1}{2}p_j^2+\lambda_jq_j^2, \quad 1\leq j\leq n.$$ The latter
are independent, in involution and the system in question is
integrable.
\end{Exmp}

\section{Rotation of a solid body about a fixed point and $SU(2)$ Yang-Mills equations}

\subsection{The problem of the rotation of a solid body about a
fixed point}

One of the most fundamental problems of mechanics is the study of
the motion of rotation of a solid body around a fixed point. The
differential equations of this problem are written in the form
\begin{eqnarray}
\dot{M}&=&M\wedge \Omega +\mu g\text{ }\Gamma \wedge L,\\
\dot{\Gamma }&=&\Gamma \wedge \Omega ,\nonumber
\end{eqnarray}
where $\wedge$ is the vector product in $\mathbb{R}^3$,
$M=\left(m_{1},m_{2},m_{3}\right)$ the angular momentum of the
solid, $\Omega=(\frac{m_{1}}{I_1}, \frac{m_{2}}{I_2},
\frac{m_{3}}{I_3})$ the angular velocity, $I_{1},I_{2}$ and
$I_{3}$, moments of inertia,
$\Gamma=\left(\gamma_{1},\gamma_{2},\gamma_{3}\right)$ the unitary
vertical vector, $\mu$ the mass of the solid, $g$ the acceleration
of gravity, and finally, $L=\left( l_{1},l_{2},l_{3}\right)$ the
unit vector originating from the fixed point and directed towards
the center of gravity; all these vectors are considered in a
mobile system whose coordinates are fixed to the main axes of
inertia. The configuration space of a solid with a fixed point is
the group of rotations $SO(3)$. This is generated by the
one-parameter subgroup of rotations
\begin{eqnarray}
{A_{1}}&=&\left(\begin{array}{ccc}
1&0&0\\
0&\cos t&-\sin t\\
0&\sin t&\cos t
\end{array}\right),\nonumber\\
{A_{2}}&=&\left(\begin{array}{ccc}
\cos t&0&\sin t\\
0&1&0\\
-\sin t&0&\cos t
\end{array}\right),\nonumber\\
{A_{3}}&=&\left(\begin{array}{ccc}
\cos t&-\sin t&0\\
\sin t&\cos t&0\\
0&0&1
\end{array}\right).\nonumber
\end{eqnarray}
Recall that this is the group of $n\times n$ orthogonal matrices
$A$ and the motion of this solid is described by a curve on this
group. The angular velocity space of all rotations (the set of
derivatives $\left.\dot{A}(t)\right|_{t=0}$ of the differentiable
curves in $SO(3)$ passing through the identity in $t=0$ :
$A(0)=I$) is the Lie algebra of the group $SO(3)$; it is the
algebra $so(3)$ of the $3\times3$ antisymmetric matrices. This
algebra is generated as a vector space by the matrices
$$
{e_{1}}=\left.\dot{A}_1(t)\right|_{t=0}=\left(\begin{array}{ccc}
0&0&0\\
0&0&-1\\
0&1&0
\end{array}\right),
$$
$$
{e_{2}}=\left.\dot{A}_2(t)\right|_{t=0}=\left(\begin{array}{ccc}
0&0&1\\
0&0&0\\
-1&0&0
\end{array}\right),
$$
$$
{e_{3}}=\left.\dot{A}_3(t)\right|_{t=0}=\left(\begin{array}{ccc}
0&-1&0\\
1&0&0\\
0&0&0
\end{array}\right),
$$
which verify the commutation relations :
$$\left[e_{1},e_{2}\right]=e_{3}, \quad \left[e_{2},e_{3}\right]=e_{1},
\quad \left[e_{3},e_{1}\right]=e_{2}.$$ We will use the fact that
if we identify $so(3)$ to $\mathbb{R}^{3}$ by sending
$\left(e_{1}, e_{2}, e_{3}\right)$ on the canonical basis of
$\mathbb{R}^{3}$, the bracket of $so(3)$ corresponds to the vector
product. In other words, consider the application
$$\mathbb{R}^{3}\longrightarrow so(3),\text{ }a=
\left( a_{1},a_{2},a_{3}\right) \longmapsto
{A}=\left(\begin{array}{ccc}
0&-a_{3}&a_{2}\\
a_{3}&0&-a_{1}\\
-a_{2}&a_{1}&0
\end{array}\right),
$$
which defines an isomorphism between Lie algebras
$\left(\mathbb{R}^{3},\wedge \right)$ et $(so(3),[ ,])$ where
$$a\wedge b\longmapsto [A,B]=AB-BA.$$
By using this isomorphism, the system (24) can be rewritten in the
form
\begin{eqnarray}
\dot{M}&=&\left[M,\Omega\right]+\mu g\text{ }\left[\Gamma,L\right],\\
\dot{\Gamma}&=&\left[\Gamma ,\Omega \right] ,\nonumber
\end{eqnarray}
where
\begin{eqnarray}
M&=&\left( M_{ij}\right)_{1\leq i,j\leq 3}\equiv
\sum_{i=1}^{3}m_{i}e_{i}\equiv\left(\begin{array}{ccc}
0&-m_{3}&m_{2}\\
m_{3}&0&-m_{1}\\
-m_{2}&m_{1}&0
\end{array}\right)\in so\left( 3\right),\nonumber\\
\Omega &=&\left( \Omega _{ij}\right) _{1\leq i,j\leq 3}\equiv
\sum_{i=1}^{3}\omega _{i}e_{i}\equiv \left(\begin{array}{ccc}
0&-\omega_{3}&\omega_{2}\\
\omega_{3}&0&-\omega_{1}\\
-\omega_{2}&\omega_{1}&0
\end{array}\right)\in so\left( 3\right),\nonumber\\
\Gamma &=&\left(\gamma _{ij}\right) _{1\leq i,j\leq 3}\equiv
\sum_{i=1}^{3}\gamma _{i}e_{i}\equiv \left(\begin{array}{ccc}
0&-\gamma_{3}&\gamma_{2}\\
\gamma_{3}&0&-\gamma_{1}\\
-\gamma_{2}&\gamma_{1}&0
\end{array}\right)\in so\left( 3\right),\nonumber
\end{eqnarray}
and
$$L = \left(\begin{array}{ccc}
0&-l_{3}&l_{2}\\
l_{3}&0&-l_{1}\\
-l_{2}&l_{1}&0
\end{array}\right)\in so\left( 3\right).
$$
Taking into account that $M=I\Omega$, then the above equations
(25) become
\begin{eqnarray}
\dot{M}&=&\left[ M,\Lambda M\right] +\mu g\text{ }\left[
\Gamma ,L\right],\\
\dot{\Gamma}&=&\left[ \Gamma ,\Lambda M\right] ,\nonumber
\end{eqnarray}
where
$$\Lambda M=\equiv
\sum_{i=1}^{3}\lambda _{i}m_{i}e_{i}\equiv
\left(\begin{array}{ccc}
0&-\lambda _{3}m_{3}&\lambda _{2}m_{2}\\
\lambda _{3}m_{3}&0&-\lambda _{1}m_{1}\\
-\lambda _{2}m_{2}&\lambda _{1}m_{1}&0
\end{array}\right)\in so\left( 3\right),\quad
\lambda _{i}\equiv \left.\frac{1}{I_{i}}\right.$$

The resolution of this problem was analyzed first by Euler [10]
and in 1758, he published the equations (case $\mu=0$) which carry
his name. Euler's equations were integrated by Jacobi [19] in
terms of elliptic functions and around 1851, Poinsot [34] gave
them a remarkable geometric interpretation. Before, around 1815
Lagrange [21] found another case ($I_{1}=I_{2}$, $l_{1}=l_{2}=0$)
of integrability, that subsequently Poisson has examined at length
thereafter. The problem continued to attract mathematicians but
for a long time no new results could be obtained. It was then
around 1888-1989 that a memoir [20], of the highest interest,
appears containing a new case ($I_{1}=I_{2}=2I_{3}$, $l_{3}=0$) of
integrability discovered by Kowalewski. For this remarkable work,
Kowalewski was awarded the Bordin Prize of the Paris Academy of
Sciences. In fact, although Kowalewski's work is quite important,
it is not at all clear why there would be no other new cases of
integrability. This was to be the starting point of a series of
fierce research on the question of the existence of new cases of
integrability. Moreover, among the remarkable results obtained by
Poincar\'{e} [33] with the aid of the periodic solutions of the
equations of dynamics, we find the following (around 1891) : in
order to exist in the motion of a solid body around of a fixed
point, an algebraic first integral not being reduced to a
combination of the classical integrals, it is necessary that the
ellipsoid of inertia relative to the point of suspension is of
revolution. In 1896, R. Liouville (not to be confused with Joseph
Liouville, well known in complex analysis) also competed for the
Bordin prize, presented a paper [30] indicating necessary and
sufficient conditions ($I_{3}=0$, $2I_{3}/I_1=\mbox{integer}$) of
existence of a fourth algebraic integral. These conditions have
been reproduced in most conventional treatises (eg Whittaker [43])
and in scientific journals. And it was not until the year 1906,
when Husson [18], working under the direction of Appell and
Painlev\'{e}, discovered an erroneous demonstration in the work of
Liouville. Indeed, paragraphs I and III of Liouville's
dissertation devoted to the search for the necessary conditions
seem at first satisfactory, but a more careful study shows that
the demonstrations are at least insufficient and that it is
impossible to accept conclusions. In fact, although the conditions
found by Liouville are necessary, they can not be deduced from the
calculations indicated and these conditions are not sufficient.
And it was Husson who first solved completely the question of
looking for new cases of integrability. Inspired by Poincar\'{e}'s
research on the problem of the three bodies and Painlev\'{e} on
the generalization of Bruns's theorem, Husson demonstrated that
any algebraic integral is a combination of classical integrals
except in the cases of Euler, Lagrange and Kowalewski. Moreover,
the question of the existence of analytic integrals has been
studied rigorously by Ziglin [45, 46] and Holmes-Marsden [96].
Towards the end of this subsection, we will mention some special
cases : cases of Hesse-Appel'rot [16, 3], Goryachev-Chaplygin [12,
8] and Bobylev-Steklov [6, 36].

In the case of the Euler rigid body motion, we have
$l_{1}=l_{2}=l_{3}=0$, that is, the fixed point is its center of
gravity. The Euler rigid body motion [10] (also called
Euler-Poinsot motion [34] of the solid) express the free motion of
a rigid body around a fixed point. Then the motion of the body is
governed by $\dot M=[M,\Lambda M]$, and is explicitly given by
\begin{eqnarray}
\dot m_{1}&=&\left(\lambda_{3}-\lambda_{2}\right)m_{2}m_{3},\nonumber\\
\dot m_{2}&=&\left(\lambda_{1}-\lambda_{3}\right)m_{1}m_{3},\\
\dot
m_{3}&=&\left(\lambda_{2}-\lambda_{1}\right)m_{1}m_{2},\nonumber
\end{eqnarray}
and (see example 10) can be written as a Hamiltonian vector field
$$\dot x=J\frac{\partial H}{\partial x},
\quad x=\left(m_{1},m_{2},m_{3}\right)^{\intercal},$$ with the
Hamiltonian
$$H=\frac{1}{2}\left(\lambda_{1}m_{1}^{2}+\lambda_{2}m_{2}^{2}+\lambda_{3}m_{3}^{2}\right),$$
and
$$J=\left(\begin{array}{ccc}
0&-m_{3}&m_{2}\\
m_{3}&0&-m_{1}\\
-m_{2}&m_{1}&0
\end{array}\right)\in so\left( 3\right).
$$
We have $\det J=0$, so $m=2n+k$ and $m-k=\mbox{\emph{rk} }J$. Here
$m=3$ and $\mbox{\emph{rk} }J=2$, then $n=k=1$. The system (27)
has beside the energy $H_{1}=H$, a trivial invariant $H_{2}$,
i.e., such that : $J\frac{\partial H_{2}}{\partial x}=0$, or
$$\left(\begin{array}{ccc}
0&-m_{3}&m_{2}\\
m_{3}&0&-m_{1}\\
-m_{2}&m_{1}&0
\end{array}\right)
\left(\begin{array}{c}
\frac{\partial H_{2}}{\partial m_{1}}\\
\frac{\partial H_{2}}{\partial m_{2}}\\
\frac{\partial H_{2}}{\partial m_{3}}
\end{array}\right)
= \left(\begin{array}{c}
0\\
0\\
0
\end{array}\right),
$$
implying $$\frac{\partial H_{2}}{\partial m_{1}} =m_{1}, \quad
\frac{\partial H_{2}}{\partial m_{2}}=m_{2}, \quad \frac{\partial
H_{2}}{\partial m_{3}}=m_{3},$$ and
$$H_{2}=\frac{1}{2}\left(m_{1}^{2}+m_{2}^{2}+m_{3}^{2}\right).$$ The
system evolves on the intersection of the sphere $H_{1}=c_{1}$ and
the ellipsoid $H_{2}=c_{2}$. In $\mathbb{R}^{3}$, this
intersection will be isomorphic to two circles $\left( \text{with
}\frac{c_{2}}{\lambda_{3}}<c_{1}<\frac{c_{2}}{\lambda
_{1}}\right)$. According to Arnold-Liouville's theorem, we have :

\begin{Theo}
The system (27) is completely integrable and the vector
$J\frac{\partial H}{\partial x}$ gives a flow on a variety :
$$\bigcap_{i=1}^{2}\left\{x\in \mathbb{R}^{3}:H_{i}\left(x\right)
=c_{i}\right\}, \quad(\mbox{for generic }c_{i}\in \mathbb{R}),$$
diffeomorphic to a real torus of dimension $1$, that is to say a
circle.
\end{Theo}
Let us now turn to explicit resolution. We shall show that the
problem can be integrated in terms of elliptic functions, as Euler
discovered using his then newly invented theory of elliptic
integrals.  We have just seen that the system in question admits
two first quadratic integrals :
\begin{eqnarray}
H_1&=&\frac{1}{2}\left(\lambda_{1}m_{1}^{2}+\lambda_{2}m_{2}^{2}+\lambda_{3}m_{3}^{2}\right),\nonumber\\
H_2&=&\frac{1}{2}\left(m_{1}^{2}+m_{2}^{2}+m_{3}^{2}\right).\nonumber
\end{eqnarray}
We'll assume that $\lambda_1, \lambda_2, \lambda_3$ are all
different from zero (that is, the solid is not reduced to a point
and is not focused on a straight line either). Under these
conditions, $H_1=0$ implies $m_1=m_2=m_3=0$ and so $H_2=0$; the
solid is at rest. We dismiss this trivial case and now assume that
$H_1\neq 0$ et $H_2\neq 0$. When $\lambda_1=\lambda_2=\lambda_3$,
the equations (27) obviously show that $m_1$, $m_2$ and $m_3$ are
constants. Suppose for example that $\lambda_1=\lambda_2$, the
equations (27) are then written
$$
\dot{m}_{1}=\left(\lambda_{3}-\lambda _{1}\right)
m_{2}m_{3},\qquad \dot{m}_{2}=\left( \lambda _{1}-\lambda
_{3}\right) m_{1}m_{3},\qquad \dot{m}_{3}=0.
$$
We deduce then that $m_3=\mbox{constante}\equiv A$ and
$$
\dot{m}_{1}=A\left(\lambda_{3}-\lambda_{1}\right)m_{2},\qquad
\dot{m}_{2}=A\left(\lambda_{1}-\lambda_{3}\right)m_{1}.
$$
Note that $$(m_1+im_2)^.=iA(\lambda_{1}-\lambda_{3})(m_1+im_2),$$
we obtain $$m_1+im_2=Ce^{iA(\lambda_{1}-\lambda_{3})t},$$ where
$C$ is a constant and so
$$
m_1=C\cos A( \lambda_{1}-\lambda_{3})t,\qquad m_2=C\sin
A(\lambda_{1}-\lambda_{3})t.
$$
The integration of Euler's equations is delicate in the general
case where $\lambda_{1}$, $\lambda_{2}$ and $\lambda_{3}$ are all
different; the solutions are expressed in this case using elliptic
functions. In the following we will suppose that $\lambda_{1}$,
$\lambda_{2}$ and $\lambda_{3}$ are all different and we discard
the other trivial cases which pose no difficulty for solving the
equations in question. To fix the ideas we will assume in the
following that : $\lambda_1>\lambda_2>\lambda_3$. Geometrically,
the equations
\begin{eqnarray}
\lambda_{1}m_{1}^{2}+\lambda_{2}m_{2}^{2}+\lambda_{3}m_{3}^{2}&=&2H_{1},\\
m_{1}^{2}+m_{2}^{2}+m_{3}^{2}&=&2H_{2}\equiv r^2,
\end{eqnarray}
respectively represent the equations of the surface of a half axis
ellipsoid: $\sqrt{\frac{2H_1}{\lambda_1}}$ (half big axis),
$\sqrt{\frac{2H_1}{\lambda_2}}$(middle half axis),
$\sqrt{\frac{2H_1}{\lambda_3}}$(half small axis), and a sphere of
radius $r$. So the movement of the solid takes place on the
intersection of an ellipsoid with a sphere. This intersection
makes sense because by comparing (28) to (29), we see that
$\frac{2H_1}{\lambda_1}<r^2<\frac{2H_1}{\lambda_3}$, which means
geometrically that the radius of the sphere (29) is between the
smallest and largest of the half-axes of the ellipsoid (28). To
study the shape of the intersection curves of the ellipsoid (28)
with the sphere (29), set $H_1>0 $ and let the radius $r$ vary.
Like $\lambda_1>\lambda_2>\lambda_3$, the semi-axes of the
ellipsoid will be
$\frac{2H_1}{\lambda_1}>\frac{2H_1}{\lambda_2}>\frac{2H_1}{\lambda_3}$.
If the radius $r$ of the sphere is less than the half-axis
$\frac{2H_1}{\lambda_3}$ or greater than the half-axis
$\frac{2H_1}{\lambda_1}$, then the intersection in question is
empty (and no real movement corresponds to these values of $H_1$
and $r$). When the radius $r$ equals $\frac{2H_1}{\lambda_3}$,
then the intersection is composed of two points. When the radius
$r$ increases
$\left(\frac{2H_1}{\lambda_3}<r<\frac{2H_1}{\lambda_2}\right)$, we
obtain two curves around the ends of the half minor axis. Likewise
if $r=\frac{2H_1}{\lambda_1}$, we get both ends of the semi-major
axis and if $r$ is slightly smaller than $\frac{2H_1}{\lambda_1}$,
we get two closed curves near these ends. Finally, if
$r=\frac{2H_1}{\lambda_2}$ then the intersection in question
consists of two circles.

\begin{Theo}
Euler's differential equations (27) are integrated by means of
Jacobi's elliptic functions.
\end{Theo}
\emph{Proof}. From the first integrals (28) and (29), we express
$m_1$ and $m_3$ as a function of $m_2$. These expressions are then
introduced into the second equation of the system (27) to obtain a
differential equation in $m_2$ and $\frac{dm_2}{dt}$ only. In more
detail, the following relationships are easily obtained from (28)
and (29) :
\begin{eqnarray}
m_{1}^2&=&\frac{2H_1-r^2\lambda_{3}-\left(\lambda_{2}-\lambda_{3}\right)m_{2}^{2}}{\lambda_{1}-\lambda_{3}},\\
m_{3}^2&=&\frac{r^2\lambda_{1}-2H_1-\left(\lambda_{1}-\lambda_{2}\right)
m_{2}^{2}}{\lambda_{1}-\lambda_{3}}.
\end{eqnarray}
By substituting these expressions in the second equation of the
system (27), we obtain
$$\dot{m}_{2}=\sqrt{(2H_1-r^2\lambda_{3}-\left(\lambda_{2}-\lambda_{3}\right)m_{2}^{2})
(r^2\lambda_{1}-2H_1-\left(\lambda_{1}-\lambda_{2}\right)m_{2}^{2})}.$$
By integrating this equation, we obtain a function $t(m_2)$ in the
form of an elliptic integral. To reduce this to the standard form,
we can assume that $r^2>\frac{2H_1}{\lambda_2}$ (otherwise, it is
enough to invert the indices 1 and 3 in all the previous
formulas). We rewrite the previous equation, in the form
$$\frac{dm_{2}}{\sqrt{(2H_1-r^2\lambda_{3})(r^2\lambda_{1}-2H_1)}dt}=
\sqrt{(1-\frac{\lambda_{2}-\lambda_{3}}{2H_1-r^2\lambda_{3}}m_{2}^{2})
(1-\frac{\lambda_{1}-\lambda_{2}}{r^2\lambda_{1}-2H_1}m_{2}^{2})}.$$
By setting
$$\tau=t\sqrt{(\lambda_{2}-\lambda_{3})(r^2\lambda_{1}-2H_1)},
\quad
s=m_2\sqrt{\frac{\lambda_{2}-\lambda_{3}}{2H_1-r^2\lambda_{3}}},$$
we obtain
$$\frac{ds}{d\tau}=
\sqrt{(1-s^{2})\left(1-\frac{(\lambda_{1}-\lambda_{2})
(2H_1-r^2\lambda_{3})}{(\lambda_2-\lambda_3)(r^2\lambda_{1}-2H_1)}s^{2}\right)},$$
which suggests choosing elliptic functions as a module
$$k^2=\frac{(\lambda_{1}-\lambda_{2})(2H_1-r^2\lambda_{3})}{(\lambda_2-\lambda_3)
(r^2\lambda_{1}-2H_1)}.$$ Inequalities
$\lambda_1>\lambda_2>\lambda_3$,
$\frac{2H_1}{\lambda_1}<r^2<\frac{2H_1}{\lambda_3}$ and
$r^2>\frac{2H_1}{\lambda_2}$ show that $0<k^2<1$. So we get
$$\frac{ds}{d\tau}=\sqrt{(1-s^{2})(1-k^2s^{2})}.$$ This equation
admits the solution (we choose the origin of the times such that
$m_2=0$ for $t=0$):
$$\tau=\int_0^s\frac{ds}{\sqrt{(1-s^2)(1-k^2s^2)}}.$$
It is the integral of a holomorphic differential on an elliptic
curve :
$$
\mathcal{E}: w^2=(1-s^2)(1-k^2s^2).
$$
The inverse function $s(\tau)$ is one of Jacobi's elliptic
functions : $s=\mathbf{sn}\tau$, which determines $m_2$ as,
$$m_2=\sqrt{\frac{2H_1-r^2\lambda_{3}}{\lambda_2-\lambda_3}}\cdot\mathbf{sn}\tau.$$
According to the equalities (30) and (31), we know that the
functions $m_1$ and $m_3$  are expressed algebraically as a
function of $m_2$, so
$$m_1=\sqrt{\frac{2H_1-r^2\lambda_{3}}{\lambda_1-\lambda_3}}
\cdot\sqrt{1-\mathbf{sn}^2\tau},$$
$$m_3=\sqrt{\frac{r^2\lambda_{1}-2H_1}{\lambda_1-\lambda_3}}\cdot
\sqrt{1-k^2\mathbf{sn}^2\tau}.$$ Given the definition of the other
two elliptical functions :
$$\mathbf{cn}\tau=\sqrt{1-\mathbf{sn}^2\tau}, \qquad\mathbf{dn}\tau=
\sqrt{1-k^2\mathbf{sn}^2\tau},$$ and the fact that
$$\tau=t\sqrt{(\lambda_{2}-\lambda_{3})(r^2\lambda_{1}-2H_1)},$$ we
finally get the following explicit formulas :
\begin{eqnarray}
m_1&=&\sqrt{\frac{2H_1-r^2\lambda_{3}}{\lambda_1-\lambda_3}}\mathbf{cn}(t\sqrt{(\lambda
_{2}-\lambda _{3})(r^2\lambda
_{1}-2H_1)}),\nonumber\\
m_2&=&\sqrt{\frac{2H_1-r^2\lambda_{3}}{\lambda_2-\lambda_3}}\mathbf{sn}(t\sqrt{(\lambda
_{2}-\lambda _{3})(r^2\lambda
_{1}-2H_1)}),\\
m_3&=&\sqrt{\frac{r^2\lambda_{1}-2H_1}{\lambda_1-\lambda_3}}\mathbf{dn}(t\sqrt{(\lambda
_{2}-\lambda _{3})(r^2\lambda _{1}-2H_1)}).\nonumber
\end{eqnarray}
In other words, the integration of the Euler equations is done by
means of elliptic Jacobi functions and the proof is complete.

\begin{Rem}
Note that for $\lambda_1=\lambda_2$, we have $k^2=0$. In this
case, the elliptical functions $\mathbf{sn}\tau, \mathbf{cn}\tau,
\mathbf{dn}\tau$ are reduced respectively to functions $\sin\tau,
\cos\tau, 1$. From the system (32), we easily obtain the
expressions
\begin{eqnarray}
m_1&=&\sqrt{\frac{2H_1-r^2\lambda_{3}}{\lambda_1-\lambda_3}}
\cos\sqrt{(\lambda_1-\lambda_3)(r^2\lambda_{1}-2H_1)}t,\nonumber\\
m_2&=&\sqrt{\frac{2H_1-r^2\lambda_{3}}{\lambda_1-\lambda_3}}
\sin\sqrt{(\lambda_1-\lambda_3)(r^2\lambda_{1}-2H_1)}t,\nonumber\\
m_3&=&\sqrt{\frac{r^2\lambda_{1}-2H_1}{\lambda_1-\lambda_3}}.\nonumber
\end{eqnarray}
We find the solutions established previously where the constants
$A$ and $B$ are
$A=\sqrt{\frac{r^2\lambda_{1}-2H_1}{\lambda_1-\lambda_3}}$ and
$C=\sqrt{\frac{2H_1-r^2\lambda_{3}}{\lambda_1-\lambda_3}}$.
\end{Rem}

In the case of the Lagrange top, we have $I_{1}=I_{2}$,
$l_{1}=l_{2}=0$, i.e., the Lagrange top [21] is a rigid body, in
which two moments of inertia are the same and the center of
gravity lies on the symmetry axis. In other words, the Lagrange
top is a symmetric top with a constant vertical gravitational
force acting on its center of mass and leaving the base point of
its body symmetry axis fixed. As in the case of Euler, we show
that in this case also the problem is solved by elliptic
integrals. Or what amounts to the same, the integration is done
using elliptic functions.

The Kowalewski top [20] is special symmetric top with a unique
ratio of the moments of inertia satisfy the relation :
$I_{1}=I_{2}=2I_3$, $l_{3}=0$; in which two moments of inertia are
equal, the third is half as large, and the center of gravity is
located in the plane perpendicular to the symmetry axis (parallel
to the plane of the two equal points). Moreover, we may choose
$l_2=0$, $\mu gl_1=l$ and $I_3=1$. After the substitution
$t\rightarrow 2t$ the system (26) is written explicitly in the
form
\begin{eqnarray}
\overset{.}{m}_{1}&=&m_{2}m_{3},\nonumber\\
\overset{.}{m}_{2}&=&-m_{1}m_{3}+2\gamma_{3},\nonumber\\
\overset{.}{m}_{3}&=&-2\gamma_{2},\\
\overset{.}{\gamma}_{1}&=&2m_{3}\gamma_{2}-m_{2}\gamma_{3},\nonumber\\
\overset{.}{\gamma}_{2}&=&m_{1}\gamma_{3}-2m_{3}\gamma_{1},\nonumber\\
\overset{.}{\gamma}_{3}&=&m_{2}\gamma_{1}-m_{1}\gamma_{2}.\nonumber
\end{eqnarray}
These equations are written in the form (see example 12) of a
Hamiltonian vector field $$\dot{x}=J\frac{\partial H}{\partial x},
\quad
x=(m_{1},m_{2},m_{3},\gamma_{1},\gamma_{2},\gamma_{3})^{\top},$$
where
$$H=\frac{1}{2}\left( m_{1}^{2}+m_{2}^{2}\right) +m_{3}^{2}+2\gamma_{1},$$
is the Hamiltonian and
$$
J=\left(\begin{array}{cccccc}
0&-m_{3}&m_{2}&0&-\gamma_{3}&\gamma _{2}\\
m_{3}&0&-m_{1}&\gamma_{3}&0&-\gamma _{1}\\
-m_{2}&m_{1}&0&-\gamma _{2}&\gamma _{1}&0\\
0&-\gamma_{3}&\gamma_{2}&0&0&0\\
\gamma_{3}&0&-\gamma_{1}&0&0&0\\
-\gamma _{2}&\gamma_{1}&0&0&0&0
\end{array}\right).
$$
The above system admits four first integrals :
\begin{eqnarray}
H_{1}&\equiv& H,\nonumber\\
H_{2}&=&m_{1}\gamma_{1}+m_{2}\gamma_{2}+m_{3}\gamma_{3},\\
H_{3}&=&\gamma _{1}^{2}+\gamma_{2}^{2}+\gamma_{3}^{2},\nonumber\\
H_{4}&=&\left(\left(\frac{m_{1}+im_{2}}{2}\right)^{2}-\left(\gamma_{1}+
i\gamma_{2}\right)\right)\left(\left(\frac{m_{1}-im_{2}}{2}\right)^{2}
-\left(\gamma_{1}-i\gamma_{2}\right)\right).\nonumber
\end{eqnarray}
A second flow commuting with the first flow is regulated by the
equations : $$\dot{x}=J\frac{\partial H_{4}}{\partial x}, \quad
x=(m_{1},m_{2},m_{3},\gamma_{1},\gamma _{2},\gamma_{3})^{\top}.$$
The first integrals $H_{1}$ and $H_{4}$ are in involution,
$$\{H_1,H_4\}=\left\langle\frac{\partial H_1}{\partial
x},J\frac{\partial H_4}{\partial x}\right\rangle=0,$$ while
$H_{2}$ and $H_{3}$ are trivial, $$J\frac{\partial H_2}{\partial
x}=J\frac{\partial H_3}{\partial x}=0.$$ Let $\mathcal{A}$ be the
complex affine variety defined by the intersection of the
constants of the motion
\begin{equation}\label{eqn:euler}
\mathcal{A}=\bigcap_{k=1}^{4}\left\{x:H_{k}(x)=c_{k}\right\},
\end{equation}
where $c=(c_1,c_2,c_3=1,c_4)$ is not a critical value. We will
explain how the affine variety $\mathcal{A}$ and vector-fields
behave after the quotient by some natural involution on
$\mathcal{A}$ and how these vector-fields become well defined when
we take Kowalewski's variables. We show that these variables are
naturally related to the so-called Euler's differential equations
and can be seen as the addition-formula for the Weierstrass
elliptic function. In the theorem below (for further information,
see also [20, 22]), we will use with Kowalewski the following
notations : $c_1=6h_1$, $c_2=2h_2$ et $c_4=k^2$.

\begin{Theo}
$a)$ Let $$(m_{1},m_{2},m_{3},\gamma_{1},\gamma_{2},\gamma_{3})
\longmapsto(x_{1},x_{2},m_{3},y_{1},y_{2},\gamma_{3}),$$ be a
birationally map on the variety $\mathcal{A}$(35)  where $x_1$,
$x_2$, $y_1$, $y_2$ are defined as
\begin{eqnarray}
x_{1}&=&\frac{1}{2}(m_{1}+im_{2}),\nonumber\\
x_{2}&=&\frac{1}{2}(m_{1}-im_{2}),\\
y_{1}&=&x_{1}^{2}-\left(\gamma_{1}+i\gamma_{2}\right),\nonumber\\
y_{2}&=&x_{2}^{2}-\left(\gamma_{1}-i\gamma_{2}\right).\nonumber
\end{eqnarray}
Then, the quotient $K\equiv\mathcal{A}/\sigma$ by the involution
\begin{equation}\label{eqn:euler}
\sigma : M_{c}\longrightarrow
M_{c}\text{}\left(x_{1},x_{2},m_{3},y_{1},y_{2},\gamma_{3}\right)
\longmapsto\left(x_{1},x_{2},-m_{3},y_{1},y_{2},-\gamma_{3}\right),
\end{equation}
is a Kummer surface
\begin{equation}\label{eqn:euler}
K : \left\{\begin{array}{rl}
&y_{1}y_{2}=k^2,\\
&y_{1}R(x_{2})+y_{2}R(x_{1})+R_{1}\left( x_{1},x_{2}\right)
+k^2(x_{1}-x_{2})^{2}=0,
\end{array}\right.
\end{equation}
where
\begin{equation}\label{eqn:euler}
R(x)=-x^{4}+6h_1x^{2}-4h_{2}x+1-k^2,
\end{equation}
is a polynomial of degree $4$ in $x$ and
\begin{eqnarray}
R_{1}(x_{1},x_{2})&=&-6h_1x_{1}^{2}x_{2}^{2}+
4h_{2}x_{1}x_{2}\left(
x_{1}+x_{2}\right)\\
&&-\left(1-k^2\right) \left(
x_{1}+x_{2}\right)^{2}+6h_1\left(1-k^2\right)-4h_2^2,\nonumber
\end{eqnarray}
is another polynomial of degree $2$ in $x_{1}$, $x_{2}$. The
ramification points of $\mathcal{A}$ on $K$ are given by the $8$
fixed points of the involution $\sigma$.

$b)$ The surface $K$ is a double cover of plane $(x_{1}, x_{2})$,
ramified along two elliptic curves intersecting exactly each other
at the $8$ fixed points of the involution $\sigma$. These curves
give rise to the Euler differential equation\index{Euler
differential equation}
$$\frac{\dot x_1}{\sqrt{R(x_1)}}\pm \frac{\dot
x_2}{\sqrt{R(x_2)}}=0,$$ to which are connected the variables of
Kowalewski
$$s_{1}=\frac{R\left(x_{1},x_{2}\right)
-\sqrt{R(x_{1})}\sqrt{R(x_{2})}}{\left( x_{1}-x_{2}\right)
^{2}}+3h_1,$$
$$s_{2}=\frac{R\left( x_{1},x_{2}\right)
+\sqrt{R(x_{1})}\sqrt{R(x_{2})}}{\left( x_{1}-x_{2}\right)
^{2}}+3h_1,$$ where
\begin{equation}\label{eqn:euler}
R(x_1,x_2)\equiv-x_1^2x_2^2+6h_1x_1x_2-2h_{2}(x_1+x_2)+1-k^2,
\end{equation}
and can be seen as addition formulas for the Weierstrass elliptic
function.

$c)$ In terms of the variables $s_{1}$ and $s_{2}$, the system of
differential equations  (33) is reduced to the system
\begin{eqnarray}
\frac{ds_{1}}{\sqrt{P_{5}(s_{1})}}\pm
\frac{ds_{2}}{\sqrt{P_{5}(s_{2})}}&=&0,\nonumber\\
\frac{s_{1}ds_{1}}{\sqrt{P_{5}(s_{1})}}\pm
\frac{s_{2}ds_{2}}{\sqrt{P_{5}(s_{2})}}&=&dt,\nonumber
\end{eqnarray}
where
$P_{5}(s)$ is a fifth-degree polynomial and the problem can be
integrated in terms of genus two hyperelliptic functions.
\end{Theo}
\emph{Proof}.  $a)$ Using the change of variables (36), with
$t\rightarrow it$, equations (33) and (34) become
\begin{eqnarray}
\overset{.}{x}_{1}&=&m_{3}x_{1}-\gamma_3,\nonumber\\
\overset{.}{x}_{2}&=&-m_{3}x_{2}+\gamma_{3},\nonumber\\
\overset{.}{m}_{3}&=&-x_1^{2}+y_1+x_2^2-y_2,\\
\overset{.}{y}_{1}&=&2m_{3}y _{1},\nonumber\\
\overset{.}{y}_{2}&=&-2m_{3}y _{1},\nonumber\\
\overset{.}{\gamma}_{3}&=&x_{1}(x_2^2-y_2)-x_{2}(x_1^2-y_1),\nonumber
\end{eqnarray}
and
\begin{eqnarray}
y_1y_2&=&k^2,\nonumber\\
m_3^2&=&6h_1+y_1+y_2-(x_1+x_2)^2,\\
m_3\gamma_3&=&2h_{2}+x_1y_2+x_2y_1-x_1x_2(x_1+x_2),\nonumber\\
\gamma_3^2&=&1-k^2+x_1^2y_2+x_2^2y_1-x_1^2x_2^2.\nonumber
\end{eqnarray}
It's obvious that $\sigma$(37) is an automorphism of
$\mathcal{A}$, of order two. The quotient $\mathcal{A}/\sigma$ by
the involution $\sigma$ is a Kummer $K$ defined by (38). The
variety $\mathcal{A}$ is a double cover of the surface $K$
branched over the fixed points of the involution $\sigma$. To find
them, we substitute $m_3=\gamma_3=0$ in the system (34), to wit
\begin{eqnarray}
y_1y_2&=&k^2,\\
y_1+y_2&=&(x_1+x_2)^2-6h_1,\\
x_2y_1+x_1y_2&=&x_1x_2(x_1+x_2)-2h_{2},\\
x_2^2y_1+x_1^2y_2&=&x_1^2x_2^2+k^2-1.
\end{eqnarray}
Away from the $x_1^2=x_2^2$, we may solve (45) and (47) in $y_1$
and $y_2$ and substitute into the equations (44) and (46); one
then finds two curves in $x_1$ and $x_2$ whose equations are
\begin{eqnarray}
R(x_1,x_2)&\equiv&-x_1^2x_2^2+6h_1x_1x_2-2h_{2}(x_1+x_2)+1-k^2=0,\nonumber\\
S(x_1,x_2)&\equiv&
\left(x_1^4+2x_1^3x_2-6h_1x_1^2+1-k^2\right)\left(x_2^4+2x_1x_2^3-6h_1x_2^2+1-k^2\right)\nonumber\\
&&+k^2(x_1^2-x_2^2)^2=0.\nonumber
\end{eqnarray}
These curves intersect at the zeroes of the resultant
$\mbox{Res}(R,S)$ of $R$, $S$ :
\begin{equation}\label{eqn:euler}
\mbox{Res}(R,S)_{x_2}=x_{1}^2\left(x_{1}^{4}+6h_{1}x_1^{2}+k^2-1\right)^{2}P_8(x_1),
\end{equation}
where $P_8(x_1)$ is a monic polynomial of degree $8$. Since the
root $x_1$ must be excluded (it indeed implies that the leading
terms of $R$ and $S$ vanish), the possible intersections of the
curve R and S will be,

$(i)$ at the roots of
$$x_{1}^{4}+6h_{1}x_1^{2}+k^2-1=0,$$ this is unacceptable, because
then one checks that the common roots of $R$ and $S$ would have
the property that $x_1^2=x_2^2$, which was excluded.

$(ii)$ at the roots of $P_8(x_1)=0$; there, for generic $k$ and
$h$, $x_1^2\neq x_2^2$.\\ Finally, we must analyze the case
$x_1^2=x_2^2$ for which one checks that (44),...,(47) has no
common roots. Consequently the involution $\sigma$ has $8$ fixed
points on the affine variety $\mathcal{A}$. Clearly the vector
field (42) vanishes at the fixed points of the involution
$\sigma$.

$b)$ From equations (38), we deduce
$$ y_1=\frac{-1}{2R(x_2)}\left(R_1(x_1,x_2)+k^2(x_1-x_2)^2+\Delta\right),$$
$$ y_2=\frac{-1}{2R(x_1)}\left(R_1(x_1,x_2)+k^2(x_1-x_2)^2-\Delta\right),$$
where
$$\Delta^2=\left(R_1(x_1,x_2)+k^2(x_1-x_2)^2\right)^2-4k^2R(x_1)R(x_2)\equiv
P(x_1,x_2).$$ Therefore, the surface $K$ is a double cover of
$\mathbb{C}^2$, ramified along the curve $\mathcal{C}:
P(x_1,x_2)=0$. This equation is reducible and can be written as
the product $$P(x_1,x_2)=P_1(x_1,x_2).P_2(x_1,x_2),$$ of two
symmetric polynomials (in $x_1, x_2$) of degree two in each one of
the variables $x_1, x_2$, i.e.,
$$
P_1(x_1,x_2)=a(x_1)x_2^2+2b(x_1)x_2-c(x_1)=a(x_2)x_1^2+2b(x_2)x_1-c(x_2),
$$
where
\begin{eqnarray}
a(x)&=&-2(k+3h_1)x^2+4h_2x-1,\nonumber\\
b(x)&=&2h_2x^2+(2k(k+3h_1)-1)x-2h_2k,\nonumber\\
c(x)&=&x^2+4h_2kx+2(k^2-1)(k+3h_1)+4h_2^2,\nonumber
\end{eqnarray}
while the polynomial $P_2(x_1, x_2)$ is obtained from $P_1(x_1,
x_2)$ after replacing $k$ with $-k$. Note that the curve
$\mathcal{C}_1 : P_1(x_1,x_2)=0$, is elliptic :
\begin{eqnarray}
x_1&=&\frac{-b(x_2)\pm\sqrt{2(k+3h_1)-4h_2^2}\sqrt{R(x_2)}}{a(x_2)},\nonumber\\
x_2&=&\frac{-b(x_1)\pm\sqrt{2(k+3h_1)-4h_2^2}\sqrt{R(x_1)}}{a(x_1)},\nonumber
\end{eqnarray}
where $R(x)$ is given by (39). Similarly, the curve $\mathcal{C}_2
: P_2(x_1, x_2)=0$, is elliptic and we notice that the two curves
$\mathcal{C}_1$ and $\mathcal{C}_2$ intersect exactly at the $8$
fixed points of involution $\sigma$, because (see (48)) :
$$\mbox{Res}(P_1, P_2)_{x_2}=16k^{2}P_8(x_1).$$ Differentiating the
symmetric equation $P_1(x_1,x_2)=0$ (or $P_2(x_1,x_2)=0$) with
regard to $t$, one finds
$$\pm
2\sqrt{2(k+3h_1)-4h_2^2}\sqrt{R(x_2)}\dot x_1\pm
2\sqrt{2(k+3h_1)-4h_2^2}\sqrt{R(x_1)}\dot x_2 =0.$$ Hence
\begin{equation}\label{eqn:euler}
\frac{\dot x_1}{\sqrt{R(x_1)}}\pm \frac{\dot
x_2}{\sqrt{R(x_2)}}=0.
\end{equation}
Since $R(x_1)$ and $R(x_2)$ are two polynomials of the fourth
degree in $x_1$ and $x_2$ respectively and having the same
coefficients, then (49) is the so-called Euler's equation. The
reader is referred to Halphen [15] and Weil [38] for this theory
that we summarize here as follows : let
$$F(x)=a_0x^4+4a_lx^3+6a_2x^2+4a_3x+a_4,$$ be a polynomial of the
fourth degree. The general integral of Euler's equation
$$\frac{\dot x}{\sqrt{F(x)}}\pm \frac{\dot y}{\sqrt{F(y)}}=0,$$ can
be written in two different ways :
$$F_1(x, y)+2sF(x, y)-s^2(x-y)^2=0,$$
where
$$F(x,y)=a_0x^2y^2+2a_1xy(x+y)+3a_2(x^2+y^2)+2a_3(x+y)+a_4,$$ and
$$F_1(x,y)=\frac{F(x)F(y)-F^2(x,y)}{(x-y)^2},$$ or in an irrational
form $$\frac{F(x,y)\mp \sqrt{F(x)}\sqrt{F(y)}}{(x-y)^2}=s,$$ which
can be seen as the addition-formula for the Weierstrass elliptic
function
$$2\wp
(u+v)=\frac{(\wp(u)+\wp(v))\left(2\wp(u)\wp(v)-\frac{1}{2}g_2\right)-g_3-\wp'(u)\wp'(v)}
{\left(\wp(u)+\wp(v)\right)^2},$$
$$\wp'^2(u)=\left(\frac{d\wp}{du}\right)^2=4\wp^3-g_2\wp-g_3,$$
$\wp(u)=x$,  $\wp(v)=y$, $F(x)=4x^3-g_2x-g_3$, $\wp'^2(u)=F(x)$,
$\wp'^2(v)=F(y)$, $2\wp(u+v)=s$ and $g_2$, $g_3$ are constants. We
apply these facts to Kowalewski's problem with $F(x)=R(x)$,
$F(x_1,x_2)=R(x_1,x_2)+3h_1(x_1-x_2)^2$, and $a_0=-1$, $a_1=0$,
$a_2=h_1$, $a_3=-h_2$, $a_4=1-k^2$and $s=k+3h_1$. So the
polynomial $P_1(x_1,x_2)$ which can also be regarded as a solution
of (49), can also be written as
$$R_1(x_1,x_2)+2sR(x_1,x_2)-s^2(x_1-x_2)^2=0,$$
where $R_1(x_1,x_2)$ is given by (40) and has the form
$$R_1(x_1,x_2)=\frac{R(x_1)R(x_2)-R^2(x_1,x_2)}{(x_1-x_2)^2}.$$
Remember that $R(x_1,x_2)$ is given by (41). The solution of (49)
can also be expressed
\begin{equation}\label{eqn:euler}
\frac{R(x_1,x_2)\mp
\sqrt{R(x_1)}\sqrt{R(x_2)}}{(x_1-x_2)^2}+3h_1=s.
\end{equation}

$c)$ Let us carry out the calculations, assuming the polynomial
$R(x)$ reduced to the form $4x^3-g_2x-g_3$ and call $s_1$ (resp.
$s_2$) the relation (50) with the sign -(resp. +). Now, outside
the branch locus of $K$(38) over $\mathbb{C}^2$, the equation (49)
is not identically zero and may be written in the form
\begin{eqnarray}
\frac{\dot x_1}{\sqrt{R(x_1)}}+\frac{\dot
x_2}{\sqrt{R(x_2)}}&=&\frac{\dot
s_1}{\sqrt{4s_1^3-g_2s_1-g_3}}\neq 0,\\
\frac{\dot x_1}{\sqrt{R(x_1)}}-\frac{\dot
x_2}{\sqrt{R(x_2)}}&=&\frac{\dot
s_2}{\sqrt{4s_2^3-g_2s_2-g_3}}\neq 0.\nonumber
\end{eqnarray}
where $g_2=k^2-1+3h_1^2$ and $g_3=h_1(k^2-1-h_1^2)+h_2^2$. After
some algebraic manipulation we deduce from (43),
\begin{eqnarray}
(m_3x_1-\gamma_3)^2&=&R(x_1)+(x_1-x_2)^2y_1,\nonumber\\
(m_3x_2-\gamma_3)^2&=&R(x_2)+(x_1-x_2)^2y_2,\nonumber\\
(m_3x_1-\gamma_3)(m_3x_2-\gamma_3)&=&R(x_1,x_2),\nonumber
\end{eqnarray}
and from (42),
\begin{eqnarray}
\dot{x}_1^2&=&R(x_1)+(x_1-x_2)^2y_1,\nonumber\\
\dot{x}_2^2&=&R(x_2)+(x_1-x_2)^2y_2.\nonumber
\end{eqnarray}
This together with (38) and (51) implies that
\begin{eqnarray}
\frac{\dot{s}_1^2}{4s_1^3-g_2s_1-g_3}&=&\left(\frac{\dot
x_1}{\sqrt{R(x_1)}}+\frac{\dot
x_2}{\sqrt{R(x_2)}}\right)^2,\nonumber\\
&=&\frac{(x_1-x_2)^4}{R(x_1)R(x_2)}
\left[\left(\frac{R(x_1,x_2)-\sqrt{R(x_1)}\sqrt{R(x_2)}}{(x_1-x_2)^2}\right)-k^2\right],\nonumber\\
&=&4\frac{(s_1-3h_1)^2-k^2}{(s_1-s_2)}.\nonumber
\end{eqnarray}
In the same way, we find
$$
\frac{\dot{s}_2^2}{4s_2^3-g_2s_2-g_3}=\left(\frac{\dot
x_1}{\sqrt{R(x_1)}}-\frac{\dot x_2}{\sqrt{R(x_2)}}\right)^2
=4\frac{(s_2-3h_1)^2-k^2}{(s_2-s_1)}.
$$
In terms of the variables $s_1$ and $s_2$, the system (33) becomes
\begin{eqnarray}
\frac{\dot{s}_1}{\sqrt{P(s_1)}}+\frac{\dot{s}_2}{\sqrt{P(s_2)}}&=&0,\nonumber\\
\frac{s_1\dot{s}_1}{\sqrt{P(s_1)}}+\frac{s_2\dot{s}_2}{\sqrt{P(s_2)}}&=&i,\nonumber
\end{eqnarray}
where
$$P_5(s)=\left((s-3h_1)^2-k^2\right)\left(4s^3-g_2s-g_3\right),$$ is
a polynomial of degree $5$. As known, such integrals are called
hyperelliptic integrals and the problem can be integrated in terms
of genus two hyperelliptic functions of time. More precisely,
these equations are integrable by the transformation of Abel
$$\mathcal{H}\longrightarrow Jac(\mathcal{H})=\mathbb{C}^{2}/\Lambda ,
\text{ }p\longmapsto \left(\int_{p_{0}}^{p}\theta_{1},
\int_{p_{0}}^{p}\theta_{2}\right),$$ where $\mathcal{H}$ is the
hyperelliptic curve of genus $2$ associated with equation :
$w^{2}=P_{5}(s)$, $\Lambda$ is the lattice generated by the
vectors $n_{1}+\Omega_{\mathcal{H}}n_{2},\left(
n_{1},n_{2}\right)\in\mathbb{Z}^{2}$, $\Omega_{\mathcal{H}}$ is
the matrix of periods of $\mathcal{H}$, $(\theta_{1},\theta_{2})$
is a basis of holomorphic differentials on $\mathcal{H}$, i.e.,
$$\theta_{1}=\frac{ds}{\sqrt{P_{5}(s)}},
\qquad \theta_{2}=\frac{sds}{\sqrt{P_{5}(s)}},$$ and $p_{0}$ is a
fixed point on $\mathcal{H}$. The theorem is thus proved.
$\square$

We mentioned previously that any algebraic integral of equations
(24) is a combination of classical integrals except in the cases
of Euler, Lagrange and Kowalewski and that there could therefore
be no first algebraic integral other than those highlighted in
these three cases. In addition, there are a few special cases :

- The case of Hesse-Appel'rot [16, 3] :
$$l_2=0,\quad l_1\sqrt{I_1(I_2-I_3)}+l_3\sqrt{I_3(I_1-I_2)}=0.$$
In this case, equation $l_1m_1+l_3m_3=0$ represents a particular
first integral obtained by Hesse and the integration is carried
out using elliptic functions.

- The case of Goryachev-Chaplygin [12, 8]: $I_1=I_2=4I_3$,
$l_2=l_3=0$. In this case, the system (22) admits the first
integral
$$\lambda_3m_3(\lambda_1^2+m_1^2+\lambda_2^2m_2^2)+\mu
gl_1\lambda_1\lambda_3m_1\gamma_3=g, \quad
\lambda_i=\frac{1}{I_i}, i=1,2,3$$ and integration is carried out
using hyperelliptic functions of genus $2$.

- The case of Bobylev-Steklov [6, 36]: $I_2=2I_1$, $l_1=l_3=0$.
The integration of the equations in this case is easy, using
elliptic functions.

\subsection{Yang-Mills field with gauge group $SU(2)$}

We begin by introducing some notions related to Yang-Mills field
[44] with $SU(2)$ as gauge group. For general considerations and
the details of certain notions, one can consult for example [9].
Consider the special unitary group $SU(2)$ of degree $2$, i.e.,
the set of $2\times 2$ unitary matrices with determinant $1$. This
is a real Lie group of dimension three. It is compact, simply
connected, simple and semi-simple. The group $SU(2)$ is isomorphic
to the group of quaternions of norm one and is diffeomorphic to
the $3$-sphere $S^3$. It is well known that the quaternions
represent the rotations in $3$-dimensional space and hence there
exists a surjective homomorphism of $SU(2)$ on the rotation group
$SO (3)$ whose kernel is $\{+I, -I\}$ (the identical application
and its opposite). Recall also that $SU(2)$ is identical to one of
the symmetry spinor groups, $Spin(3)$, that enables a spinor
presentation of rotations. The Lie algebra $su(2)$ corresponding
to $SU(2)$ consists of the $2\times 2$ antihermitian complex
matrices with null trace, the standard commutator serving as a Lie
bracket. It is a real algebra. The algebra $su(2)$ is isomorphic
to the Lie algebra $so(3)$.

We consider the Yang-Mills field $F_{kl}$ as a vector field with
values in the algebra $su (2)$. It is a local expression of the
gauge field or connection defining the covariant derivative of
$F_{kl}$ in the adjoint representation of $su(2)$. To determine
this expression, note that each Lorentz component of the
Yang-Mills field develops on a basis ($\sigma_1, \sigma_2,
\sigma_3$) de $su(2)$, $A_k=A_k^\alpha \sigma_\alpha$,
$\alpha=1,2,3$, $k=1,2,3,4$, where the $\sigma_\alpha$ are the
matrices of Pauli
$$
{\sigma_1}=\left(\begin{array}{cc}
0&-i\\
i&0
\end{array}\right),\quad
{\sigma_2}=\left(\begin{array}{cc}
0&1\\
1&0
\end{array}\right),\quad
{\sigma_3}=\left(\begin{array}{cc}
1&0\\
0&-1
\end{array}\right).
$$
These matrices are often used in quantum mechanics to represent
the spin of particles. Moreover, the group $SU(2)$ is associated
with gauge symmetry in the description of the weak or weak force
interaction (one of the four fundamental forces of nature) and is
therefore of particular importance in the physics of particles. In
fact, it was only in the late 1960s that the importance of the
Yang-Mills equations became apparent, especially when the concept
of the gauge fields was defined as the one of the four fundamental
physical interactions (gravitational, electromagnetic, weak and
strong interactions). The dynamics of the Yang-Mills theory is
determined by the Lagrangian density $$\mathcal{L}=-\frac{1}{2}{
Tr}\{F_{kl}F^{kl}\}, \quad 1\leq k,l\leq 4$$ where
$$F_{kl}=\frac{\partial A_{l}}{\partial\tau_{k}}-\frac{\partial
A_{k}}{\partial\tau_{l}}+\left[A_{k},A_{l}\right],$$ is the
expression of the anti-symmetric Faraday tensors with values in
$su(2)$. These tensors are not invariant under gauge
transformations. On the other hand, we verify that
${Tr}\{F_{kl}F^{kl}\}$ is actually gauge invariant. The trace
relates to the internal space $su(2)$. The equations of the motion
are given by
$$D _{k}F^{kl}=\frac{\partial F^{kl}}{\partial \tau _{k}}+
\left[ A_{k},F^{kl}\right]=0,\qquad F_{kl},A_{k}\in su(2),\quad
1\leq k,l\leq 4,$$ with $D_k$ the covariant derivative in the
adjoint representation of the algebra $su(2)$ and in which
$\left[A_{k},F^{kl}\right]$ is the crochet of the two fields in
$su(2)$. The Yang-Mills theory extends the principle of gauge
invariance of electromagnetism to other groups of continuous Lie
transformations. Thus the $F_{kl}$ tensor generalizes the
electromagnetic field and the Yang-Mills equations are the
non-commutative generalization of the Maxwell equations. The
self-dual Yang-Mills (SDYM) equations is an universal system for
which some reductions include all classical tops from Euler to
Kowalewski (0+1-dimensions), KdV, Nonlinear Schr\"{o}dinger,
Sine-Gordon, Toda lattice and N-waves equations (1+1-dimensions),
KP and D-S equations (2+1-dimensions). There is a vast literature
devoted to the study of the above equations (see monograph [35]
for many references concerning both theoretical and practical
results).

We are interested here in the field of homogeneous
double-component field. In this case, we have $\frac{\partial
A_{l}}{\partial\tau_{k}}=0, (k\neq1)$, $A_{1}=A_{2}=0$,
$A_{3}=n_{1}U_{1}\in su(2)$, $A_{4}=n_{2}U_{2}\in su(2)$, where
$n_{i}$ are $su(2) $-generators (i.e., they satisfy commutation
relations : $n_1=[n_2,[n_1,n_2]], n_2=[n_1,[n_2,n_1]]$). The
system becomes
\begin{eqnarray}
\frac{\partial {2}U_{1}}{\partial
t^{2}}+U_{1}U_{2}^{2}&=&0,\nonumber\\ \frac{\partial
^{2}U_{2}}{\partial t^{2}}+U_{2}U_{1}^{2}&=&0,\nonumber
\end{eqnarray}
with $t=\tau _{1}$. By setting $U_{1}=q_{1}$, $U_{2}=q_{2}$,
$\frac{\partial U_{1}}{\partial t}=p_{1}$, $\frac{\partial
U_{2}}{\partial t}=p_{2}$, Yang-Mills equations are reduced to
Hamiltonian system $$\dot{x}=J\frac{\partial H}{\partial x}, \quad
x=(q_{1},q_{2},p_{1},p_{2})^{\intercal}, \quad
J=\left(\begin{array}{cc}
O&-I\\
I&O
\end{array}\right),$$
where $$H=\frac{1}{2}\left(
p_{1}^{2}+p_{2}^{2}+q_{1}^{2}q_{2}^{2}\right),$$ is the
Hamiltonian. Note that the symplectic transformation :
$$p_{1}=\frac{\sqrt{2}}{2}\left(x_{1}+x_{2}\right),
\qquad p_{2}=\frac{\sqrt{2}}{2}\left(x_{1}-x_{2}\right),$$
$$q_{1}=\frac{1}{2}\left(\root{4}\of{2}\right)^{3}\left(
y_{1}+iy_{2}\right), \qquad q_{2}=\frac{1}{2}\left(
\root{4}\of{2}\right)^{3}\left(y_{1}-iy_{2}\right),$$ takes this
Hamiltonian into
$$H=\frac{1}{2}\left( x_{1}^{2}+x_{2}^{2}\right)
+\frac{1}{4}\left(y_{1}^{2}+y_{2}^{2}\right)^{2}.$$ The
Hamiltonian dynamical system associated with $H$ is written
\begin{eqnarray}
\dot{y}_{1}&=&x_{1},\nonumber\\
\dot{y}_{2}&=&x_{2},\\
\dot{x}_{1}&=&-\left( y_{1}^{2}+y_{2}^{2}\right)y_1,\nonumber\\
\dot{x}_{2}&=&-\left(y_{1}^{2}+y_{2}^{2}\right)y_2.\nonumber
\end{eqnarray}
These equations give a vector field on $\mathbb{R}^{4}$. The
existence of a second independent first integral in involution
with $H_{1}\equiv H$, is enough for the system to be completely
integrable. The above differential system implies
\begin{eqnarray}
\ddot{y}_1+\left( y_{1}^{2}+y_{2}^{2}\right)y_1&=&0,\nonumber\\
\ddot{y}_2+\left( y_{1}^{2}+y_{2}^{2}\right)y_2&=&0.\nonumber
\end{eqnarray}
Obviously, the moment : $H_2=x_{1}y_{2}-x_{2}y_{1}$, is a first
integral, $H_1$ and $H_2$ are in involution $\{H_1,H_2\}=0$ and
$H_2$ determines with $H_1$ an integrable system. Let
$$\mathcal{M}_c=\left\{x\equiv (y_1,y_2,x_1,x_2)\in
\mathbb{R}^4:H_1(x)=c_1,H_2(x)=c_2\right\},$$ be the invariant
surface (where $c=(c_1,c_2)$ is not a critical value).
Substituting $y_1=r\cos\theta$, $y_2=r\sin\theta$, in equations
\begin{eqnarray}
H_1&=&\frac{1}{2}\left(x_{1}^{2}+x_{2}^{2}\right)
+\frac{1}{4}\left(y_{1}^{2}+y_{2}^{2}\right)^{2}=c_1 \nonumber\\
H_2&=&x_{1}y_{2}-x_{2}y_{1}=c_2,\nonumber
\end{eqnarray}
we obtain
$$
\frac{1}{2}\left(\dot{r}^2+(r\dot{\theta})^2\right)+\frac{1}{4}
r^{2}=c_{1},\qquad r^{2}\dot{\theta}=-c_{2}.
$$
Hence
$$\left(r\dot{r}\right)^{2}+\frac{1}{2}r^4-2c_{1}r^{2}+c_{2}^{2}=0,$$
and $$w^{2}+P(z)=0,$$ where $w\equiv r\displaystyle{\dot{r}}$,
$z\equiv r^{2}$, and $P(z)=\frac{1}{2} z^3-2c_{1}z+c_{2}^{2}$. The
polynomial $P(z)$ is of degree $3$, the Riemann surface
$\mathcal{C}$
\begin{equation}\label{eqn:euler}
\mathcal{C}=\overline{\{(w,z):w^{2}+P(z)=0\}},
\end{equation}
is of genus $g=1$ (an elliptic curve). We thus have a single
holomorphic differential $$\omega=\frac{dz}{\sqrt{P(z)}},$$ and
the linearization occurs on the elliptic curve $\mathcal{C}$.
Although the variety $\mathcal{M}_c$ has dimension $2$, here we
have a reduction of dimension $1$ and, consequently, we get the
following result :

\begin{Theo}
The differential system (52) is completely linearized on the
Jacobian variety of $\mathcal{C}$, i.e. on the elliptic curve
$\mathcal{C}$(53).
\end{Theo}

\begin{Rem}
Further information and methods for resolving the Yang-Mills
system will be found for example in [35, 25] and references
therein.
\end{Rem}

\emph{\textbf{Conclusion}} : We do not consider here the solution
techniques based on the important notion of dynamical systems that
are algebraically completely integrable (the survey of these
results, as well as, the extensive list of references can be found
in [2, 27]). Let's just mention that the concept of algebraic
complete integrability is quite effective in small dimensions and
has the advantage to lead to global results, unlike the existing
criteria for real analytic integrability, which, at this stage are
perturbation results (in fact, the perturbation techniques
developed in that context are of a totally different nature).
However, besides the fact that many Hamiltonian dynamical
integrable systems possess this structure, another motivation for
its study is due to the fact that algebraic completely integrable
systems come up systematically whenever you study the isospectral
deformation of some linear operator containing a rational
indeterminate (indeed a theorem of Adler-Kostant-Symes [2] applied
to Kac-Moody algebras provides such systems which, by a theorem of
van Moerbeke-Mumford [2] are algebraic completely integrable). In
recent years, other important results have been obtained following
studies on the KP and KdV hierarchies (we refer the interested
reader for example to [24] for an exposition and a survey of the
results in this field, as well as, a list of references). In fact,
many problems related to algebraic geometry, combinatorics,
probabilities and quantum gauge theory,..., have been solved
explicitly by methods inspired by techniques from the study of
dynamical integrable systems.

\end{document}